\documentclass[]{revtex4-1}  

\draft 

\usepackage[]{graphicx}
\graphicspath{{figs/}}
\usepackage{float}
\usepackage{amsmath,amsfonts,amssymb,nicefrac,esint,bbm}
\DeclareMathAlphabet{\mathpzc}{OT1}{pzc}{m}{it}
\usepackage{xcolor}

\newcommand\nc{\newcommand}

\nc\wrt{w.r.t.\ }  \nc\cf{cf.\ }
\nc\lat\textit
\nc\ie{\lat{i.e.,\ }}  \nc\etal{\lat{et al.\ }}  \nc\eg{\lat{e.g.,\ }}
\nc\re[1]{(\ref{#1})}

\nc\FTminus{\hskip-0.4ex}
\nc\FTapprox{\mathord{\approx}\!\:}
\nc\FThcite{\,}

\nc\tensorr\mathbf
\nc\Tensor\boldsymbol
\nc\ident{\tensorr 1}
\nc\zero{\tensorr 0}
\nc\Sym[1]{#1^{\rm Sym}}
\nc\Nicefrac[2]{\nicefrac{#1}{#2}}

\nc\ee{{\rm e}}
\nc\ii{\mathbbm{i}}

\nc\hx{\hat{x}}
\nc\hht{\hat{t}}
\nc\cx{\check{x}}
\nc\ct{\check{t}}

\nc\dd{{\rm d}}
\nc\pd{\partial}
\nc\dt[1]{\frac{\dd #1}{\dd t}}
\nc\pdt[1]{\frac{\pd #1}{\pd t}}
\nc\ppdt[2]{\frac{\pd^#2 #1}{\pd t^#2}}
\nc\nab{\nabla}
\nc\pder[3]{\left. \frac{\pd #1}{\pd #2} \right|_{#3}}
\nc\Lapl{\Delta}
\nc\pdx[1]{\frac{\pd #1}{\pd x}}
\nc\dht[1]{\frac{\dd #1}{\dd \hht}}
\nc\pdht[1]{\frac{\pd #1}{\pd \hht}}
\nc\pdhx[1]{\frac{\pd #1}{\pd \hx}}
\nc\pdct[1]{\frac{\pd #1}{\pd \tau}}
\nc\ppdct[1]{\frac{\pd^2 #1}{\pd \tau^2}}

\nc\aT{a_T}
\nc\as{a_s}
\nc\cpp{c_p}
\nc\cv{c_v}
\nc\qqJ{\tensorr{J}}
\nc\qqP{\tensorr{P}}
\nc\qqk{\tensorr{k}}
\nc\qv{{\mathsf{v}}}
\nc\qqv{{\tensorr{v}}}
\nc\qqom{{\Tensor{\Omega}}}
\nc\qq{\dot{q}}
\nc\qqq{\dot{\tensorr{q}}}
\nc\Dv{\mathcal{D}_{\qqv}}
\nc\qqr{\tensorr{r}}

\nc\bp{\beta_p}
\nc\kT{\kappa_T}
\nc\qrho{\varrho}
\nc\qqPi{\Tensor{\Pi}}

\nc\brho{\qrho^0}
\nc\hrho{\hat{\qrho}}
\nc\crho{\qrho_{\rm c}}
\nc\bT{T^0}
\nc\hT{\hat{T}}
\nc\cT{T_{\rm c}}
\nc\hv{\hat{\qv}}
\nc\vpar{\qqv^*}
\nc\vper{\qqv^\circ}
\nc\bpp{p^0}
\nc\bbp{\bp^0}
\nc\bkT{\kT^0}
\nc\baT{\aT^0}
\nc\bcv{\cv^0}
\nc\bas{\as^0}
\nc\bcp{\cpp^0}
\nc\ba{\alpha^0}
\nc\bgam{\gamma^0}
\nc\bnu{\nu^0}
\nc\hPi{\hat{\Pi}}
\nc\hq{\hat{\qq}}

\nc\vph{\qv_{\rm ph}}

\nc\qKn{\mathpzc{Kn}}
\nc\qMa{\mathpzc{Ma}}
\nc\qPr{\mathpzc{Pr}}
\nc\qRe{\mathpzc{Re}}
\nc\qB{\mathpzc{B}}

\nc\qPea{\mathpzc{Pe}_{\rm a}}
\nc\qRea{\mathpzc{Re}_{\rm a}}
\nc\qEca{\mathpzc{Ec}_{\rm a}}

\nc\Dt{\Delta t}
\nc\hDt{\Delta \hat{t}}
\nc\hDx{\Delta \hat{x}}
\nc\half{\nicefrac{1}{2}}
\nc\quart{\nicefrac{1}{4}}
\nc\tquart{\nicefrac{3}{4}}

\nc\reta{\mathcal{R}_{\eta}}
\nc\bReta{\mathcal{R}_{\eta}^0}
\nc\bbeta{\eta^0}
\nc\blam{\lambda^0}

\nc\X{\mathsf{X}}
\nc\Xrev{\X_{\mathrm{rev}}}
\nc\Xirr{\X_{\mathrm{irr}}}
\nc\phirev{\Phi_{\mathrm{rev}}}
\nc\psiirrev{\Phi_{\mathrm{irr}}}

\nc\pDhx[1]{\frac{#1}{\hDx}}

\nc{\leftf}{\mathopen{}\mathclose\bgroup\left}  
\nc{\rightf}{\aftergroup\egroup\right}   

\nc\Exp{\operatorname{Exp}}
\nc\Expof[1]{\Exp\leftf(#1\rightf)}
\nc\tX{\tilde{\X}}
\nc\tXrevthr{\tX_{\mathrm{rev},3}}
\nc\tXirrthr{\tX_{\mathrm{irr},3}}
\nc\tXone{\tX_{1}}
\nc\tXtwo{\tX_{2}}
\nc\Ord[1]{\mathcal{O}\leftf(#1\rightf)}
\nc\Qper{/}
\nc\nuR{\bnu \left( \bReta + \frac{4}{3} \right)}
\nc\Qds{\displaystyle}

\nc\VV{\vphantom{\frac{\hT_{n+\half}\hq_{n}}{\hT_{n+\half}}}}
\nc\VVV{\VV \\ \VV}

\begin{document}


\setlength{\tabcolsep}{0.45em}  

\title{The piston effect in supercritical fluids investigated via a reversible–irreversible vector field splitting-based explicit time integration scheme} 



\author{Donát M. Takács}
\affiliation{Department of Energy Engineering, Faculty of Mechanical Engineering, Budapest University of Technology and Economics, Műegyetem rkp.\ 3., H-1111 Budapest, Hungary}
\affiliation{Montavid Thermodynamic Research Group, Society for the Unity of Science and Technology, Lovas út 18., H-1012 Budapest, Hungary}

\author{Tamás Fülöp}
\affiliation{Department of Energy Engineering, Faculty of Mechanical Engineering, Budapest University of Technology and Economics, Műegyetem rkp.\ 3., H-1111 Budapest, Hungary}
\affiliation{Montavid Thermodynamic Research Group, Society for the Unity of Science and Technology, Lovas út 18., H-1012 Budapest, Hungary}

\author{Róbert Kovács}
\affiliation{Department of Energy Engineering, Faculty of Mechanical Engineering, Budapest University of Technology and Economics, Műegyetem rkp.\ 3., H-1111 Budapest, Hungary}
\affiliation{Department of Theoretical Physics, Institute for Particle and Nuclear Physics, HUN-REN Wigner Research Centre for Physics, Konkoly-Thege Miklós út 29--33., H-1121 Budapest, Hungary}
\affiliation{Montavid Thermodynamic Research Group, Society for the Unity of Science and Technology, Lovas út 18., H-1012 Budapest, Hungary}

\author{Mátyás Szücs}
\email[]{szucs.matyas@gpk.bme.hu}
\affiliation{Department of Energy Engineering, Faculty of Mechanical Engineering, Budapest University of Technology and Economics, Műegyetem rkp.\ 3., H-1111 Budapest, Hungary}
\affiliation{Department of Theoretical Physics, Institute for Particle and Nuclear Physics, HUN-REN Wigner Research Centre for Physics, Konkoly-Thege Miklós út 29--33., H-1121 Budapest, Hungary}
\affiliation{Montavid Thermodynamic Research Group, Society for the Unity of Science and Technology, Lovas út 18., H-1012 Budapest, Hungary}


\date{\today}

\begin{abstract}
    In the vicinity of the liquid--vapor critical point, supercritical fluids behave strongly compressibly and, in parallel, thermophysical properties have strong state dependence. These lead to various peculiar phenomena, one of which being the piston effect where a sudden heating induces a mechanical pulse. The coupling between thermal and mechanical processes, in the linear approximation, yields a non-trivially rich thermoacoustics. The numerous applications of supercritical fluids raise the need for reliable yet fast and efficient numerical solution for thermoacoustic time and space dependence in this sensitive domain. Here, we present a second-order accurate, fully explicit staggered space-time grid finite difference method for such coupled linear thermoacoustic problems. Time integration is based on the splitting of the state space vector field representing the interactions that affect the dynamics into reversible and irreversible parts, which splitting procedure leads to decoupled wave and heat equations. The former is a hyperbolic partial differential equation, while the latter is a parabolic one, therefore, different time integration algorithms must be amalgamated to obtain a reliable, dispersion error-free, and dissipation error-free numerical solution. Finally, the thermoacoustic approximation of the supercritical piston effect is investigated via the developed method.
\end{abstract}

\pacs{05.70.Ln, 44.05.+e,  47.11.Bc}

\maketitle 

\section{Introduction}

\begin{figure*}[!htb]
    \begin{minipage}[c]{0.32\textwidth}
        \includegraphics[height=3.7cm]{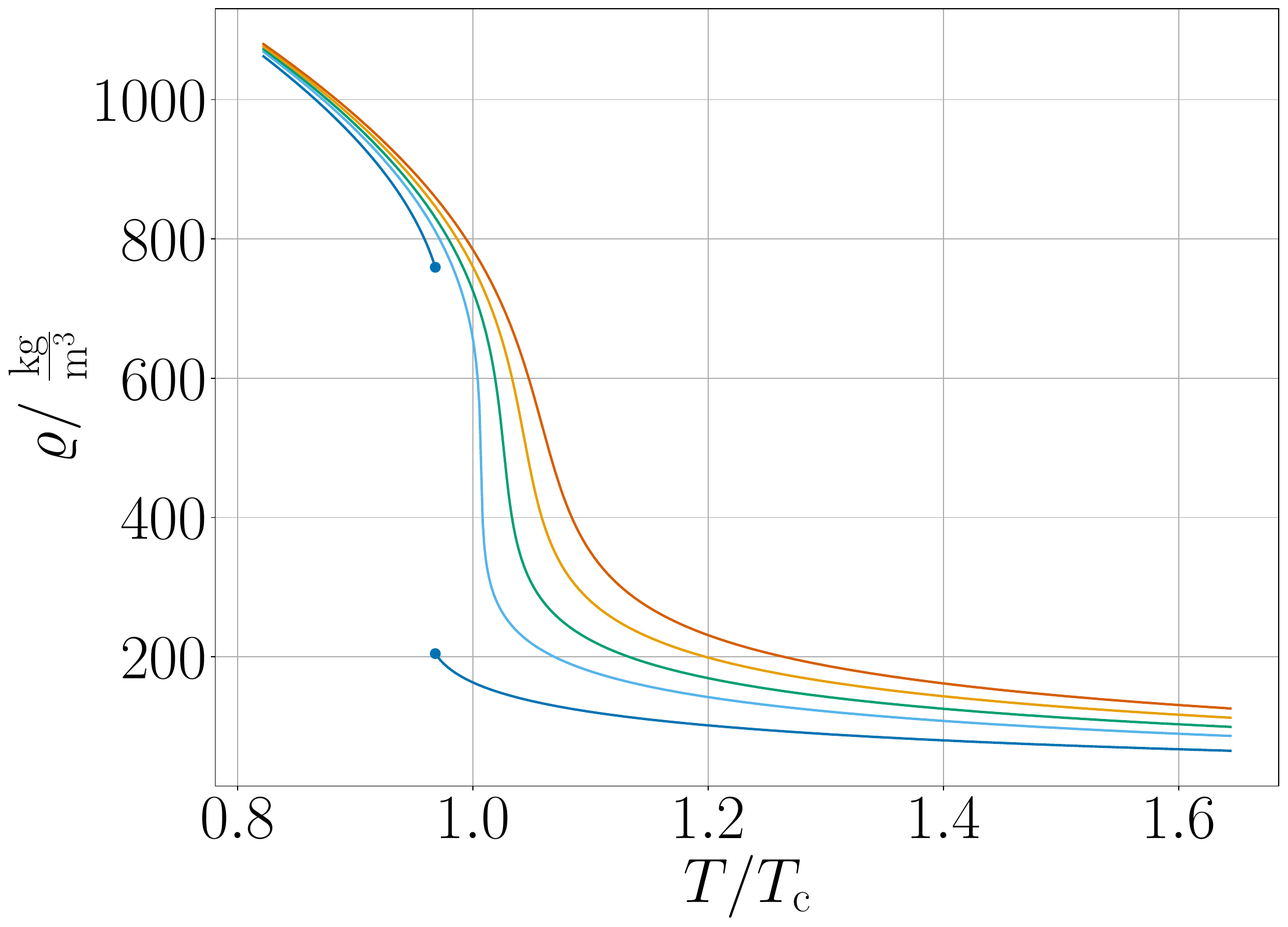}
    \end{minipage}%
    \hspace{.007\textwidth}%
    \begin{minipage}[c]{0.32\textwidth}
        \includegraphics[height=3.7cm]{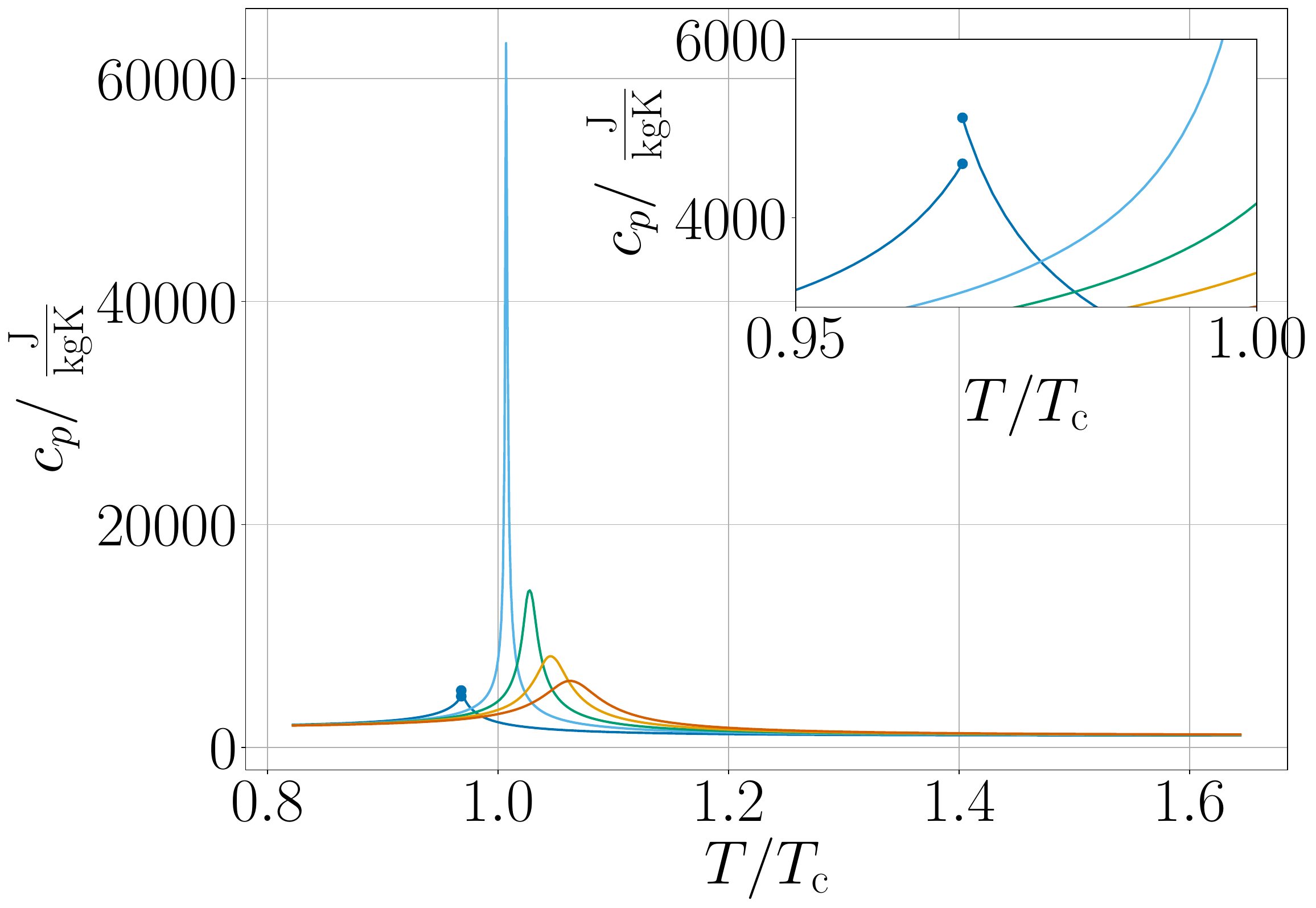}
    \end{minipage}%
    \hfill  
    \begin{minipage}[c]{0.32\textwidth}
        \includegraphics[height=3.7cm]{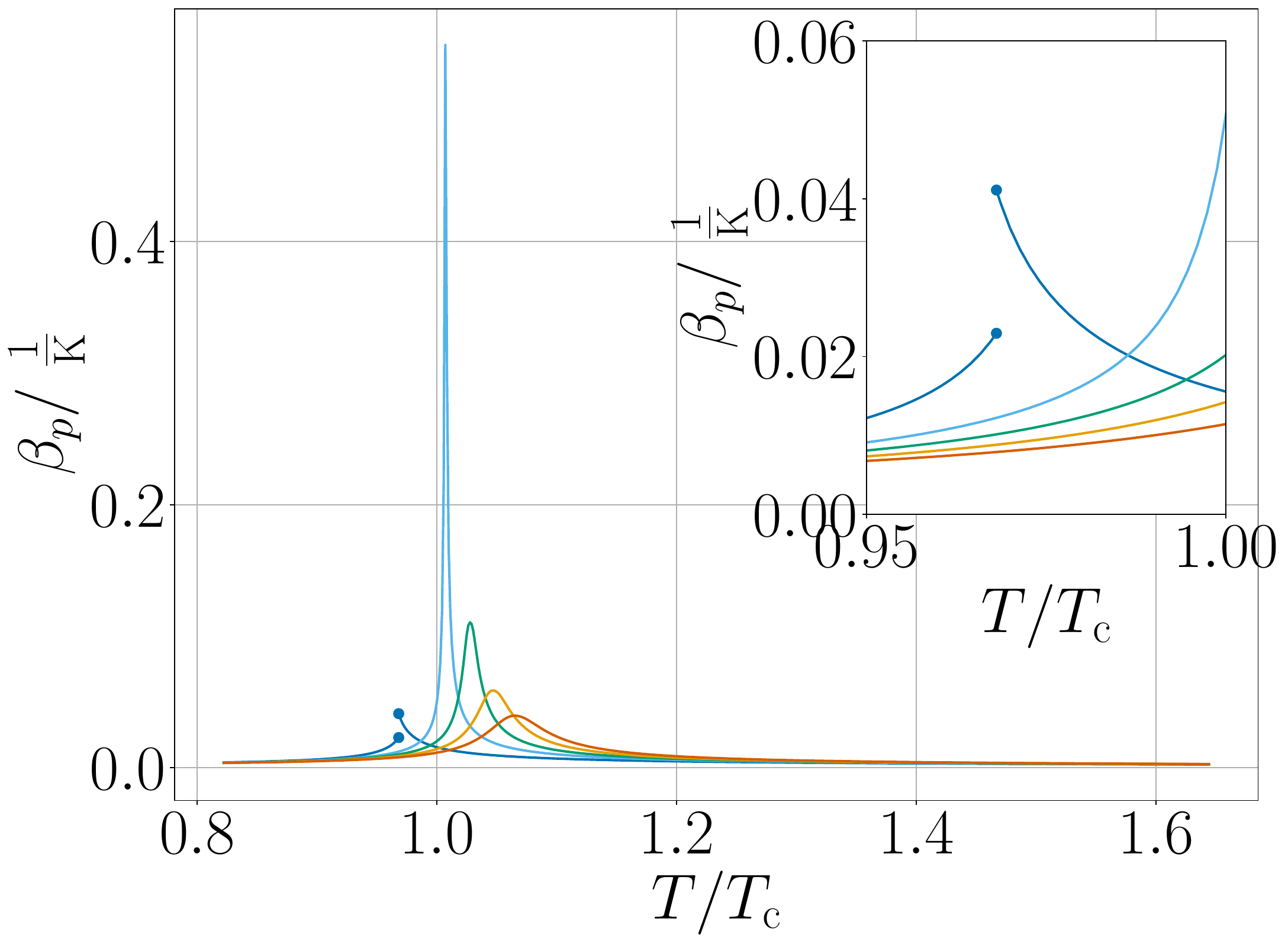}
    \end{minipage}%
    \newline\null
    \hspace{.011\textwidth}%
    \begin{minipage}[c]{0.32\textwidth}
        \includegraphics[height=3.7cm]{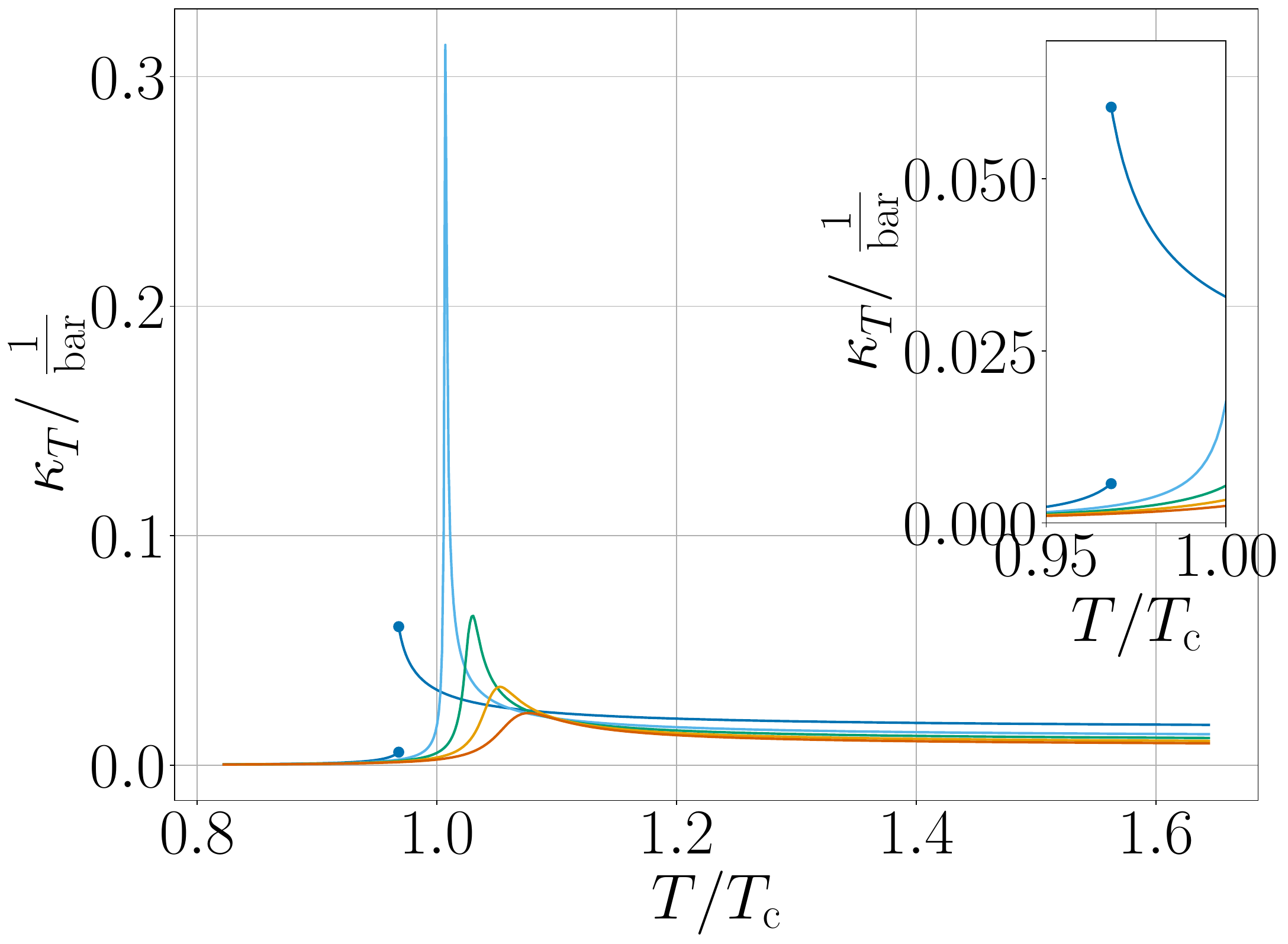}
    \end{minipage}%
    \hspace{.016\textwidth}%
    \begin{minipage}[c]{0.32\textwidth}
        \includegraphics[height=3.7cm]{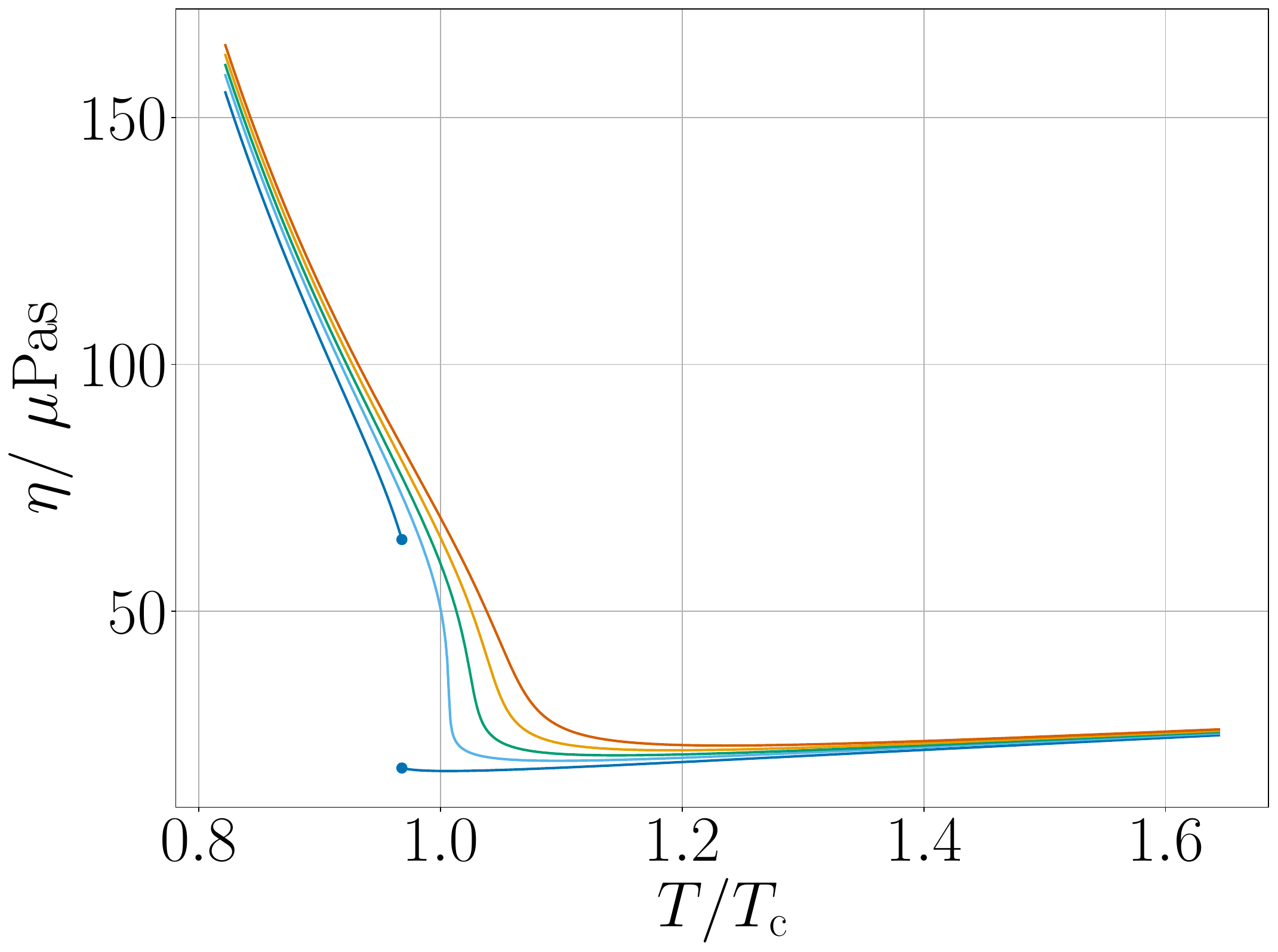}
    \end{minipage}%
    \begin{minipage}[c]{0.32\textwidth}
        \includegraphics[height=3.7cm]{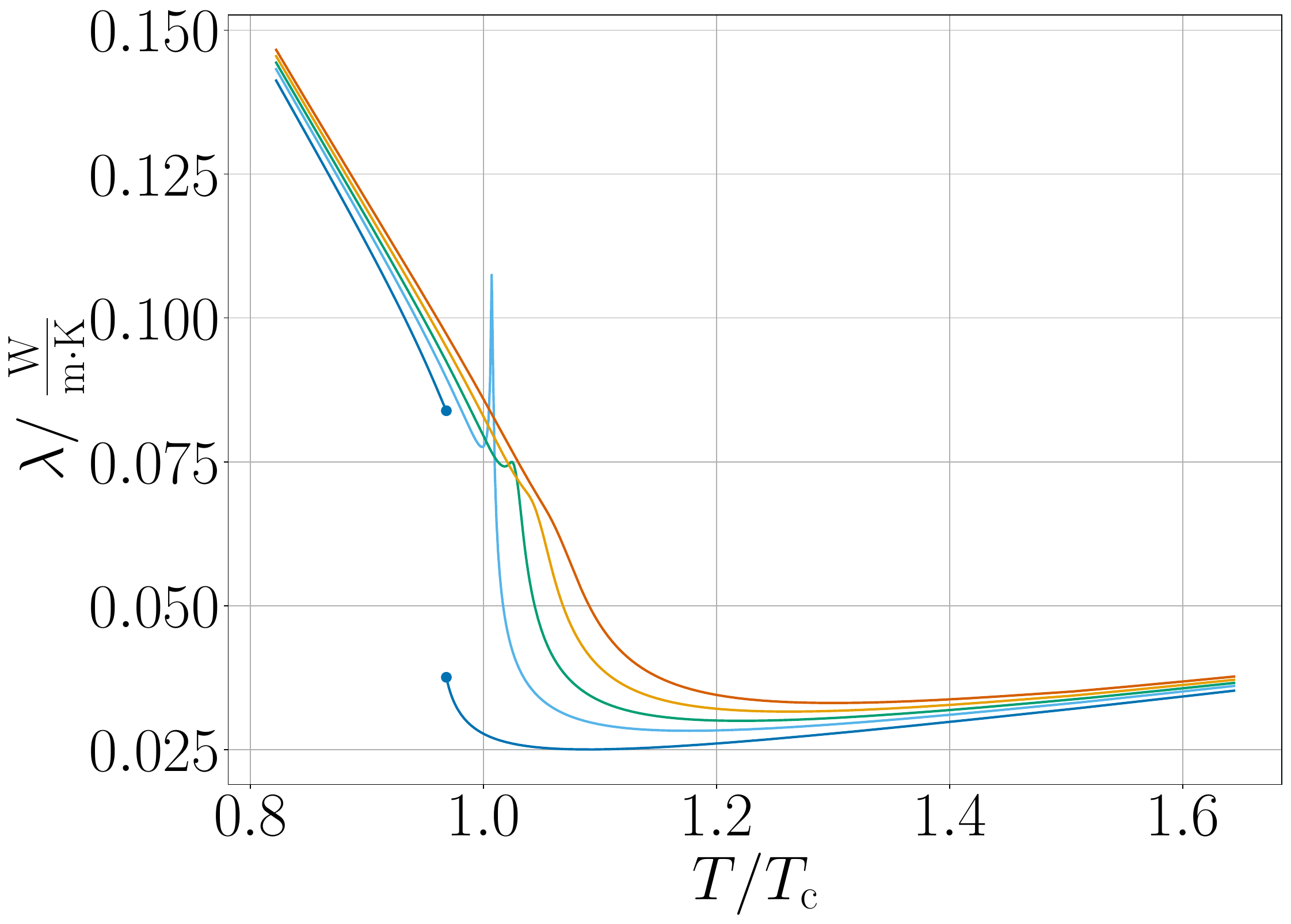}
    \end{minipage}%
    \caption{\emph{From left to right:} Temperature dependence of density, isobaric specific heat capacity, isobaric volume thermal expansion coefficient, isothermal compressibility, shear viscosity, and thermal conductivity for carbon dioxide along the isobaric curves $ p/p_{\rm c} = 0.8 , \ 1.05 , \ 1.2 , \ 1.35 , \ 1.5 $ (from dark blue to red, respectively). Critical pressure and temperature of carbon dioxide are $ p_{\rm c} = 7.377 \ {\rm MPa} $ and $ T_{\rm c} = 304.128 \ {\rm K} $. Below the critical pressure (\ie in the figures $ p/p_{\rm c} = 0.8 $) jump occurs in the thermophysical properties, where left and right sides belong to the liquid and gas phase, respectively, this jump disappears for pressures above the critical pressure. Numerical data are taken from the NIST database \cite{lemmon2024thermophysical}.}
    \label{fig:SCCO2}
\end{figure*}

A substance is called a supercritical fluid if it has a pressure and temperature above its liquid--vapor critical point. In the supercritical region, fluids are strongly compressible, liquid and gas characteristics are indistinguishable, and thermophysical properties such as specific heat capacities, thermal expansion coefficient, compressibility, thermal conductivity, or viscosity, undergo rapid changes \cite{carles2010brief,imre2019anomalous}. This behavior is even more visible in a small vicinity of the critical point, where the former quantities can change by several orders of magnitude due to the effect of a temperature difference at the order of one Kelvin. The above statements are demonstrated in Figure \ref{fig:SCCO2} for carbon dioxide. Note that we provide values of state variables and material parameters at three decimal places as supercritical fluids are notably sensitive to minor changes. Since supercritical fluids are strongly compressible, bulk viscosity also plays an important role in thermomechanical processes, about which information is only partially available, \eg in \cite{carles1998effect,onuki1997dynamic,hasan2012thermoacoustic}.

Nowadays, the technical potential inherent in supercritical fluids is not only exploited by the chemical industry (such as in extraction, drying, cleaning, sterilization or polymer processing) \cite{knez2014industrial}, but it also plays an important role in more and more applications in energy engineering. Supercritical carbon dioxide is used as a working fluid for enhanced geothermal systems \cite{reinsch2017utilizing,dobson2017supercritical}. Additionally, a possible sustainable, economical, and eco-friendly solution for power supply can be provided by a concept of generation IV nuclear systems, which use supercritical water (at a state of about 25~MPa and 500~$^\circ$C) both as coolant and as a working fluid \cite{rahman2020design,wu2022review}, called Supercritical Water Reactors. Although supercritical power cycles have several advantages, \eg the absence of boiling -- hence technical difficulties like `boiling crisis' \cite{theofanous2002boiling1,theofanous2002boiling2} are avoidable --, the presence of supercritical fluid presents other issues such as heat transfer deterioration \cite{longmire2022onset}.

In supercritical fluids, the relatively high value of the thermal expansion coefficient \footnote{which is correlated to high compressibility, via the Imre ellipse \cite{takacs2024leading}} leads to the strong coupling of thermal and flow processes, therefore, the coupling of thermal and mechanical processes is inevitable. Due to the intense thermal expansion, pressure waves can be induced in the fluid not just by mechanical changes, but also by thermal effects. This attribute is hard to observe and explain without a proper coupled thermoacoustic model \cite{carles2006thermoacoustic}.

\begin{figure}[!htb]
    \includegraphics[height=3.7cm]{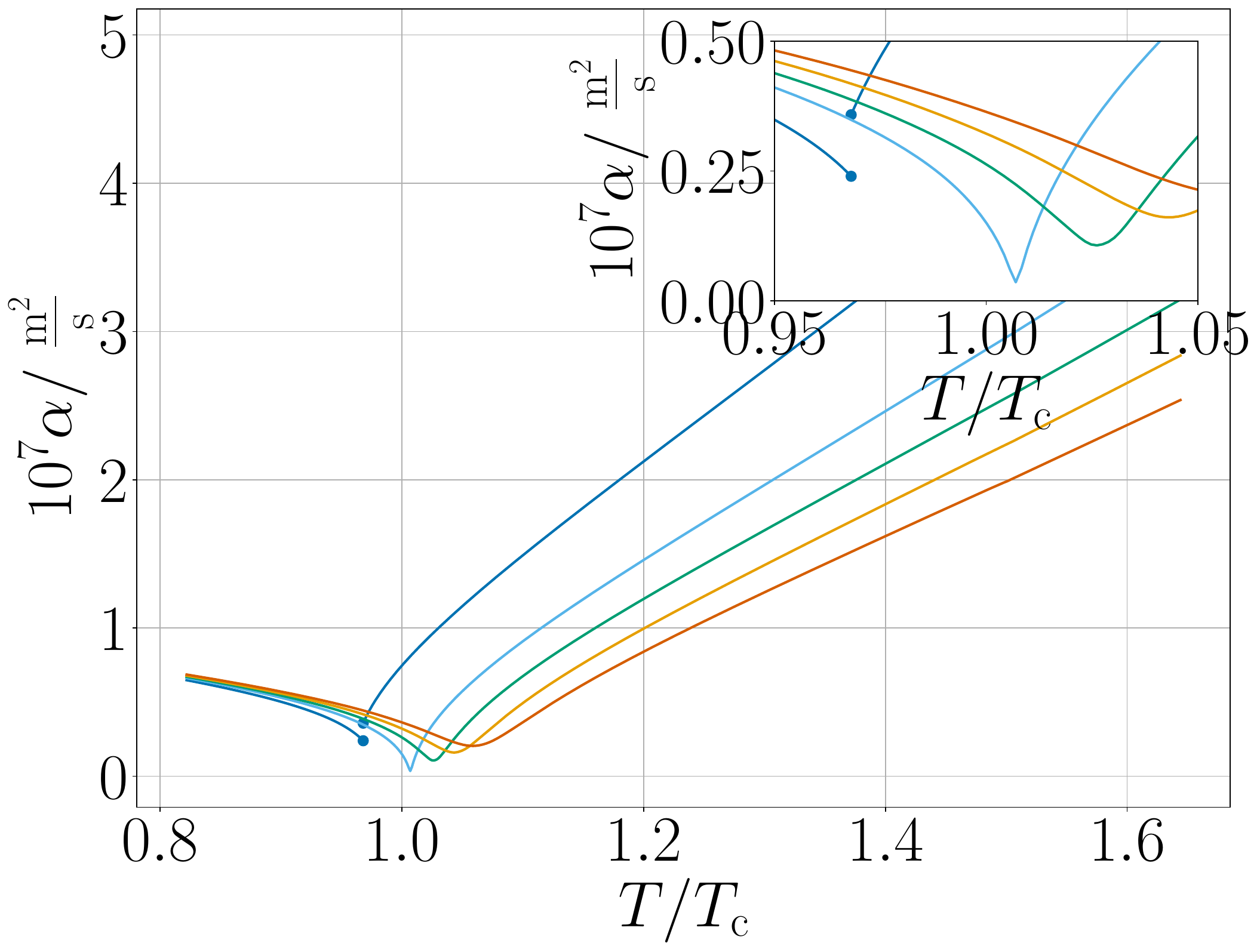} \centering
    \caption{Temperature dependence of thermal diffusivity for carbon dioxide along the isobaric curves $ p/p_{\rm c} = 0.8 , \ 1.05 , \ 1.2 , \ 1.35 , \ 1.5 $ (from dark blue to red, respectively). Numerical data are calculated from the NIST database \cite{lemmon2024thermophysical}.}
    \label{fig:SCCO2-as-alp}
\end{figure}

The coupling of thermal and mechanical processes is well observed in the vicinity of the liquid--vapor critical point and is sometimes referred to as the `fourth mechanism of heat transport' \cite{zappoli2003nearcritical,zappoli2015heat}. Regarding heat conduction, thermal equilibration near the critical point is expected to be very slow since thermal diffusivity tends to zero in the vicinity of the critical point (see Figure \ref{fig:SCCO2-as-alp}). This effect is usually called the `critical slowing down'. In contrast, a fast temperature equilibration is observed in experiments, which under normal gravity circumstances is explained by the mixing of the fluid induced by buoyancy convection. However, measurements performed under microgravity circumstances presented the same fast equilibration of temperature \cite{straub1995dynamic}, but density homogenization was still very slow. The fast temperature equilibration is called `critical speeding up' phenomenon, which is explained by Onuki \etal with an isentropic -- wave-like -- propagation of temperature \cite{onuki1990fast}.

In detail, due to surface heating near the heated wall, a thermal boundary layer occurs in the fluid, which suddenly expands due to the relatively large value of the isobaric thermal expansion coefficient and compresses the rest of the fluid just like a piston. That is why this phenomenon is usually called the piston effect \cite{zappoli2003nearcritical,zappoli1990anomalous}. In parallel to the sudden expansion, a decrease of density and an increase of pressure is implied, which disturbances propagate through the fluid with the speed of sound, and therefore, at the same time, the propagation of temperature at the speed of sound can be observed.

The previous statements enlighten that there exist two typical time scales in the piston effect, an \textit{acoustic} time scale corresponding to the isentropic speed of sound and a \textit{thermal} time scale corresponding to thermal diffusivity \cite{garrabos1998relaxation}. In contrast to liquids and gases, the transport of thermal energy in supercritical fluids simultaneously occurs at these two time scales initially. Due to the compression caused by thermal expansion, a sudden increase in temperature is observed, while the homogenization caused by diffusion in the case of an adiabatic isolated sample requires a longer time.

The piston effect has been investigated in many works, \eg Boukari \etal neglected acoustic propagation and solved Fourier's heat conduction equation together with pressure changes \cite{boukari1990critical}, Zappoli \etal presented the solution of the full Navier--Stokes---Fourier equations with constant thermophysical parameters \cite{zappoli1990anomalous}, while Wagner \etal considered state dependence of thermophysical properties, too \cite{wagner2001variable}. 

Customary space-time discretizations and numerical solutions of the linearized approximation of the compressible Navier--Stokes---Fourier equations typically lead to significant numerical artifacts, including dispersion error (artificial oscillations) and dissipation error (artificial decrease or increase in the amplitude).
In parallel, many of the common numerical methods used for time integration in commercial software are known to introduce significant numerical artifacts when applied to to wave propagation phenomena, particularly when coupled with dissipative processes\cite{fulop2020thermodynamical,pozsar2020four}: resulting in dispersion and dissipation errors, especially over longer time scales. These studies have highlighted that reliable -- (nearly) numerical artifact-free -- solutions for dissipative and dispersive wave propagation phenomena require special algorithms.

To achieve a reliable and accurate time integration algorithm, we took a physics-inspired approach and developed a fully explicit, second-order accurate time integration scheme based on a reversible--irreversible splitting procedure motivated by the GENERIC (General Equation for the Non-Equilibrium Reversible--Irreversible Coupling) framework \cite{grmela1997dynamics}.

GENERIC considers the time evolution of state variables via a sum of two vector fields: one conserves entropy, \ie connects to reversible time evolution, while the other increases entropy, \ie describes the dissipative effects. A direct numerical exploitation of GENERIC is realized by Shang and Öttinger \cite{shang2020structurepreserving}, however, they only consider finite degree-of-freedom systems, which are described by ordinary differential equations. The GENERIC integrators introduced by them are based on a special operator splitting procedure: namely, a given time step is split along the entropy-conserving and entropy-increasing contributions. 

Here, we follow a similar idea for the time integration after the spatial semi-discretization of the governing equations. We split the vector field\footnote{In order to avoid creating the impression that this splitting is applicable to linear equations only, instead of the frequently used term `operator splitting' we refer to this approach as \emph{vector field splitting}.} representing the interactions into two parts according to their contribution to entropy production. As we will later show, in linearized thermoacoustics, the entropy-conserving contribution leads to a wave equation, while the entropy-increasing contribution leads to the heat equation. In a previous work \cite{fulop2020thermodynamical}, we have presented a second-order accurate, explicit quasi-symplectic numerical scheme for the wave equation, which -- with the appropriate choice of the Courant number -- results in a dispersion- and dissipation error-free numerical solution. Unfortunately, this scheme can not be applied directly to solving the heat equation. However, here, along the lines of the reversible--irreversible splitting procedure, the quasi-symplectic scheme is applied to solve the reversible contribution, while a different scheme (here, the explicit midpoint method) is applied to solve the irreversible contribution. Although splitting-based numerical methods are not unknown in fluid mechanics -- see, e.g.\, the Alternating-Direction Implicit (ADI) method or the Pressure-Implicit with Splitting of Operators (PISO) algorithm \cite{anderson1995computational,versteeg2007introduction} -- our explicit method provides a novel strategy for numerical techniques to treat dispersive and dissipative wave  ropagation phenomena while producing fast-obtainable robust and reliable results.

The outline of the present article is as follows. In Sec.~\ref{sec:th-ac}, we briefly review the governing equations of thermofluid dynamics and its linear approximation, and the corresponding entropy production is also investigated. Sec.~\ref{sec:num-scheme} presents the reversible--irreversible vector field splitting-based numerical scheme for one spatial dimensional linear thermoacoustic problems. In Sec.~\ref{sec:piston}, numerical experiments are performed to analyze the schemes, and then the piston effect is investigated with the developed method in carbon dioxide, both in a nearly-ideal gas state and in a supercritical fluid state. Finally, in Sec.~\ref{sec:concl}, conclusions are drawn, and possibilities for the further development of the scheme are outlined.

\section{Linear thermoacoustics}\label{sec:th-ac}

Thermomechanical processes of single-phase pure fluids are described by the balance equations of mass, linear momentum, and internal energy formulated w.r.t. an inertial reference frame as 
\begin{align}
    \label{eq:bal-m}
    \Dv \qrho + \qrho \nab \cdot \qqv &= 0 , \\
    \label{eq:bal-v}
    \qrho \Dv \qqv &= - \qqP \cdot \nab ,
    \vphantom{\qqq \Sym{\left( \qqv \otimes \nab \right)}}  
    \\
    \label{eq:bal-u}
    \qrho \Dv u &= - \nab \cdot \qqq - \qqP : \Sym{\left( \qqv \otimes \nab \right)} 
\end{align}
with (mass) density $ \qrho $, velocity $ \qqv $, (mass) specific internal energy $ u $, symmetric pressure tensor $ \qqP $ and heat current density $ \qqq $, which all are fields, \ie functions of time $ t $ and spatial coordinates $ \qqr $, furthermore, $ \Sym{} $ denotes the symmetric part of a second-order tensor, \ie if $ \tensorr{A} $ denotes an arbitrary second-order tensor then its symmetric part is $ \Sym{ \tensorr{A} } = \frac{1}{2} \left( \tensorr{A} + \tensorr{A}^{\rm T} \right) $ with $ ^{\rm T} $ denoting the transpose of the tensor \cite{anderson2001fundamentals}. In equations \re{eq:bal-m}--\re{eq:bal-u}, we have neglected volumetric source terms. The material time derivative
\begin{align}
    \label{eq:Dv}
    \Dv \bullet = \pdt \bullet + \left( \bullet \otimes \nab \right) \cdot \qqv
\end{align}
expresses the temporal change of any physical property at a fixed material point, which travels with velocity $ \qqv $, here and later on $ \pdt{} $ denotes the partial time derivative and $ \nab $ is the nabla operator representing gradient, divergence or curl depending on the tensorial multiplication%
\footnote{Let us make a short remark on the tensorial notations. If $ \Tensor{\varphi} $ denotes a quantity with tensorial order $ N \ge 0 $ then its components \wrt the corresponding Cartesian basis vectors $ \tensorr{e}_j , \ j = 1,2,3 $ are denoted as $ \varphi_{i_{1} \dots i_{N}} , \ i_{1} , \dots , i_{N} = 1,2,3 $.  In accordance with the usual notation and usage of the nabla operator in continuum mechanics, \textit{right} gradient, \textit{right} divergence, and \textit{right} curl (indicated by subscript R) are defined as
\begin{align*}
	\operatorname{grad}_{\rm R} \Tensor{\varphi} &:= \Tensor{\varphi} \otimes \nabla = \pd_j \varphi_{i_{1} \dots i_{N}} \tensorr{e}_{i_1} \otimes \dots \otimes \tensorr{e}_{i_N} \otimes \tensorr{e}_{j}, \\
	\operatorname{div}_{\rm R} \Tensor{\varphi} &:= \Tensor{\varphi} \cdot \nabla = \pd_j \varphi_{i_{1} \dots i_{N - 1} j} \tensorr{e}_{i_1} \otimes \dots \otimes \tensorr{e}_{i_{N-1}} , \\
	\operatorname{curl}_{\rm R} \Tensor{\varphi} &:= \Tensor{\varphi} \times \nabla  = \epsilon_{i_{N} j k} \pd_j \varphi_{i_{1} \dots i_{N}} \tensorr{e}_{i_1} \otimes \dots \otimes \tensorr{e}_{i_{N -1}} \otimes \tensorr{e}_{k} ,
\end{align*}
where Einstein summation is applied and
\begin{align*}
    \epsilon_{ijk} =
    \begin{cases}
        1 , & \text{if } \left( i , j , k \right) \text{ is an even permutation} , \\
        - 1 , & \text{if } \left( i , j , k \right) \text{ is an odd permutation} , \\
        0 , & \text{otherwise}
    \end{cases}
\end{align*}
is the Levi-Civita symbol. Especially, for an arbitrary scalar field $ a $
\begin{align*}
	a \otimes \nabla = \nabla \otimes a =: \nabla a ;
\end{align*}
in this case, divergence and curl are not interpreted. For an arbitrary vector field $ \tensorr{v} $
\begin{align*}
	\tensorr{v} \otimes \nabla &= \left( \nabla \otimes \tensorr{v} \right)^{\rm T} , &
	\tensorr{v} \cdot \nabla &= \nabla \cdot \tensorr{v} , &
	\tensorr{v} \times \nabla &= - \nabla \times \tensorr{v} .
\end{align*}
In what follows, wherever needed, parentheses indicate whether $\nabla$ acts to the left or to the right [$(a \nabla) b$ vs. $a (\nabla b)$].
}.

The realized process has to fulfill the second law of thermodynamics, which is expressed by the balance of entropy
\begin{align}
    \label{eq:bal-s}
    \qrho \Dv s &= - \nab \cdot \qqJ + \Sigma , 
\end{align}
where $ s $, $ \qqJ $, and $ \Sigma $ are the specific entropy, the entropy current density, and the positive semidefinite entropy production rate density, respectively \cite{degroot1962nonequilibrium,gyarmati1970nonequilibrium}. Therefore, the balance equations \re{eq:bal-m}--\re{eq:bal-u} are equipped with the consistency condition \re{eq:bal-s}. Specific entropy can be obtained as the potential function of the (co)vector field $ \begin{pmatrix} \frac{1}{T} & \frac{p}{T} \end{pmatrix} $ interpreted on the thermodynamic state space spanned by $ \left( u , v \right) $, where $ v = \frac{1}{\qrho} $ is the specific volume and, $ T $ and $ p $ denote the (absolute) temperature and hydrostatic pressure, respectively. This property is expressed through the Gibbs relation \cite{bejan2016advanced}
\begin{align}
    \label{eq:Gibbs}
    \dd s &= \frac{1}{T} \left( u , v \right) \dd u + \frac{p}{T} \left( u , v \right) \dd v = \frac{1}{T} \left( u , \qrho \right) \dd u - \frac{1}{\qrho^2} \frac{p}{T} \left( u , \qrho \right) \dd \qrho ,
\end{align}
which expresses the functional relationships among density, temperature, pressure, and specific internal energy as
\begin{align}
    \label{eq:eos-Gibbs}
    &&
    \pder{s}{u}{\qrho} &= \frac{1}{T} \left( u , \qrho \right) , &
    \pder{s}{\qrho}{u} &= - \frac{1}{\qrho^2} \frac{p}{T} \left( u , \qrho \right) .
    &&
\end{align}
Assuming invertibility, these relationships are usually given by the thermal and caloric equations of state
\begin{align}
    \label{eq:eos}
    &&
    p &= p ( T , \qrho ) , & 
    u &= u ( T , \qrho ) .
    &&
\end{align}

Expressing the l.h.s. of the entropy balance \re{eq:bal-s} with $ u $ and $ \qrho $ through the Gibbs relation, applying \re{eq:bal-m} and \re{eq:bal-u} as constraints and separating a full divergence term one obtains
\begin{align}
    \label{eq:ds/dt}
    \begin{split}
    \qrho \Dv s &= \frac{1}{T} \qrho \Dv u - \frac{1}{\qrho} \frac{p}{T} \Dv \qrho = - \frac{1}{T}
    ( \nab \cdot \qqq )
    - \frac{1}{T} \qqP : \Sym{\left( \qqv \otimes \nab \right)} + \frac{p}{T}
    ( \nab \cdot \qqv )
    \\
    &= - \nab \cdot \left( \frac{1}{T} \qqq \right) - \frac{\qqq}{T^2}
    ( \nab T ) - \frac{1}{T} \left( \qqP - p \ident \right) : \Sym{\left( \qqv \otimes \nab \right)} ,
    \end{split}
\end{align}
where $ \ident $ denotes the identity tensor. Comparing \re{eq:bal-s} and \re{eq:ds/dt}, entropy current density and entropy production rate density are identified as
\begin{align}
    \qqJ &= \frac{1}{T} \qqq , \\
    \label{eq:ent-pr}
    \Sigma &= - \frac{\qqq}{T^2}
    ( \nab T )
    - \frac{1}{T} \left( \qqP - p \ident \right) : \Sym{\left( \qqv \otimes \nab \right)} \ge 0 .
\end{align}
Introducing the viscous pressure tensor $ \qqPi := \qqP - p \ident $ and assuming isotropic fluid (therefore, the different tensorial orders and characters do not couple), the positive semidefiniteness of \re{eq:ent-pr} is ensured via the linear equations
\begin{align}
    \label{eq:Fou}
    \qqq &= - \lambda 
    ( \nab T )
    , \\
    \label{eq:New}
    \qqPi &= - \left( \eta_{\rm Vol} - \frac{2}{3} \eta \right) \left( \nabla \cdot \qqv \right) \ident - 2 \eta \Sym{\left( \qqv \otimes \nab \right)} ,
\end{align}
where $ \lambda \ge 0 $, $ \eta_{\rm Vol} \ge 0 $ and $ \eta \ge 0 $ are the thermal conductivity, the volume, and shear viscosities, respectively, which in general are state-dependent thermophysical parameters. Note that \re{eq:Fou} is Fourier's law of heat conduction and \re{eq:New} is Newton's law of viscosity generalized incorporating the volume viscosity.

To summarize, the balance equations \re{eq:bal-m}--\re{eq:bal-u}, the thermal and caloric equations of state \re{eq:eos}, the constitutive equations on the heat current density and viscous pressure tensor \re{eq:Fou} and \re{eq:New} together with appropriate initial and boundary conditions form a closed system on the variables $ T $, $ \qrho $, $ p $, $ u $, $ \qqv $, $ \qqq $ and $ \qqPi $.

The thermal and caloric equations of state establish relationships among density, temperature, pressure, and specific internal energy, therefore, specific internal energy and pressure can be eliminated from \re{eq:bal-v} and \re{eq:bal-u}, \ie
\begin{align}
    \label{eq:NS}
    \qrho \Dv \qqv &= - \pder{p}{T}{\qrho}
    ( \nab T )
    - \pder{p}{\qrho}{T}
    ( \nab \qrho )
    - \qqPi \cdot \nab , \\
    \label{eq:heat-eq-Trho}
    \begin{split}
        \qrho \pder{u}{T}{\qrho} \Dv T &= - \nab \cdot \qqq - \qrho \pder{u}{\qrho}{T} \Dv \qrho \\
        & \hskip 2.5ex - p \left( \nab \cdot \qqv \right) - \qqPi : \Sym{\left( \qqv \otimes \nab \right)} .
    \end{split}
\end{align}
The appearing partial derivatives of the thermal and caloric equations of state are in strong connections with material properties, namely, isochoric specific heat capacity, isothermal compressibility, and isobaric volumetric thermal expansion coefficient are defined as
\begin{align}
    \label{eq:cv}
    \cv &:= \pder{u}{T}{\qrho}
    >
    0 , \\
    \label{eq:kT}
    \kT &:= \frac{1}{\qrho} \frac{1}{\pder{p}{\qrho}{T}}
    >
    0 , \\
    \label{eq:bp}
    \bp &:= \frac{1}{\qrho} \frac{\pder{p}{T}{\qrho}}{\pder{p}{\qrho}{T}} ,
\end{align}
where positivity of $ \cv $ and $ \kT $ ensure material stability, which actually realize the Le Ch\^{a}telier--Braun principle \cite{grigull1964prinzip,matolcsi2004ordinary}. If the above three material properties \eqref{eq:cv}--\eqref{eq:bp} are known then all others can be determined via these; for instance, isobaric specific heat capacity and specific heat ratio can be expressed as
\begin{align}
    \label{eq:cp}
    \cpp &= \pder{u}{T}{\qrho} + \frac{T}{\qrho^2} \frac{\left( \pder{p}{T}{\qrho} \right)^2}{\pder{p}{\qrho}{T}} = \cv + \frac{T}{\qrho} \frac{\bp^2}{\kT} \ge \cv , \\
    \label{eq:gam}
    \gamma &= \frac{\cpp}{\cv} = 1 + \frac{T}{\qrho^2} \frac{\left( \pder{p}{T}{\qrho} \right)^2}{\pder{u}{T}{\qrho} \pder{p}{\qrho}{T}} = 1 + \frac{T}{\qrho} \frac{\bp^2}{\cv \kT} \ge 1 ;
\end{align}
in both cases, equality holds if thermal expansion is neglected, \ie $ \bp = 0 $. Thanks to the entropic property \eqref{eq:Gibbs}, the relationship
\begin{align}
    \label{eq:urT-pTr}
    \pder{u}{\qrho}{T} = - \frac{1}{\qrho^2} \left( T \pder{p}{T}{\qrho} - p \right)
\end{align}
holds, therefore, \re{eq:heat-eq-Trho} can be further simplified to
\begin{align}
    \label{eq:heat-eq}
    \qrho \pder{u}{T}{\qrho} \Dv T &= - \nab \cdot \qqq - T \pder{p}{T}{\qrho} ( \nab \cdot \qqv ) - \qqPi : \Sym{\left( \qqv \otimes \nab \right)} .
\end{align}

{Limiting ourselves to small perturbations around a homogeneous static equilibrium state characterized by $ \bT = {\rm const.} $ and $ \brho = {\rm const.} $, and, assuming negligible second-order terms, the linearization of equations \re{eq:bal-m}, \re{eq:NS} and \re{eq:heat-eq} via substituting equations \re{eq:Fou} and \re{eq:New} are}
\begin{align}
    \label{eq:llin-1}
    &\pdt{\qrho} + \brho \nab \cdot \qqv = 0 , \\
    \label{eq:llin-2}
    \begin{split}
        &\brho \pdt{\qqv} = - \pder{p}{T}{\qrho}^0 \nab T - \pder{p}{\qrho}{T}^0 \nab \qrho \\
        & \hskip 2.5ex + \bbeta \left[ \left( \bReta + \frac{4}{3} \right) \nab \left( \nab \cdot \qqv \right) - \nab \times \nab \times \qqv \right],
    \end{split} \\
    \label{eq:llin-3}
    &\brho \pder{u}{T}{\qrho}^0 \pdt{T} = \blam \Lapl T - \bT \pder{p}{T}{\qrho}^0 \nab \cdot \qqv ,
\end{align}
where $ ^0 $ denotes the value of a quantity in the equilibrium state and $ \reta = \frac{\eta_{\rm Vol}}{\eta} $ is the ratio of the viscosities.

Via Helmholtz decomposition, the velocity field is given as a sum of an irrotational (\ie curl-free) and a solenoidal (\ie divergence-free) vector fields, \ie
\begin{align}
    \label{eq:Helmholtz}
    \begin{split}
        & \qqv = \vpar + \vper = \nabla \varphi + \nabla \times \Tensor{\psi} \\
        & \text{with} \qquad
        \nab \cdot \qqv = \nab \cdot \vpar = \Lapl \varphi \\
        & \text{and} \qquad
        \nab \times \qqv = \nab \times \vper = - \Lapl \Tensor{\psi} ,
        \end{split}
\end{align}
via scalar and vector potentials $ \varphi $ and $ \Tensor{\psi} $. Introducing the vorticity
\begin{align}
    \qqom = \nab \times \qqv = \nab \times \vper ,
\end{align}
the decomposed equations corresponding to equations \re{eq:llin-1}--\re{eq:llin-3} are
\begin{align}
    \label{eq:lin-H-1}
    &\pdt{\qrho} + \brho \Lapl \varphi = 0 , \\
    \label{eq:lin-H-2}
    &\brho \pdt{\varphi} = - \pder{p}{T}{\qrho}^0
    {
    \left( T - \bT \right)} - \pder{p}{\qrho}{T}^0
    {
    \left( \qrho - \brho \right)} + \bbeta \left( \bReta + \frac{4}{3} \right) \Lapl \varphi , \\
    \label{eq:lin-H-3}
    &\brho \pder{u}{T}{\qrho}^0 \pdt{T} = \blam \Lapl T - \bT \pder{p}{T}{\qrho}^0 \Lapl \varphi , \\
    \label{eq:lin-H-4}
    &\brho \pdt{\qqom} = \bbeta \Lapl \qqom .
\end{align}
One can observe that equation \re{eq:lin-H-4} -- which expresses the linearized vorticity equation -- is decoupled from equations \re{eq:lin-H-1}--\re{eq:lin-H-3}, therefore, in the linear approximation, the dynamics of $ \qqom $ does not affect the dynamics of $ \qrho $, $ \varphi $ and $ T $, and vice versa. Inserting \re{eq:lin-H-1} into the time derivative of \re{eq:lin-H-2}, the generalized wave equation
\begin{align}
    \label{eq:gen-wave}
    \brho \ppdt{\varphi}{2} &= \brho \pder{p}{\qrho}{T}^0 \Lapl \varphi - \pder{p}{T}{\qrho}^0 \pdt{T} + \bbeta \left( \bReta + \frac{4}{3} \right) \Lapl \pdt{\varphi}
\end{align}
is obtained. Therefore, the linear approximation of the coupled thermomechanical process is fully characterized by the coupled system of the generalized heat conduction equation \re{eq:lin-H-3} and the generalized wave equation \re{eq:gen-wave}.

\subsection{On thermodynamic compatibility} \label{sec:therm-comp}

In general, balance equation of internal energy \re{eq:bal-u} -- via Gibbs relation \re{eq:Gibbs} and balance equation of mass \re{eq:bal-m} -- can be reformulated as
\begin{align}
    \label{eq:u->s}
    \qrho T \Dv s = - \nab \cdot \qqq - \qqPi : \Sym{\left( \qqv \otimes \nab \right)} .
\end{align}
The partial derivatives of $ s ( T , \qrho ) $ are
\begin{align}
    \pder{s}{T}{\qrho} &= \frac{1}{T} \pder{u}{T}{\qrho} , &
    \pder{s}{\qrho}{T} &= - \frac{1}{\qrho^2} \pder{p}{T}{\qrho} ,
\end{align}
which can be simply shown by considering $ s ( T , \qrho ) = s \big( u ( T , \qrho ) , \qrho \big) $ and applying \re{eq:Gibbs} and \re{eq:urT-pTr}. As a consequence, the linearized heat conduction equation \re{eq:lin-H-3} can be expressed via specific entropy as
\begin{align}
    \label{eq:lin-H-3-s}
    \brho \bT \pdt{s} =  - \nab \cdot \left( - \blam \nab T \right) .
\end{align}
Comparing \re{eq:u->s} to \re{eq:lin-H-3-s} points out that via linearization, the entropy production contribution of viscous effects is neglected, hence, in the linear approximation, viscous force density attenuates velocity amplitude, however, it preserves entropy, which, from the point of view of thermodynamics, is a contradiction. 
{Nevertheless, this is a side effect of linearization, a widespread and established method. In the case of the piston effect -- the phenomenon at the center of our study --, the absence of this damping term may predict a 
somewhat different temperature rise at the beginning of the process. Meanwhile, in pressure-driven problems, it may produce more significant deviations.}

\subsection{Isothermal and isentropic wave propagation}

In case we neglect all dissipative effects, \ie heat conduction and viscous momentum transport, taking $ \blam = 0 $ and $ \bbeta = 0 $, then, substituting \re{eq:lin-H-3} into \re{eq:gen-wave}, the wave equation 
\begin{align}
    \label{eq:wave-eq}
    \ppdt{\varphi}{2} &= \pder{p}{\qrho}{T}^0 \left[ 1 + \frac{\bT}{\left( \brho \right)^2} \frac{\left( \pder{p}{T}{\qrho}^0 \right)^2}{\pder{u}{T}{\qrho}^0 \pder{p}{\qrho}{T}^0} \right] \Lapl \varphi 
\end{align}
can be recognized, where the expression in the bracket is actually the specific heat ratio defined in \re{eq:gam}, hence, the wave propagation velocity is
\begin{align}
    \label{eq:as}
    \bas = \sqrt{\bgam \pder{p}{\qrho}{T}^0} \stackrel{\re{eq:kT}}{=} \sqrt{\frac{\bgam}{\brho \bkT}}
\end{align}
and is called the \textit{isentropic} speed of sound, since \re{eq:lin-H-3} with neglected heat conduction prescribes an isentropic process [\cf \re{eq:lin-H-3-s}]. When $ \pder{p}{T}{\qrho}^0 = 0 $, then thermal expansion is neglected [\cf \re{eq:bp}], hence \re{eq:lin-H-3} with neglected heat conduction characterizes an isothermal process and the wave propagation velocity is the \textit{isothermal} speed of sound
\begin{align}
    \label{eq:aT}
    \baT := \sqrt{ \pder{p}{\qrho}{T}^0 } .
\end{align}
Therefore, the relationship among the isentropic and isothermal wave propagation velocities is
\begin{align}
    \frac{\left( \bas \right)^2}{\left( \baT \right)^2} = \bgam \ge 1 .
\end{align}

\subsection{Coupling of acoustic and thermal processes -- the effect of thermal expansion}

We proceed by identifying the relevant combinations of constants parametrizing the linearized set of equations. With the help of the material properties, equations \re{eq:gen-wave} and \re{eq:lin-H-3} can be reformulated as
\begin{align}
    \label{eq:lin-wave-mat}
    \ppdt{\varphi}{2} &= \frac{\left( \bas \right)^2}{\bgam} \Lapl \varphi - \bbp \frac{\left( \bas \right)^2}{\bgam} \pdt{T} + \bnu \left( \bReta + \frac{4}{3} \right) \Lapl \pdt{\varphi} , \\
    \label{eq:lin-heat-mat}
    \pdt{T} &= \bgam \ba \Lapl T - \bT \frac{\bbp \left( \bas \right)^2}{\bcp} \Lapl \varphi ,
\end{align}
where $ \bnu = \frac{\bbeta}{\brho} $ and $ \ba = \frac{\blam}{\brho \bcp} $ are the kinematic viscosity and the thermal diffusivity, respectively.

Neglected thermal expansion leads to the decoupling of equations \re{eq:lin-wave-mat} and \re{eq:lin-heat-mat} and yields the viscous damped acoustic wave equation and the heat equation. In such a case, since $ \bbp = 0 $, we have $ \bgam = 1 $ and $ \bcv = \bcp $. It is then apparent that, in the linear approximation, thermal expansion is the only mechanism that can lead to the coupling of acoustic and thermal processes. If neither viscous and thermal attenuations nor thermal expansion are neglected, then, rearranging equation \re{eq:lin-heat-mat} into the operator form
\begin{align}
    \left( \pdt{} - \bgam \ba \Lapl \right) T = - \bT \frac{\bbp \left( \bas \right)^2}{\bcp} \Lapl \varphi
\end{align}
and acting on \re{eq:lin-wave-mat} with the operator $ \pdt{} - \bgam \ba \Lapl $, one obtains 
{
\small  
\begin{align}
    \ppdt{\varphi}{3} -
    \FTminus  
    \left[ \bgam \ba + \bnu
    \FTminus  
    \left( \bReta + \frac{4}{3} \right)
    \FTminus  
    \right]
    \FTminus  
    \Lapl \ppdt{\varphi}{2} -
    \left( \bas \right)^2 \Lapl \pdt{\varphi} + \bgam \ba \bnu
    \FTminus  
    \left( \bReta + \frac{4}{3} \right)
    \FTminus  
    \Lapl \Lapl \pdt{\varphi} + \left( \bas \right)^2 \ba \Lapl \Lapl \varphi = 0 ,
\end{align}
}
which -- after rearranging -- results in a hierarchical\footnote{Such a hierarchy is a typical outcome whenever one eliminates some of the degrees of freedom, see, \eg \cite{van2015thermodynamic,fulop2018emergence}.} wave equation of modified heat conduction equations, \ie
{
\small 
\begin{align}
    \label{eq:th-ac-wave-eq}
    \ppdt{}{2} \left\{ \pdt{\varphi} - \ba \left[ \bgam + \frac{\bnu}{\ba} \left( \bReta + \frac{4}{3} \right) \right] \Lapl \varphi\right\} = \left( \bas \right)^2 \Lapl \left\{ \pdt{\varphi} - \ba \left[ 1 + \bgam \frac{\bnu}{\left( \bas \right)^2} \left( \bReta + \frac{4}{3} \right) \pdt{} \right] \Lapl \varphi \right\} .
\end{align}
}
In equation \re{eq:th-ac-wave-eq}, both hyperbolic and parabolic characteristics are present, thus the usually applied numerical time integration methods may fail or may provide solutions with strong numerical artifacts. Furthermore, realizing boundary conditions prescribed on thermal quantities is complicated, and vice versa, \ie the same equation as \re{eq:th-ac-wave-eq} can be derived for the temperature field, and then realizing flow boundary conditions is difficult.

The quantitative analysis of coupling through thermal expansion can be performed via rescaling equations \re{eq:lin-wave-mat} and \re{eq:lin-heat-mat}. Via the thermal velocity of particle motion $ \bar{\qv}_{\rm th}^0 = \sqrt{\bcp \bT} $, the mean free path related characteristic length scale $ \frac{\bnu}{\bar{\qv}_{\rm th}^0} $ and the corresponding characteristic time scale $ \frac{\bnu}{\left( \bar{\qv}_{\rm th}^0 \right)^2} $ are introduced, via which the non-dimensional time and space variables
and the corresponding non-dimensional derivatives are
\begin{align}
&&&&
    \tau &= \frac{t}{\frac{\bnu}{\left( \bar{\qv}_{\rm th}^0 \right)^2}} , &
    \Tensor{\rho} &= \frac{\qqr}{\frac{\bnu}{\bar{\qv}_{\rm th}^0}} ,
&&&&
\\
&&&&
    \pdt{} &= \frac{\left( \bar{\qv}_{\rm th}^0 \right)^2}{\bnu} \pdct{} , &
    \nabla{} &= \frac{\bar{\qv}_{\rm th}^0}{\bnu} \tilde{\nabla} .
&&&&
\end{align}
In parallel, non-dimensional velocity potential and temperature are defined as
\begin{align}
    \phi &= \frac{\varphi}{\bnu} , &
    \theta &= \frac{T}{\bT} .
\end{align}
Correspondingly, the dimensionless equivalents of \re{eq:lin-wave-mat} and \re{eq:lin-heat-mat} are
\begin{align}
    \label{eq:lin-wave-nondim}
    \ppdct{\phi} &= \frac{\Gamma^0}{\bgam \qB} \tilde{\Lapl} \phi - \frac{\Gamma^0}{\bgam} \pdct{\theta} + \left( \bReta + \frac{4}{3} \right) \tilde{\Lapl} \pdct{\phi} , \\
    \label{eq:lin-heat-nondim}
    \pdct{\theta} &= \frac{\bgam}{\qPr} \tilde{\Lapl} \theta - \Gamma^0 \tilde{\Lapl} \phi ,
\end{align}
where $ \qB = \bbp \bT $ is the non-dimensional thermal expansion coefficient, $ \Gamma^0 = \frac{\bbp \left( \bas \right)^2}{\bcp} $ is the Grüneisen parameter, and $ \qPr = \frac{\bnu}{\ba} $ is the Prandtl number. Thermophysical parameters and calculated dimensionless coefficients appearing in \re{eq:lin-wave-nondim} and \re{eq:lin-heat-nondim} for water and carbon dioxide in some unique thermodynamic states are given in Table~\ref{tab} and \ref{tab2}. 
{Although the non-dimensional parameters characterizing the dynamics vary from material to material and from state to state, one may expect similar magnitudes for their ratios in the same states of matter of different materials, due to the common properties of liquids vs. gases (\eg ``small'' vs.~``large'' compressibility, ``large'' vs.~``small'' thermal conductivity). This statement is related to the ``Theorem of corresponding states''. Therefore, we present these parameters in Tables \ref{tab} and \ref{tab2} for liquid, gas, saturated liquid, saturated steam, and supercritical fluid states, respectively. Not the same material is displayed for these five cases because, experimentally, it is easy to study certain states with water and other states with carbon dioxide. Finally, we note that our numerical investigations (see Sec.~\ref{sec:piston}) were performed for supercritical and nearly ideal gas state of carbon dioxide for the values reported in Table~\ref{tab} and \ref{tab2}.} Since viscous momentum transport is 
uninteresting
from the point of view of thermal expansion coupling, it is neglected in what follows.

A first level of estimating the remarkableness of the various terms in \re{eq:lin-wave-nondim}--\re{eq:lin-heat-nondim} is when one compares the emerged dimensionless coefficient combination factors. Starting with \re{eq:lin-wave-nondim}, contrasting $ \frac{\Gamma^0}{\bgam \qB} $ with $ \frac{\Gamma^0}{\bgam} $, it can be stated from the value of the dimensionless thermal expansion coefficient $ \qB $ that both in liquid and saturated liquid states the first, mechanical, term of \eqref{eq:lin-wave-nondim} is expected to dominate over the second, thermal expansion, one. In such situations, the effect of thermal expansion on mechanical processes can be neglected. On the other side, for vapor and saturated vapor states, $ \qB $ is around 1, thus the two terms in \re{eq:lin-wave-nondim} appear to be comparable, \ie both terms should be considered to determine the dynamics. Third, near the liquid-vapor critical point, thermal expansion dominates. This latter statement explains why considering only thermal expansion and neglecting acoustic propagation results a good approximation for heat conduction near the critical point, as done in \cite{boukari1990critical}. Turning towards the other equation, the non-dimensional heat equation \re{eq:lin-heat-nondim}, heat conduction appears to always dominate over the effect of thermal expansion.

\begin{table*}[!t]  
    \caption{Thermophysical parameters for water and carbon dioxide in some unique thermodynamic states. Numerical values -- where reference is not indicated -- are taken from the NIST database \cite{lemmon2024thermophysical}. Legend: L -- liquid, V -- vapor, sL -- saturated liquid, sV -- saturated vapor, SC -- supercritical.}
    \label{tab}
    \begin{center}
        {
            \begin{tabular}{ l l c r r r r r l r r r}
                \rule[-3ex]{0em}{1ex}
                Material & \multicolumn{1}{c}{%
                \hspace{-.5em}
                $ \Qds \frac{\bT}{ \mathrm{K} } $} &
                \multicolumn{1}{c}{%
                \makebox[1.5em][c]
                {$ \Qds \frac{\bpp}{ \mathrm{MPa} } $}} &
                \multicolumn{1}{c}{$ \Qds \frac{\bcp}{ \Nicefrac{\mathrm{J}}{\mathrm{kg
                \,
                K}}}  $} &
                \multicolumn{1}{c}{$ \Qds \frac{\bas}{ \Nicefrac{\mathrm{m}}{\mathrm{s}}}  $} &
                \multicolumn{1}{c}{$ \Qds \frac{\bbp}{ 10^{-3} \, \Nicefrac{1}{\mathrm{K}}}  $} &
                \multicolumn{1}{c}{$ \Qds \bgam $} &
                \multicolumn{1}{c}{%
                \hspace{-.5em}
                $ \Qds \frac{\bnu}{ 10^{-7} \, \Nicefrac{\mathrm{m}^2}{\mathrm{s}}} $} &
                \multicolumn{1}{c}{%
                \hspace{-1.em}
                $ \Qds \bReta $
                } &
                \multicolumn{1}{c}{%
                \hspace{-.4em}
                $ \Qds \frac{\ba}{ 10^{-7} \, \Nicefrac{\mathrm{m}^2}{\mathrm{s}}} $}
                \\  \hline
                L-H$_2$O & 280 & 0.1 & \hphantom{1}4200.945 & 1434.274 & 0.046 & $ \FTapprox 1.000 $ & 14.337 & $ \FTapprox 3 $\FThcite
                \cite{holmes2011temperature} & 1.362 \\
                sL-H$_2$O & 372.76 & 0.1 & \hphantom{1}4215.223 & 1543.501 & 0.753 & 1.118 & 2.950 & N/A & 1.676 \\
                sV-H$_2$O & 372.76 & 0.1 & \hphantom{1}2078.449 & 471.994 & 2.903 & 1.337 & 206.965 & N/A & 199.941  \\
                V-CO$_2$ & 305 & 0.1 & \hphantom{1}857.314 & 271.433 & 3.329 & 1.291 & 87.422 & $ \FTapprox 0.4 $\FThcite
                \cite{wang2019bulk} & 114.791 \\
                SC-CO$_2$ & 305 & 7.4 & 16328.205 & 184.164 & 136.873 & 12.868 & 0.729 & $ \FTapprox 6 $\FThcite
                \cite{hasan2012thermoacoustic} & 0.126
            \end{tabular}
        }
    \end{center}
\end{table*}

\begin{table}[!t]  
    \caption{Dimensionless coefficients appearing in \re{eq:lin-wave-nondim} and \re{eq:lin-heat-nondim}, calculated from the values given in Table~\ref{tab}.}
    \label{tab2}
    \begin{center}
        {
            \begin{tabular}{ l l c r r r r}
                \rule[-3ex]{0em}{1ex}
                Material & \multicolumn{1}{c}{%
                \hspace{-.5em}
                $ \Qds \frac{\bT}{ \mathrm{K} } $} &
                \multicolumn{1}{c}{%
                \makebox[1.5em][c]
                {$ \Qds \frac{\bpp}{ \mathrm{MPa} } $}} &
                \multicolumn{1}{c}{$ \qB  $} &
                \multicolumn{1}{c}{$ \Gamma^0 $} &
                \multicolumn{1}{c}{$ \qPr $} &
                \multicolumn{1}{c}{$ \frac{\bgam}{\qPr} $} 
                \\  \hline
                L-H$_2$O & 280 & 0.1 & 0.013 & 0.023 & 10.529 & 0.095 \\
                sL-H$_2$O & 372.76 & 0.1 & 0.281 & 0.426 & 1.671 & 0.635 \\
                sV-H$_2$O & 372.76 & 0.1 & 1.082 & 0.311 & 1.016 & 1.316 \\
                V-CO$_2$ & 305 & 0.1 & 1.015 & 0.286 & 0.769 & 1.679 \\
                SC-CO$_2$ & 305 & 7.4 & 41.746 & 0.284 & 5.794 & 2.221
            \end{tabular}
        }
    \end{center}
\end{table}

However, this first-level analysis does not incorporate the possible further time scale(s) introduced by external effects, manifesting themselves in the form of time-dependent boundary conditions. Especially, abrupt changes can induce large time derivatives, and then the picture becomes more complex, the degrees of freedom become considerably intertwined, and all the analyzed terms have to be kept. Consequently, a purely mechanical excitation can induce nontrivial thermal changes, and a thermal shock can originate a spectacular mechanical implication.

Our current goal is to investigate linearized heat conduction in fluids by taking into account the effect of thermal expansion as well. We would like to track both large and small time scales, without assuming that either acoustic or thermal processes dominate. From numerical point of view, this task requires a reliable time-integration algorithm since wave-like and diffusive signal propagation -- two rather different kinds of phenomena -- take place simultaneously. Our developed numerical scheme presented in the next section addresses these demands appropriately. In details, the splitting of the semi-discretized equations into reversible and irreversible parts requires the time integration of wave and heat equations in different sub time steps. Therefore, choosing for the relationship between the time step and the space step of the scheme for stability and  minimized numerical artifacts can be inherited from the classical methods for hyperbolic and parabolic equations.

\section{The numerical scheme} \label{sec:num-scheme}

Let us investigate the thermoacoustic processes along the axial direction of a pipe. Neglecting heat conduction and fluid flow in the cross-section directions, the problem can be treated in a single spatial dimension. The corresponding Cartesian coordinate is denoted by $ x $. Instead of directly solving equations \re{eq:lin-wave-mat} and \re{eq:lin-heat-mat} for the scalar potential and temperature, the whole system of linearized balance equations -- given together with Fourier's heat conduction law and Newton's law of viscosity -- will be solved for all fields. In this problem, there exists a macroscopic length scale, namely, the pipe length $ X $, via which and the isentropic speed of sound, one can introduce the non-dimensional time and space coordinates and the corresponding non-dimensional derivatives as
\begin{align}
    \hht &= \frac{t}{\frac{X}{\bas}} , &
    \hx &= \frac{x}{X} , \\
    \pdt{} &= \frac{\bas}{X} \pdht{} , &
    \pdx{} &= \frac{1}{X} \pdhx{} .
\end{align}
Moreover, the non-dimensional fields are introduced as
\begin{align}
    \begin{split}
        \hT &= \frac{T}{\cT} , \qquad
        \hrho = \frac{\qrho}{\crho} , \qquad
        \hv= \frac{\qv}{\bas} , \\
        \hq &= \bgam \frac{\qq}{\crho \cT \bcp \bas} , \hskip 8ex
        \hPi = \frac{\Pi}{\crho \left( \bas \right)^2} ,
    \end{split}
\end{align}
where $ \cT $ and $ \crho $ denote the critical temperature and density of the investigated fluid. Therefore, the non-dimensional one spatial dimensional thermoacoustic equations are
\begin{align}
    \pdht{\hrho} + \hrho^0 \pdhx{\hv} &= 0 , \label{eq:dimless1}\\
    \pdht{\hv} &= - \frac{\qB}{\bgam \hT^0} \pdhx{\hT} - \frac{1}{\bgam \hrho^0} \pdhx{\hrho} - \frac{1}{\hrho^0} \pdhx{\hPi} , \\
    \pdht{\hT} &= - \frac{1}{\hrho^0} \pdhx{\hq} - \qB \qEca \pdhx{\hv} , \\
    \label{eq:dimless4}
    \hq &= - \bgam \hrho^0 \frac{1}{\qPr \qRea} \pdhx{\hT} , \\
    \hPi &= - \hrho^0 \frac{1}{\qRea} \left( \bReta + \frac{4}{3} \right) \pdhx{\hv} , \label{eq:dimless5}
\end{align}
where $ \qRea = \frac{\bas X}{\bnu} $ is the acoustic Reynolds number (which together with the Prandtl number $ \qPr $ also determines the acoustic P\'eclet number $ \qPea = \qPr \qRea = \frac{\bas X}{\ba} $) and $ \qEca = \frac{\left( \bas \right)^2}{\bcp \cT} $ is the acoustic Eckert number related to the critical temperature. Non-dimensional pressure is calculated from the determined temperature and density values through the linearized thermal equations of state 
\begin{align}
    \label{eq:p-Trho}
    \hat{p} = \frac{p}{\crho \left( \bas \right)^2} = \hat{p}_0 + \frac{\hrho^0 \qB}{\bgam \hT^0} \left( \hT - \hT^0 \right) + \frac{1}{\bgam} \left( \hrho - \hrho^0 \right) ,
\end{align}
with $ \hat{p}_0 = \frac{1}{\crho \left( \bas \right)^2} p \left( T_0 , \qrho_0 \right) $. Investigated initial conditions correspond to the homogeneous static equilibrium state.

To solve equations \eqref{eq:dimless1}--\eqref{eq:dimless5} numerically, we propose a finite-difference scheme on a staggered space-time grid. The discretization proceeds in two steps: first, the spatial finite difference discretization yields the semi-discrete equations, which are then discretized with respect to time, using a second-order splitting method based on reversible--irreversible splitting.

\subsection{Spatial semi-discretization on a staggered grid}

For the spatial discretization, we use a half-step staggered grid, along the lines of \cite{fulop2020thermodynamical,pozsar2020four,takacs2024thermodynamically}, realizing all spatial derivatives in a central-difference way. This guarantees a second-order spatial accuracy, as shown in \cite{takacs2024thermodynamically}. Following this discretization, the semi-discrete equations read
\begin{align}
    \dht{\hrho_{n+\half}} &= - \hrho^0 \pDhx{\hv_{n+1} - \hv_{n}}, \label{eq:semidiscrete1}\\
    \dht{\hv_{n}} &= - \frac{\qB}{\bgam \hT^0} \pDhx{\hT_{n+\half} - \hT_{n-\half}} - \frac{1}{\bgam \hrho^0} \pDhx{\hrho_{n+\half} - \hrho_{n-\half}} - \frac{1}{\hrho^0} \pDhx{\hPi_{n+\half} - \hPi_{n-\half}} , \label{eq:semidiscrete2} \\
    \dht{\hT_{n+\half}} &= - \frac{1}{\hrho^0} \pDhx{\hq_{n+1} - \hq_{n}} - \qB \qEca \pDhx{\hv_{n+1} - \hv_{n}} ,  \label{eq:semidiscrete3}\\
    \hq_{n} &= - \bgam \hrho^0 \frac{1}{\qPr \qRea} \pDhx{\hT_{n+\half} - \hT_{n-\half}} , \label{eq:semidiscrete4}\\
    \hPi_{n+\half} &= - \hrho^0 \frac{1}{\qRea} \left( \bReta + \frac{4}{3} \right) \pDhx{\hv_{n+1} - \hv_{n}} , \label{eq:semidiscrete5}
\end{align}
where the subscript $n=1,\ldots,N-1$ denotes the index of the grid point located at $\hx=n \hDx$, while $n+\half$ or $n-\half$ denote locations shifted by $+\hDx/2$ or $-\hDx/2$, respectively. The values of $\hv$ and $\hq$ for $n=0$ and $n=N$ are specified directly as boundary conditions. If some different type of boundary condition (e.g. one involving $\hT$, such as an isothermal boundary condition) is needed, the grid can be extended virtually to $n=-\half$ and $n=N+\half$, and similar steps can be taken for more complex boundary conditions. In what follows, we will only need boundary conditions specifying $\hv$ and $\hq$, thus, for further details, we refer to \cite{pozsar2020four}.

It should be noted that, among the semi-discrete equations \eqref{eq:semidiscrete1}--\eqref{eq:semidiscrete5}, equations \eqref{eq:semidiscrete1}--\eqref{eq:semidiscrete3} are ordinary differential equations with respect to time, while \eqref{eq:semidiscrete4}--\eqref{eq:semidiscrete5} are algebraic equations that relate $\hq$ and $\hPi$ to $\hT$ and $\hv$, respectively. Thus, in what follows, we only need to perform numerical time integration for the equations \eqref{eq:semidiscrete1}--\eqref{eq:semidiscrete3}, while \eqref{eq:semidiscrete4}--\eqref{eq:semidiscrete5} can be used to update the remaining fields after each step.

\subsection{Time integration via reversible--irreversible splitting}

For obtaining a qualitatively correct and quantitatively satisfactory numerical solution of equations \eqref{eq:semidiscrete1}--\eqref{eq:semidiscrete5}, an appropriate, structure-preserving numerical method for time integration is needed. Here, we split the generating vector field as
\begin{align}
    \everymath{\displaystyle}
    \dht{}
    \begin{pmatrix}
        \hrho_{n+\half} \VVV
        \hv_{n} \VVV
        \hT_{n+\half}
    \end{pmatrix}
    =
    \underbrace{
    \begin{pmatrix}
    - \hrho^0 \pDhx{\hv_{n+1} - \hv_{n}} \VVV
     - \frac{\qB}{\bgam \hT^0} \pDhx{\hT_{n+\half} - \hT_{n-\half}} - \frac{1}{\bgam \hrho^0} \pDhx{\hrho_{n+\half} - \hrho_{n-\half}}  \VVV
     - \qB \qEca \pDhx{\hv_{n+1} - \hv_{n}}
    \end{pmatrix}
    }_{\Xrev}
    +
    \underbrace{
    \begin{pmatrix}
        0 \VVV
        - \frac{1}{\hrho^0} \pDhx{\hPi_{n+\half} - \hPi_{n-\half}} \VVV
        - \frac{1}{\hrho^0} \pDhx{\hq_{n+1} - \hq_{n}} 
    \end{pmatrix}
    }_{\Xirr} ,
    \label{eq:semidiscrete}
\end{align}
where the original vector field $\X$ is decomposed as $\X = \Xrev + \Xirr$ into a reversible part and an irreversible part,  respectively. Note that the two parts correspond to semi-discretized versions of a hyperbolic  partial differential equation and of a parabolic one, respectively.

A similar approach has already proved to be successful \cite{shang2020structurepreserving}. The main advantage of reversible--irreversible splitting is that different numerical methods can be used for integrating the different vector fields: a symplectic or quasi-symplectic one for the reversible part and a different, generic, one for the irreversible part. This allows for preserving some structure present in the original physical equations, while using an identical numerical method for the different types of processes would destroy that structure, resulting in qualitatively and quantitatively worse results. Identifying reversible and irreversible parts of a vector field is relatively straightforward, while other types of splittings are usually less generalizable (as explored in \cite{mclahlan2002splitting}). Additionally, the above-identified hyperbolic and parabolic nature of the corresponding PDE-s makes this splitting a natural choice: the splitting also preserves this distinction.

Here, we use a quasi-symplectic method (inspired by our previous work, \cite{pozsar2020four}) to integrate the reversible part, and the explicit midpoint method for the irreversible part. For a time step $\hDt$, these maps are denoted as $\phirev(\hDt)$ and $\psiirrev(\hDt)$, respectively. The quasi-symplectic method $\phirev$ preserves the phase-space structure of the reversible equations and is second-order accurate while being explicit, and the midpoint method $\psiirrev$ has been chosen as it is also second-order accurate and explicit.

As will be shown in the next subsection, the composition of substeps
\begin{gather}
    \Phi(\hDt) = \psiirrev(\hDt/2) \circ \phirev(\hDt) \circ \psiirrev(\hDt/2) \label{eq:split}
\end{gather}
is second-order accurate with respect to the time step $\hDt$. Detailed procedure of \re{eq:split} regarding the non-dimensional thermoacoustic problem \re{eq:dimless1}--\re{eq:dimless5} is presented in Appendix~\ref{sec:app-A}.

\subsection{Accuracy} \label{sec:acc}

Besides the qualitative advantages of our scheme given above, its second-order accuracy also yields quantitative advantages in the results. We have already stated that the staggered spatial discretization yields second-order accuracy, and the numerical methods used in $\phirev$ and $\psiirrev$ for the two vector fields are also second-order in time. What remains is to show that the splitting \eqref{eq:split} also yields a second-order accuracy if all substeps are of second-order.

It is well-known that a splitting \eqref{eq:split} is second-order if the steps are exact flows (this is also known as Strang splitting or Marchuk splitting) \cite{strang1968construction,marchuk1968some,leimkuhler2005simulating,hairer2006geometric}, and also if all steps are symmetric numerical schemes of second order \cite{yoshida1990construction,suzuki1990fractal,mclahlan2002splitting}. However, here, $\psiirrev$ is not symmetric. More general order conditions for such methods exist \cite{mclahlan1995numerical}, but we find it instructive to prove the case at hand concisely.

Consider the numerical methods, or maps, $\phirev(\Dt)$ and $\psiirrev(\Dt)$, which are second-order accurate with respect to $\Dt$. Through backward error analysis \cite{reich1999backward}, for any map, there exists a vector field that approximates the flow of that map, to an arbitrary order. As such, these second-order maps are generated from the corresponding distorted vector fields $\tX_{\ldots}$ as
\begin{align}
    \phirev(\hDt) &= \Expof{\hDt \Xrev + \hDt^3 \tXrevthr + \Ord{\hDt^4}}, \\
    \psiirrev(\hDt) &= \Expof{\hDt \Xirr + \hDt^3 \tXirrthr + \Ord{\hDt^4}},
\end{align}
where $\Exp$ is the exponential map. The vector field corresponding to the composite map $\Phi$ can be expressed as
\begin{gather}
    \Phi(\hDt) = \Expof{\hDt \tXone + \hDt^2 \tXtwo + \Ord{\hDt^3}},
\end{gather}
for which, due to \eqref{eq:split},
\begin{align}
    & \Expof{\hDt \tXone + \hDt^2 \tXtwo + \Ord{\hDt^3}}
    \nonumber \\
    & \quad
    \mathrel=
    \Expof{\frac{\hDt}{2} \Xirr + \Ord{\hDt^3}}
    \Exp\Bigl(\hDt \Xrev + \Ord{\hDt^3}\Bigr)
    \Expof{\frac{\hDt}{2} \Xirr + \Ord{\hDt^3}}\label{eq:expprod}
\end{align}
holds. If, for the first distorted term, $\tXone = \Xrev + \Xirr\equiv\X$ holds then the composition map $\Phi$ is  consistent\footnote{A numerical scheme is called consistent if, as $\hDt \to 0$, the numerical solution converges to the exact one.}. If $\tXtwo$ vanishes then the method is indeed second-order.

The right-hand side of equation \eqref{eq:expprod} can be evaluated using the Baker--Campbell--Hausdorff formula \cite{lee2012introduction,mclahlan2002splitting}. A straightforward calculation yields the terms
\begin{align}
    \tXone &= \Xrev + \Xirr \equiv \X , \\
    \tXtwo &= 
    \underbrace{
    \frac{1}{8} \left[\Xirr, \Xirr \right]
    }_{= \; 0}
    +
    \underbrace{
    \frac{1}{4} \left[\Xirr, \Xrev \right] + \frac{1}{4} \left[\Xrev, \Xirr \right]
    }_{= \; 0}\label{eq:bch2}
    \equiv 0 ,
\end{align}
where $[ \, \cdot \, , \, \cdot \, ]$ is the Lie bracket of two vector fields, and we have used its antisymmetric property in \eqref{eq:bch2}.
Hence, the method $\Phi$, given as \eqref{eq:split}, is indeed second-order, and consistent.

\section{Numerical experiments} \label{sec:piston}

Our test problem considers the piston effect in a closed pipe with length $ X $. Thermal excitation is a pulse-like heat current density with pulse duration $ t_{\rm P} $ prescribed at the boundary at $ x = 0 $. Apart from this time interval, the sample is considered to be isolated and adiabatic during the entire process. According to the former statements, the boundary conditions are
\begin{align}
    \qq ( t , x = 0 ) &=
    \begin{cases} \frac{q}{t_{\rm P}} \left[ 1 - \cos \leftf( 2 \pi \frac{t}{t_{\rm P}} \rightf) \right] & \text{if } 0 \le t \le t_{\rm P},
    \\ 0 & \text{otherwise},
    \end{cases} & 
    \qq ( t , x = X ) &= 0 , \\
    \qv ( t , x = 0 ) &= 0 , &
    \qv ( t , x = X ) &= 0 , &
\end{align}
where $ q $ is the amount of cross-section-specific heat introduced during the heat pulse. Therefore, the dimensionless form of the only non-trivial boundary condition is
\begin{align}
    \hq ( \hht , \hx = 0 ) &=
    \begin{cases} \frac{1}{\hrho^0} \frac{\hat{q}}{\hht_{\rm P}} \left[ 1 - \cos \leftf( 2 \pi \frac{\hht}{\hht_{\rm P}} \rightf) \right] & \text{if } 0 \le \hht \le \hht_{\rm P},
    \\ 0 & \text{otherwise},
    \end{cases} 
\end{align}
with
\begin{align}
    \hat{q} &= \bgam \frac{q}{\crho \cT \bcp X} , &
    \hht_{\rm P} &= \frac{t_{\rm P}}{\frac{X}{\bas}} .
\end{align}
The physically free parameters are the initial (homogeneous equilibrium) state $ \left( \bT , \brho \right) $, the length of the pipe $ X $, the cross-sectional heat $ q $ and the pulse duration $ t_{\rm P} $. These parameters -- together with the temperature and density of the critical point of the investigated material -- determine all required non-dimensional parameters.

There exist four time scales in this problem, namely, acoustic time scale $ \tau_{\rm a} = \frac{X}{\bas} $, heat conduction time scale $ \tau_{\rm hc} = \frac{X^2}{\ba} $, viscous time scale $ \tau_{\rm v} = \frac{X^2}{\bnu} $, and pulse duration $ t_{\rm P} $. The ratio of heat conduction and acoustic time scales, as well as the ratio of viscous and acoustic time scales, are related through the P\'eclet and Reynolds numbers, furthermore, heat conduction and viscous time scales are connected through the Prandtl number, \ie
\begin{align}
    \begin{split}
        \qPea &= \frac{\bas X}{\ba} = \frac{X^2}{\ba} \frac{\bas}{X} = \frac{\tau_{\rm hc}}{\tau_{\rm a}} , \\
        \qRea &= \frac{\bas X}{\bnu} = \frac{X^2}{\bnu} \frac{\bas}{X} = \frac{\tau_{\rm v}}{\tau_{\rm a}} , \\
        \qPr &= \frac{\bnu}{\ba} = \frac{\tau_{\rm hc}}{\tau_{\rm v}} .
    \end{split}
\end{align}
In the case of a millimeter- or centimeter-sized sample of CO$_2$, $ \qPea = 10^5 $--$ 10^8 $ (which order of magnitude may probably be generalized for most pure fluids). Typical values of the Prandtl number for liquids and gases vary in the order of magnitudes $ 10^{-1} $--$ 10^2 $. Therefore, in general, we can assume that $ \tau_{\rm a} \ll \tau_{\rm hc} $ and $ \tau_{\rm a} \ll \tau_{\rm v} $. Our numerical calculations presented below support that if $ t_{\rm P } \ll \tau_{\rm a} $ then temperature, density, and pressure disturbances caused by the boundary heat pulse exhibit weakly damped wave propagation. When $ \tau_{\rm a} \ll t_{\rm P } \ll \tau_{\rm hc} $ -- \eg, in case of the piston effect -- then pressure disturbance propagates through the pipe with the isentropic speed of sound, but we expect that reflection will not take place. Thus, the initial, nearly isentropic temperature and density increases on the boundary at $ x = X $ will be followed by diffusion.

In what follows, we investigate numerically the performance of the scheme. Material parameters are considered for carbon dioxide in the supercritical state characterized by $ \bT = 305 $ K and $ p_0 = 7.4 $ MPa. In this state, based on the NIST database \cite{lemmon2024thermophysical}, $ \brho = 321.083 \ {\rm \frac{kg}{m^3}} $, $ \bgam = 12.868 $, $ \qB = 41.744 $, $ \qEca = 0.007 $, $ \qPr = 5.805 $ and $\bReta \approx 6 $ \cite{hasan2012thermoacoustic}. Critical temperature, pressure and density of carbon dioxide are $ \cT = 304.128 $ K, $ p_{\rm c} = 7.377 $ MPa and $ \varrho_{\rm c} = 467.600 \ {\rm \frac{kg}{m^3}} $, respectively.

\subsection{Weakly damped wave propagation}

This section is dedicated to investigating the wave propagation induced by the boundary heat pulse. For this purpose we choose $ \hht_{\rm P} = 0.1 $ and  $ \hat{q} = 0.001 $.

First, we neglect viscosity and perform the calculations with $ \qPea = 10^5 $. The spatial domain is divided into 100 cells. According to our previous experiences for the quasi-symplectic scheme applied to wave propagation, a time step chosen to be exactly on the stability limit yielded a dispersion error-free and dissipation error-free solution, to which a Courant number $ \mathpzc{Co} = \frac{\hDt}{\hDx} = 1 $ corresponded \cite{fulop2020thermodynamical}. However, our newly developed scheme presented here is more complex, thus exact analysis of stability conditions and numerical originated errors are more complicated. Our present numerical investigations recommend a Courant number $ \mathpzc{Co} = 0.95 $, via which time stepping remains stable and dispersion errors seemingly disappear. Numerical experiments have shown that, with an increasing P\'eclet number, the Courant number resulting in a dispersion error-free solution tends to 1. {
Here we note that the Courant number based stability limit corresponds to the stability limit of hyperbolic equations, meanwhile, for small P\'eclet numbers this limit can shift to a parabolic one characterized by $ \frac{\hDt}{\hDx^2} $. Our experience has shown that, for small P\'eclet numbers, wave propagation is suppressed, diminishing any possible dispersion error as well.}

\begin{figure*}[!htb]\centering
    \includegraphics[width=0.8\textwidth]{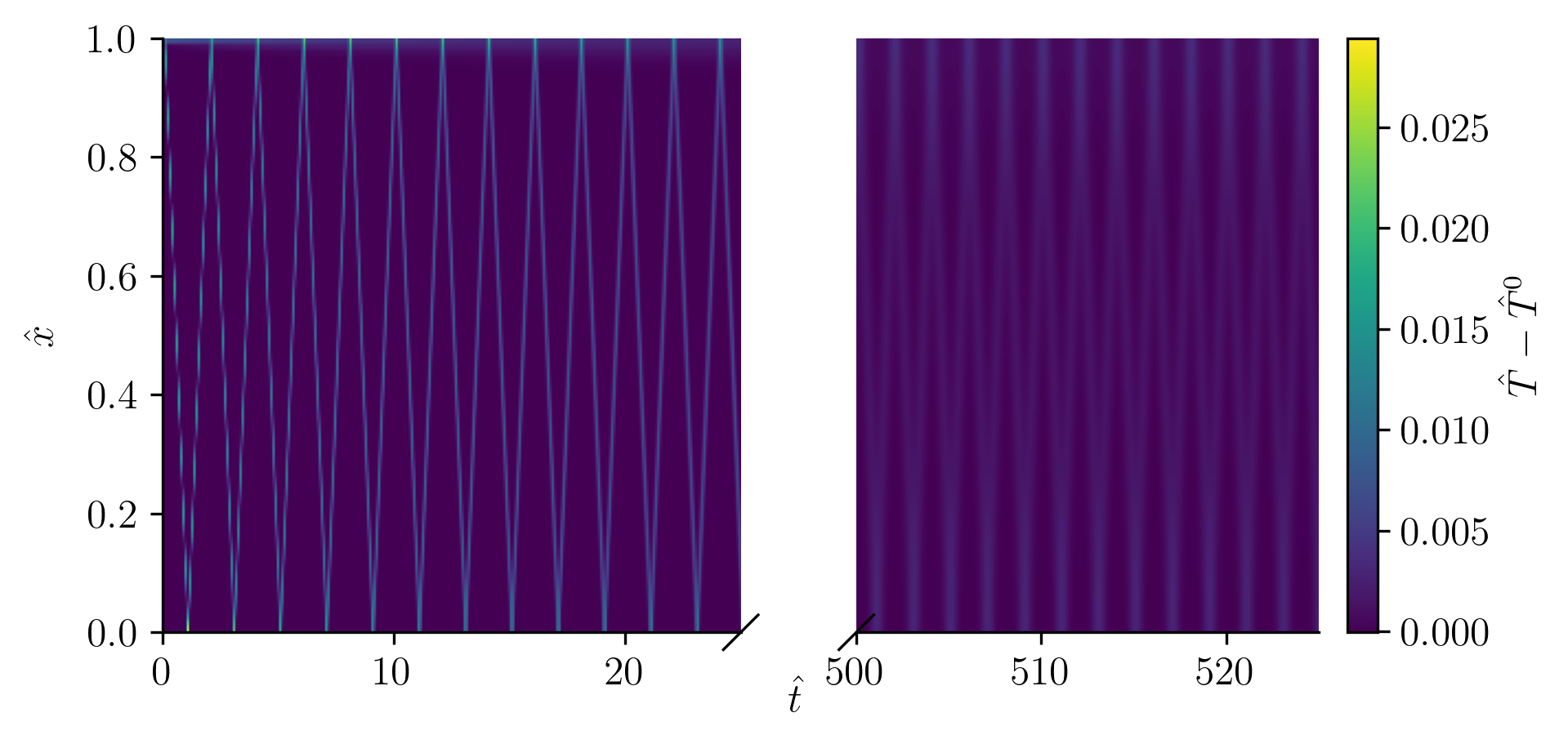}
    \caption{Spacetime plot of the temperature field.}
    \label{fig:spacetime-T}
\end{figure*}

Time evolution of the temperature field is presented as a spacetime plot in Fig.~\ref{fig:spacetime-T} for $ \qPea = 10^5 $. Initially, wave propagation dominates the process, followed by more intensive diffusion over time.

\begin{figure}[!htb]\centering
    \includegraphics[width=0.4\textwidth]{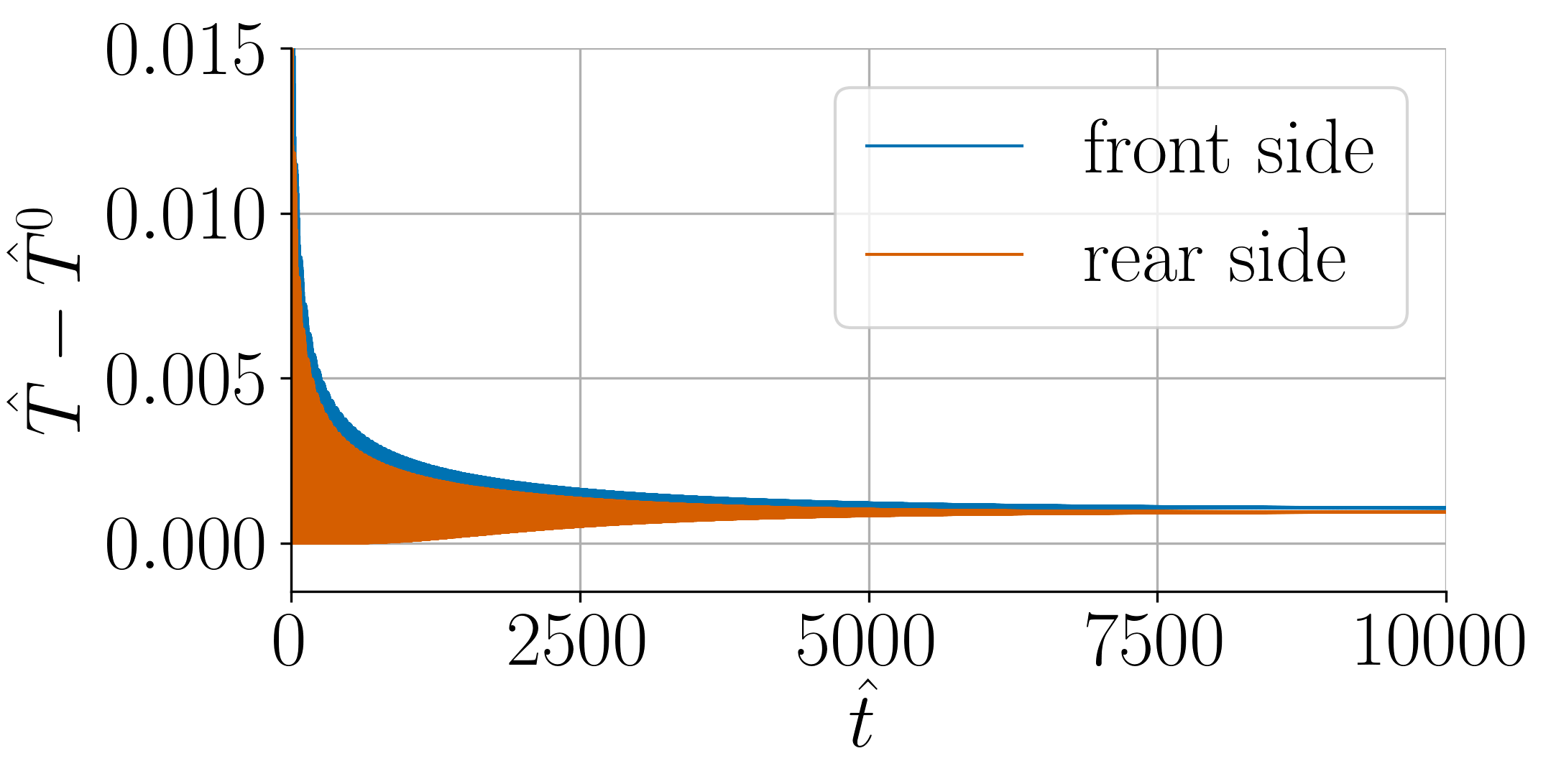}
    \caption{Calculated time evolution of temperature at the front side and at the rear side of the sample.}
    \label{fig:T-t}
\end{figure}

Time evolution of temperature at the front side (\ie where excitation occurs) and at the rear side of the sample is shown in Fig.~\ref{fig:T-t}. As we would expect based on the spacetime diagram, at the end of the investigated time interval, the temperature distribution of the sample can be considered practically homogeneous.

\begin{figure}[!htb]\centering
    \includegraphics[width=0.4\textwidth]{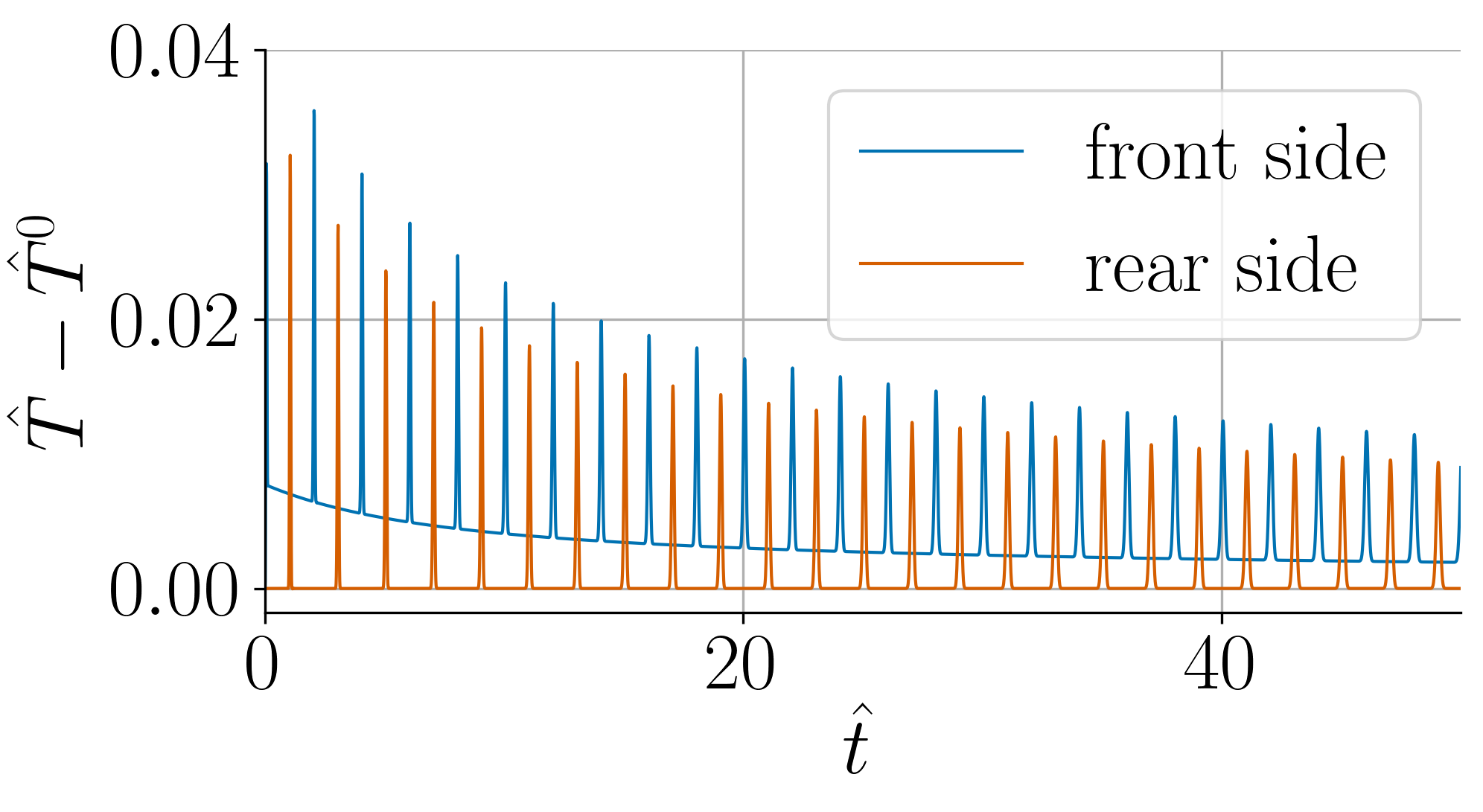}\\
    \hskip -1.3ex \includegraphics[width=0.4\textwidth]{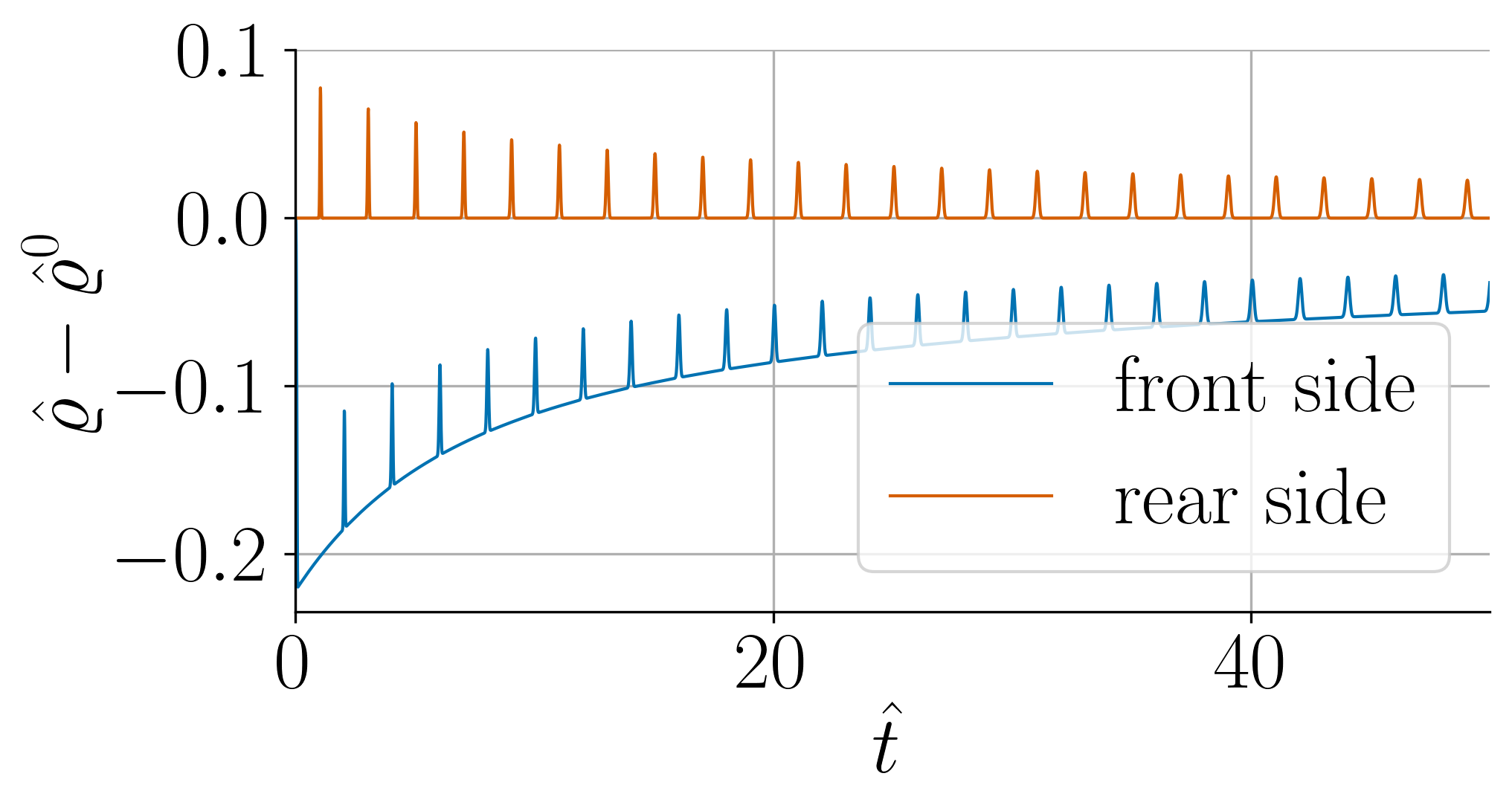}\\
    \includegraphics[width=0.4\textwidth]{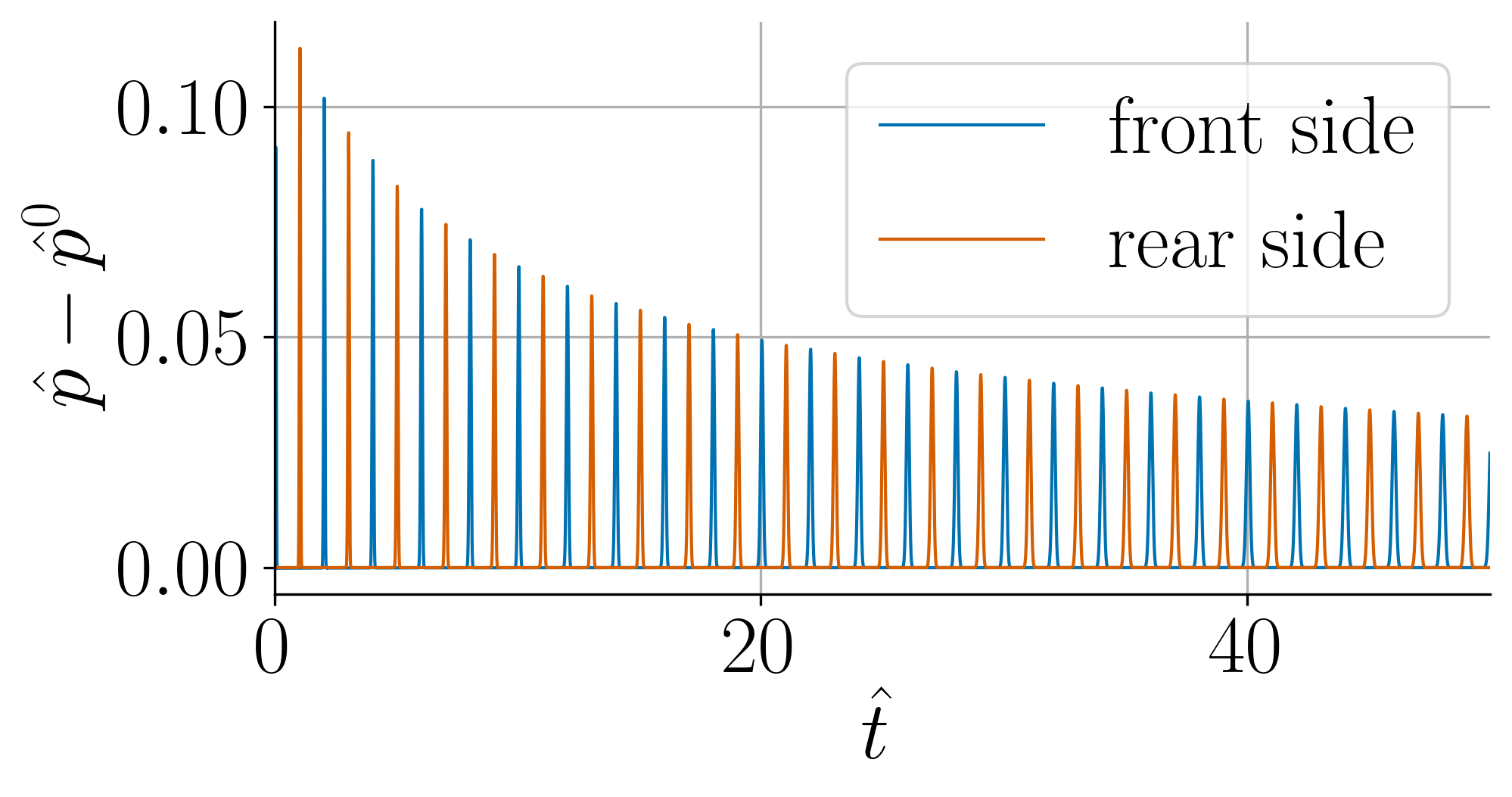}
    \caption{Time evolution of temperature, density, and pressure at the front side  and at the rear side of the sample at the beginning of the process.}
    \label{fig:Trp-t}
\end{figure}

\begin{figure}[!htb]\centering
    \includegraphics[width=0.4\textwidth]{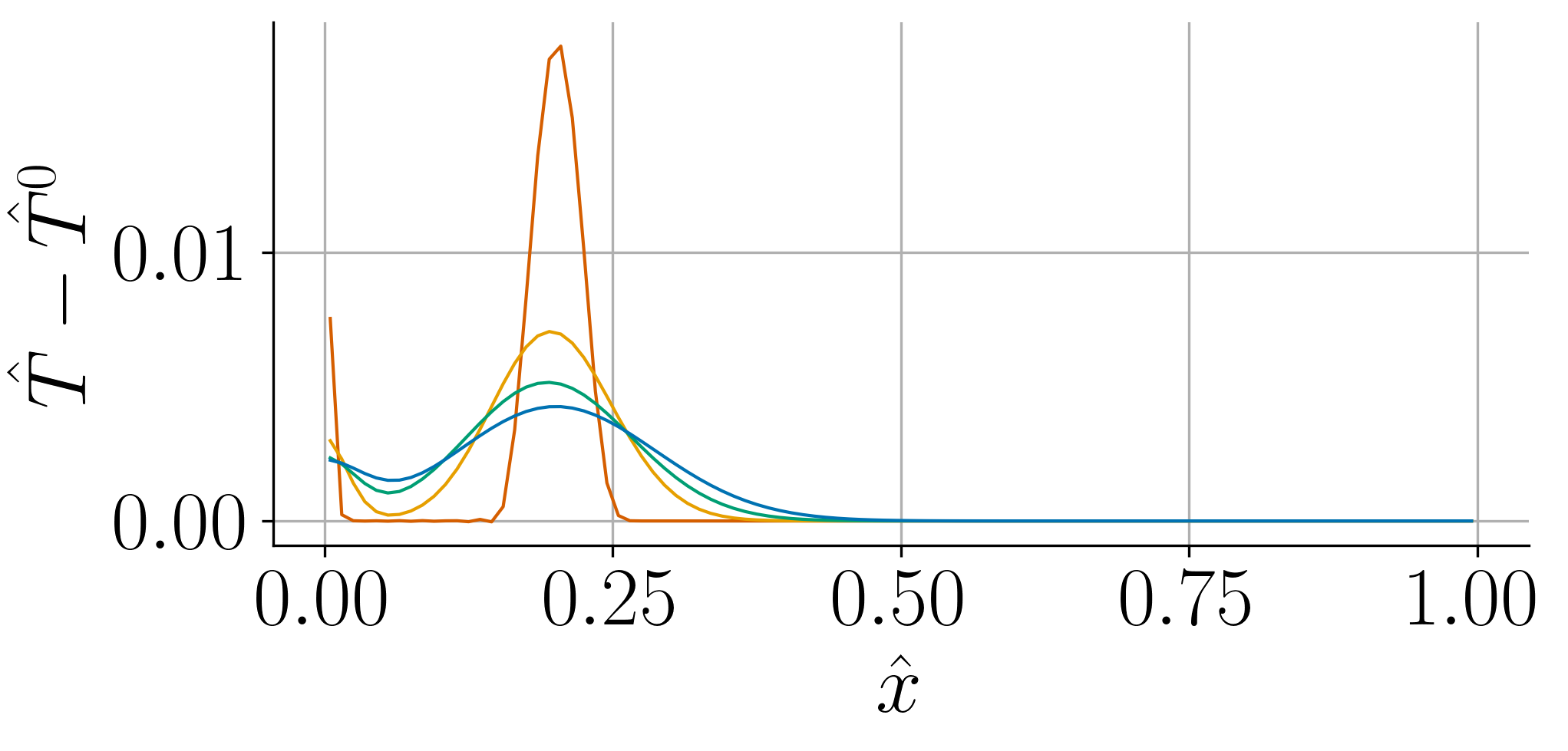}\\
    \hskip -1.3ex \includegraphics[width=0.4\textwidth]{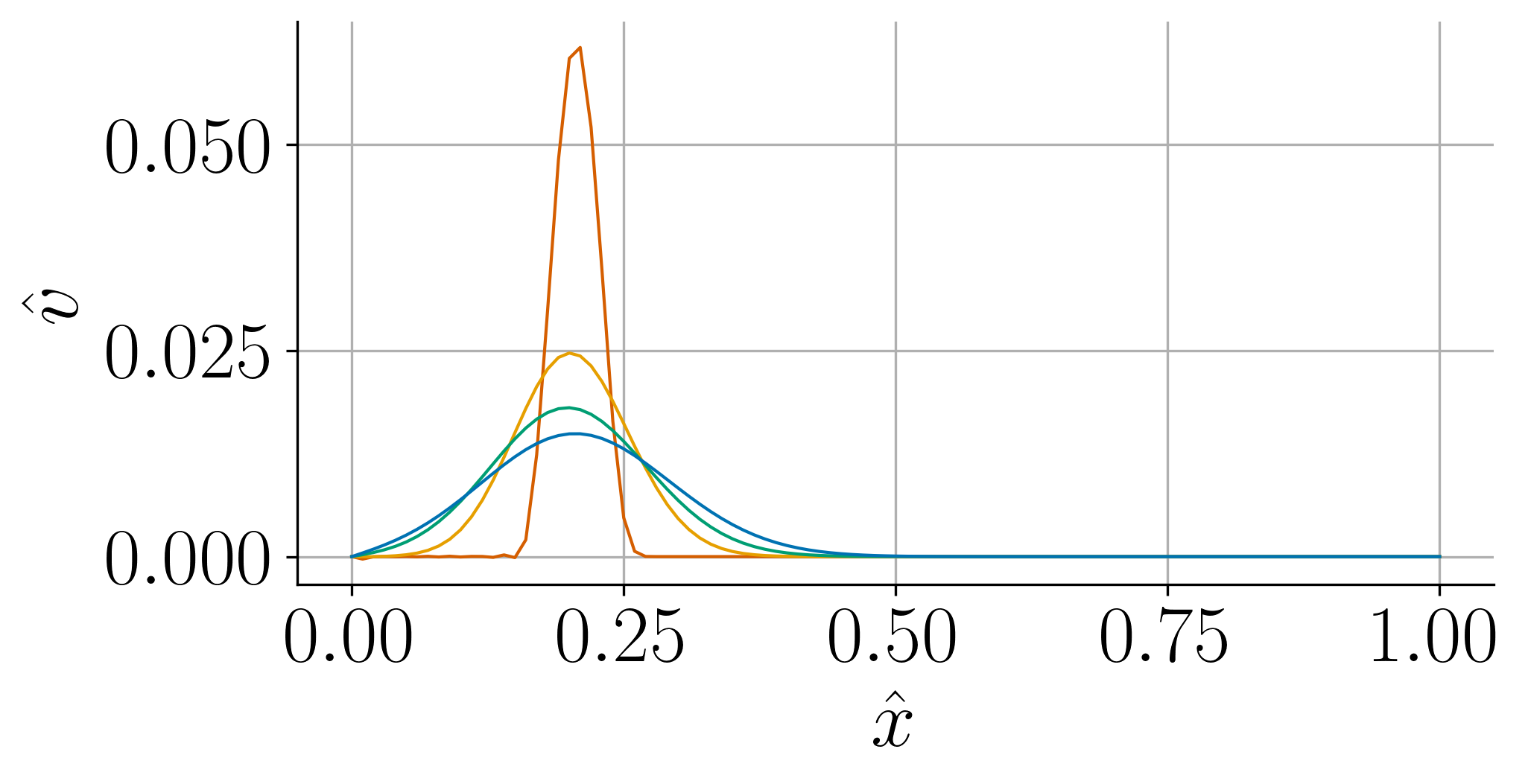}
    \caption{Spatial distributions of temperature and velocity fields along the pipe at the non-dimensional time instants 0.25 (red), 20.25 (orange), 40.25 (green), and 60.25 (blue).}
    \label{fig:Tv-x}
\end{figure}

In Fig.~\ref{fig:T-t}, individual scattering of temperature waves is not visible. Enlarging the beginning of the process, these scatterings are also visible, seemingly, the calculated temperature, density, and pressure responses (presented in Fig.~\ref{fig:Trp-t}) are almost dispersion error-free, which is even more spectacularly illustrated by the spatial distributions of the temperature and velocity field in Fig.~\ref{fig:Tv-x}. 

\begin{figure*}[!htb]
    \centering
    \begin{minipage}{0.45\textwidth}
        \hspace*{2.8ex}\includegraphics[width=0.9\textwidth]{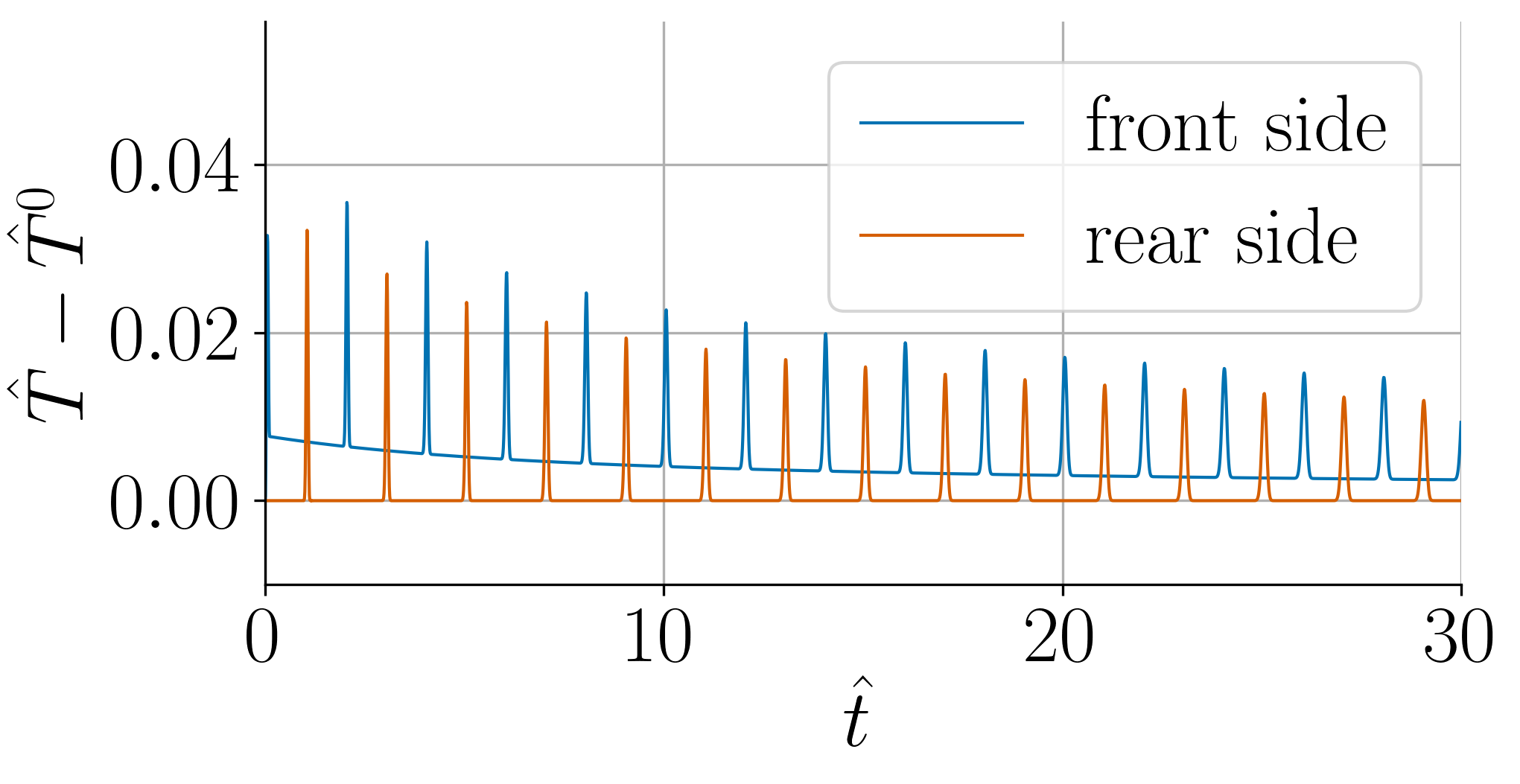}\\
        \includegraphics[width=\textwidth]{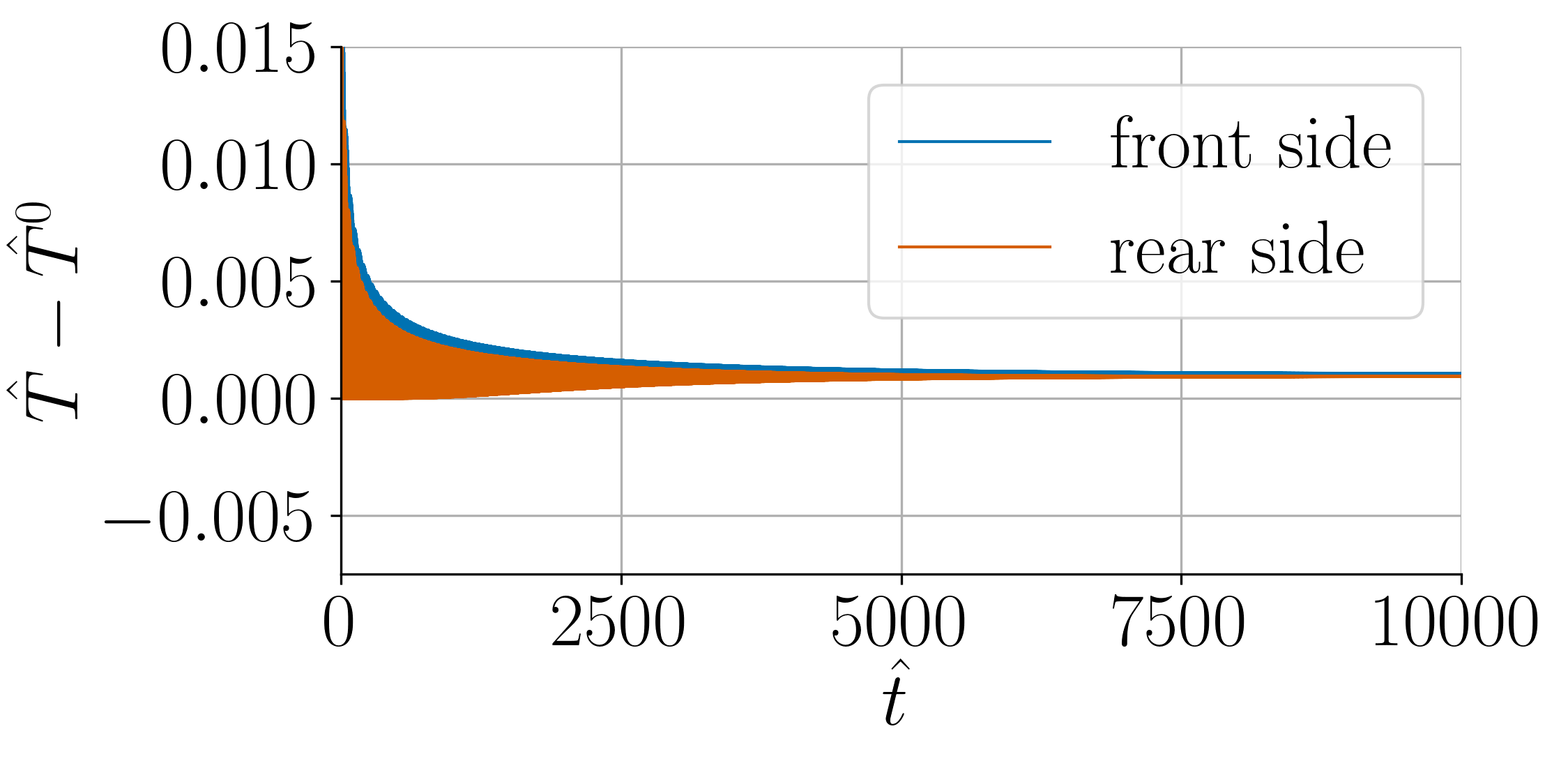}
    \end{minipage}%
    \begin{minipage}{0.45\textwidth}
        \hspace*{2.8ex}\includegraphics[width=0.9\textwidth]{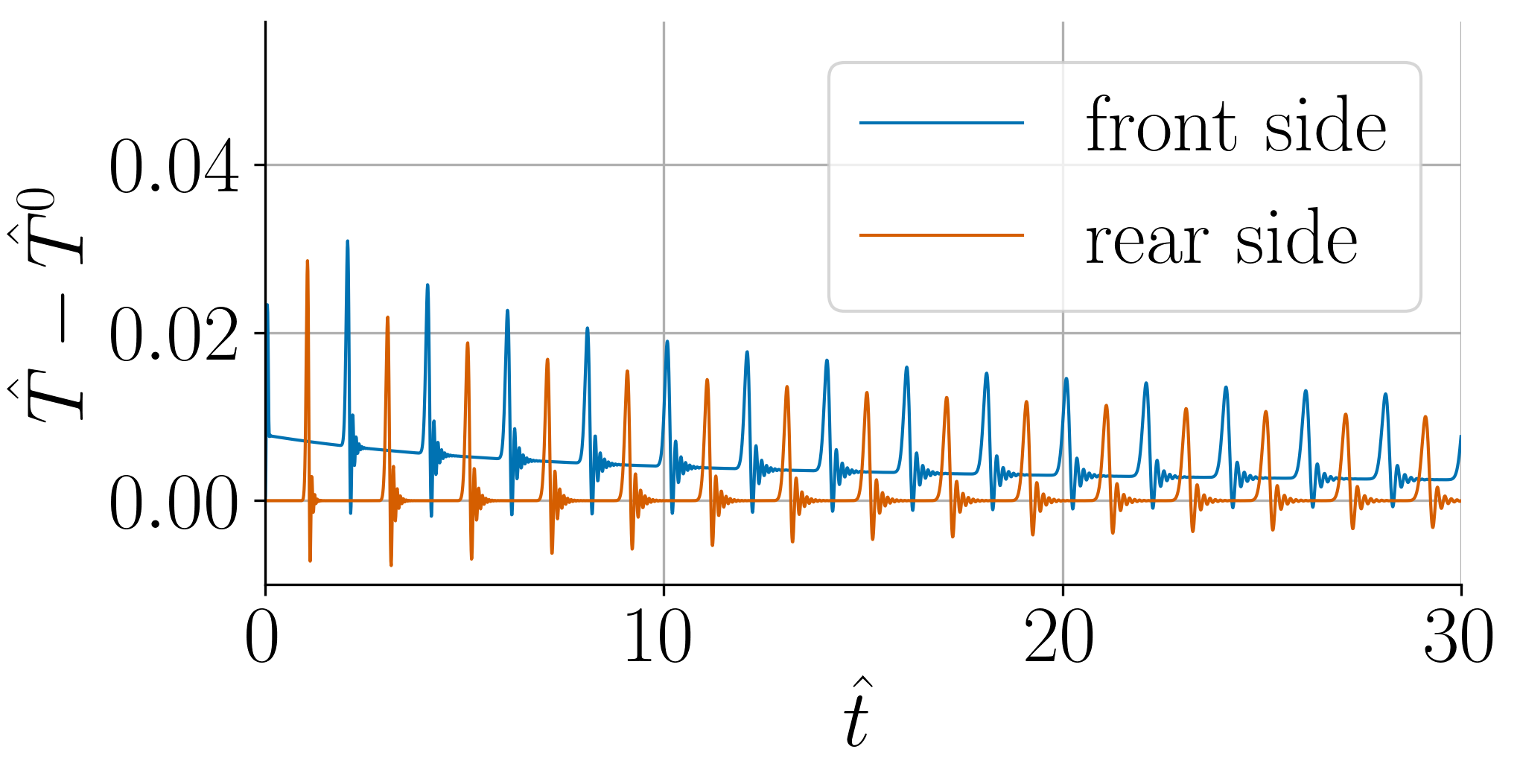}\\
        \includegraphics[width=\textwidth]{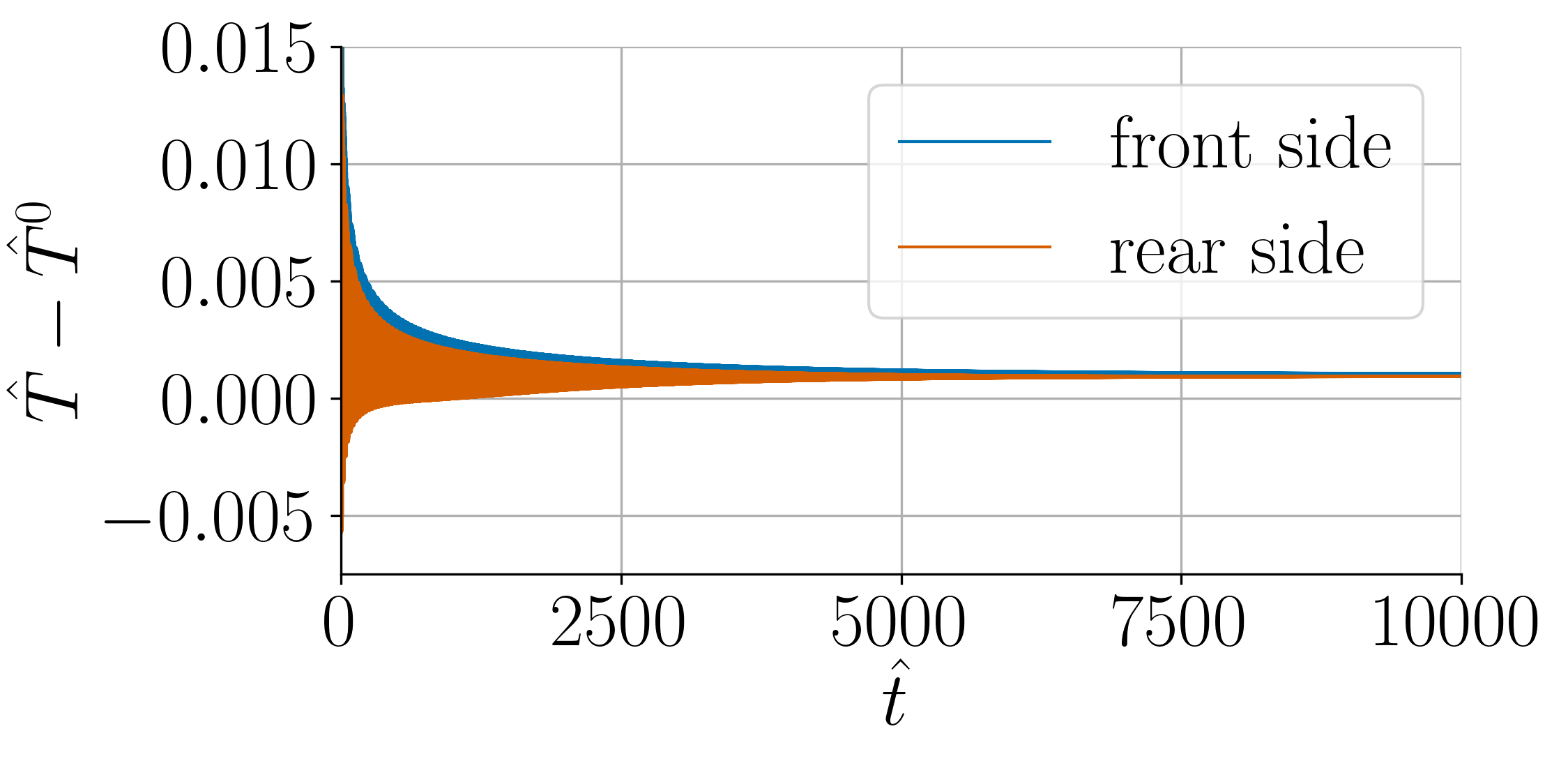}
    \end{minipage}%
    \caption{Comparing our time integration method (left) with a 4th-order Runge--Kutta method (right), for a shorter (top) and longer (bottom) timescale. In our simulation results, dispersion errors are absent, while the classical RK4 method introduces significant dispersion errors. (All settings for the two methods were identical.)}
    \label{fig:dispersion_time}
\end{figure*}

\begin{figure*}[!htb]
    \centering
    \begin{minipage}{0.4\textwidth}
        \includegraphics[width=\textwidth]{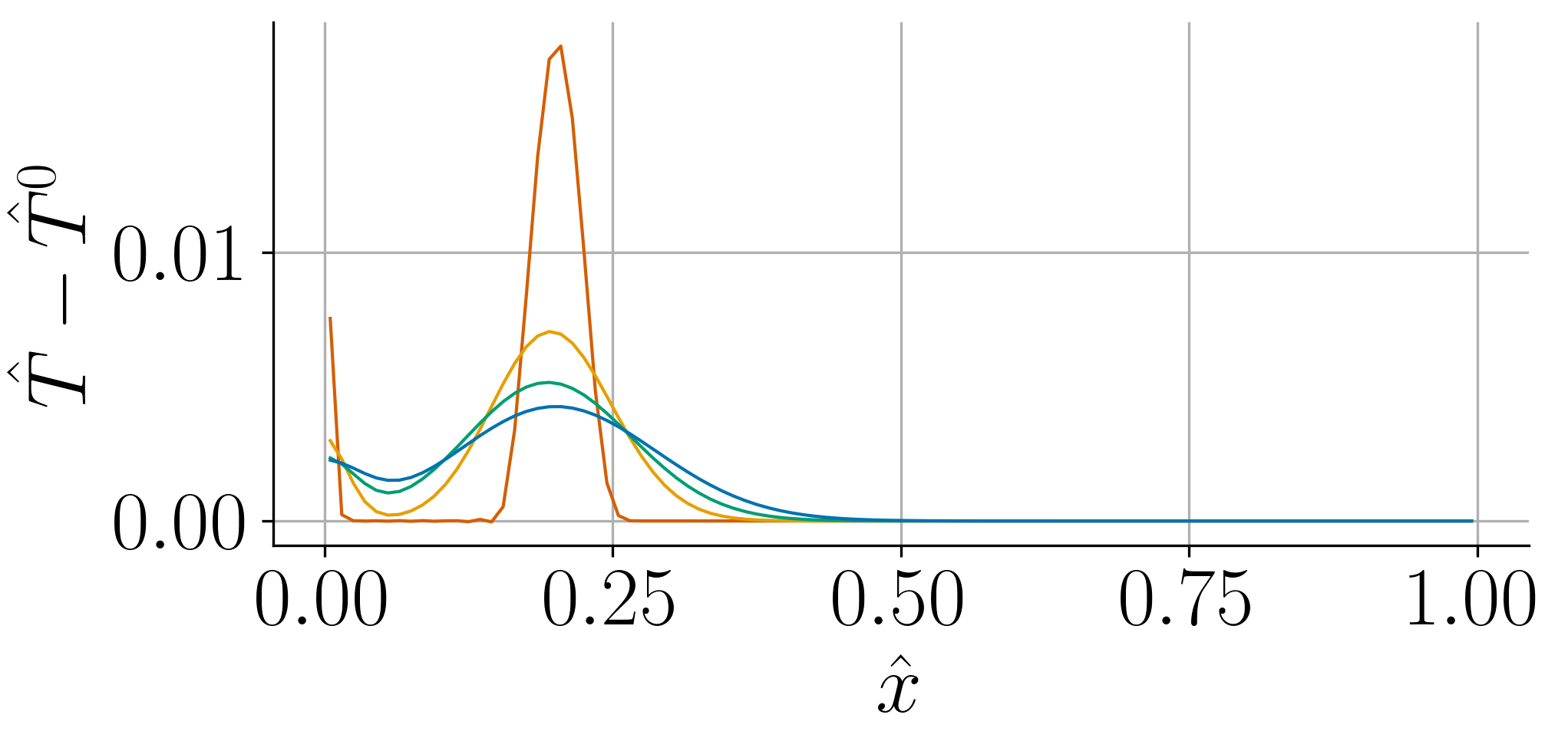}
    \end{minipage}%
    \begin{minipage}{0.4\textwidth}
        \includegraphics[width=\textwidth]{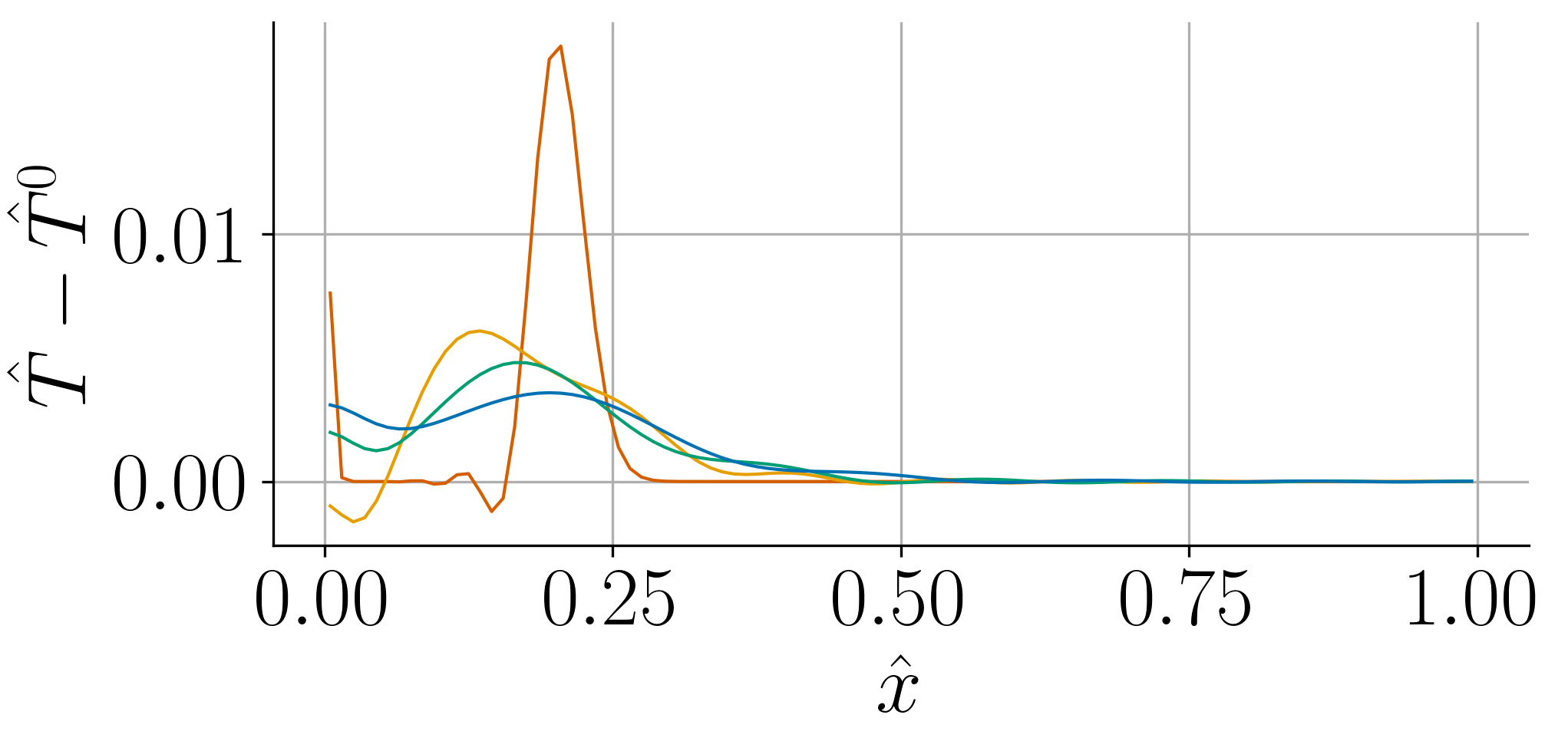}
    \end{minipage}%
    \caption{Comparing the effects of the dispersion error on the spatial distribution of the velocity fields, based on Fig.~\ref{fig:Tv-x}. Left: our time integration method, right: 4th-order Runge--Kutta method with identical settings. (Velocity fields at non-dimensional time instants 0.25 (red), 20.25 (orange), 40.25 (green), and 60.25 (blue).)}
    \label{fig:dispersion_space}
\end{figure*}

As a comparison, we have performed the same wave propagation simulation with an explicit 4th-order Runge--Kutta (RK4) time integration method, based on the semi-discrete equations \eqref{eq:semidiscrete}. All simulation settings, including the time step, have been identical to the previously shown simulation performed using the reversible--irreversible splitting method. To illustrate the advantages of this method, it is compared to the RK4 results in Fig.~\ref{fig:dispersion_time} and Fig.~\ref{fig:dispersion_space} for the temperature field (other fields show similar behaviour). In Fig.~\ref{fig:dispersion_time}, it can be clearly seen that the dispersion error introduced by the RK4 method affects the temporal solution significantly: the amplitude of the observed peaks is decreased, non-physical oscillation is introduced and as a result, the front and rear side temperatures deviate in the negative direction after each peak in an artificial way. Similarly, the dispersion error also distorts the numerical solution in the spatial sense, as shown in Fig.~\ref{fig:dispersion_space}. On the other hand, the structure-preserving composition of the reversible-irreversible splitting time integration method produces clean, sharp peaks and preserves the spatial distribution of the pulse over time. Is should be emphasized that this is true despite the fact that our method is only second-order, while the RK4 method is fourth order accurate: yet, it produces inferior results with the same settings, while also being more computationally intensive.

\subsubsection{Demonstration of second-order accuracy}
\begin{figure}[!htb]
    \centering
    \includegraphics[width=0.25\textwidth]{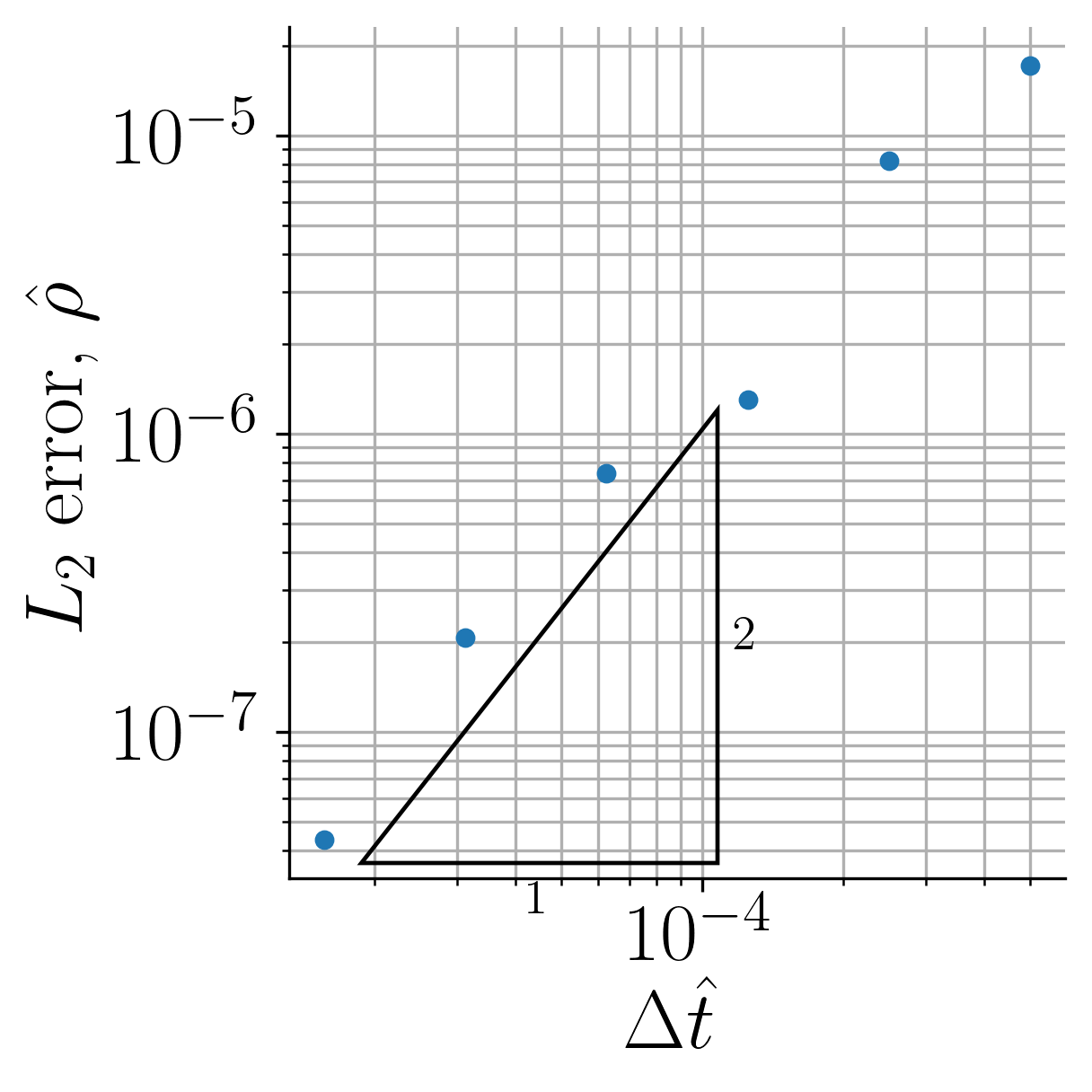}
    \includegraphics[width=0.25\textwidth]{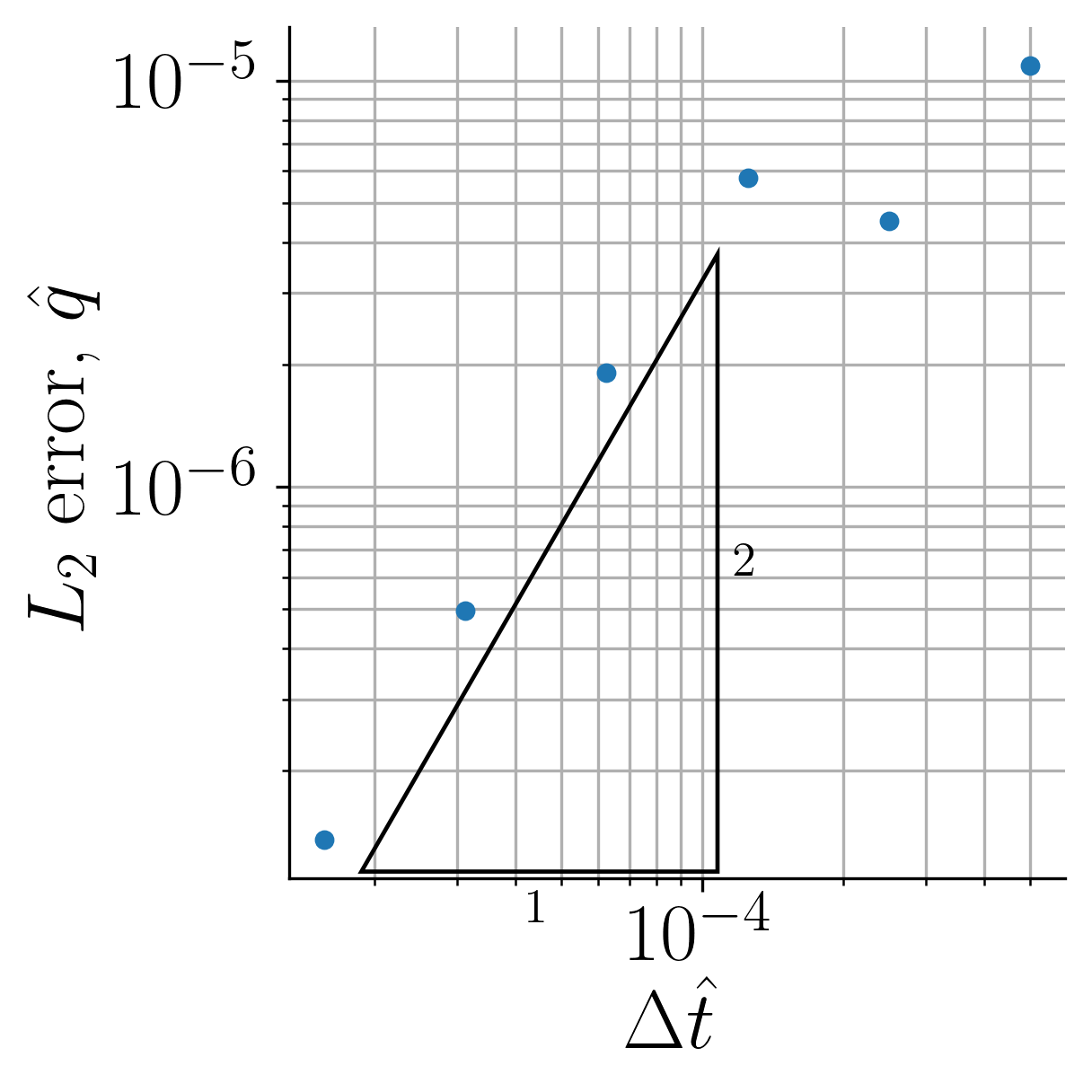}
    \caption{Second-order convergence of the numerical solution of the fields $\hrho$ and $\hq$.}
    \label{fig:numerical-convergence}
\end{figure}

In Sec.~\ref{sec:acc}, we have analytically proven that the developed scheme is second-order accurate. This result can also be verified numerically. By setting $\hDx/\hDt = \mathrm{const.}$, convergence graphs for all fields can be obtained by decreasing the time step and calculating the $L_2$ error of the solutions compared to a much finer reference solution. Fig.~\ref{fig:numerical-convergence} shows the second-order convergence of the fields $\hrho$ and $\hq$, for reference; the other fields are omitted for brevity but converge similarly.

\subsubsection{Grid independence test}
\begin{figure}[!htb]\centering
    \includegraphics[width=0.4\textwidth]{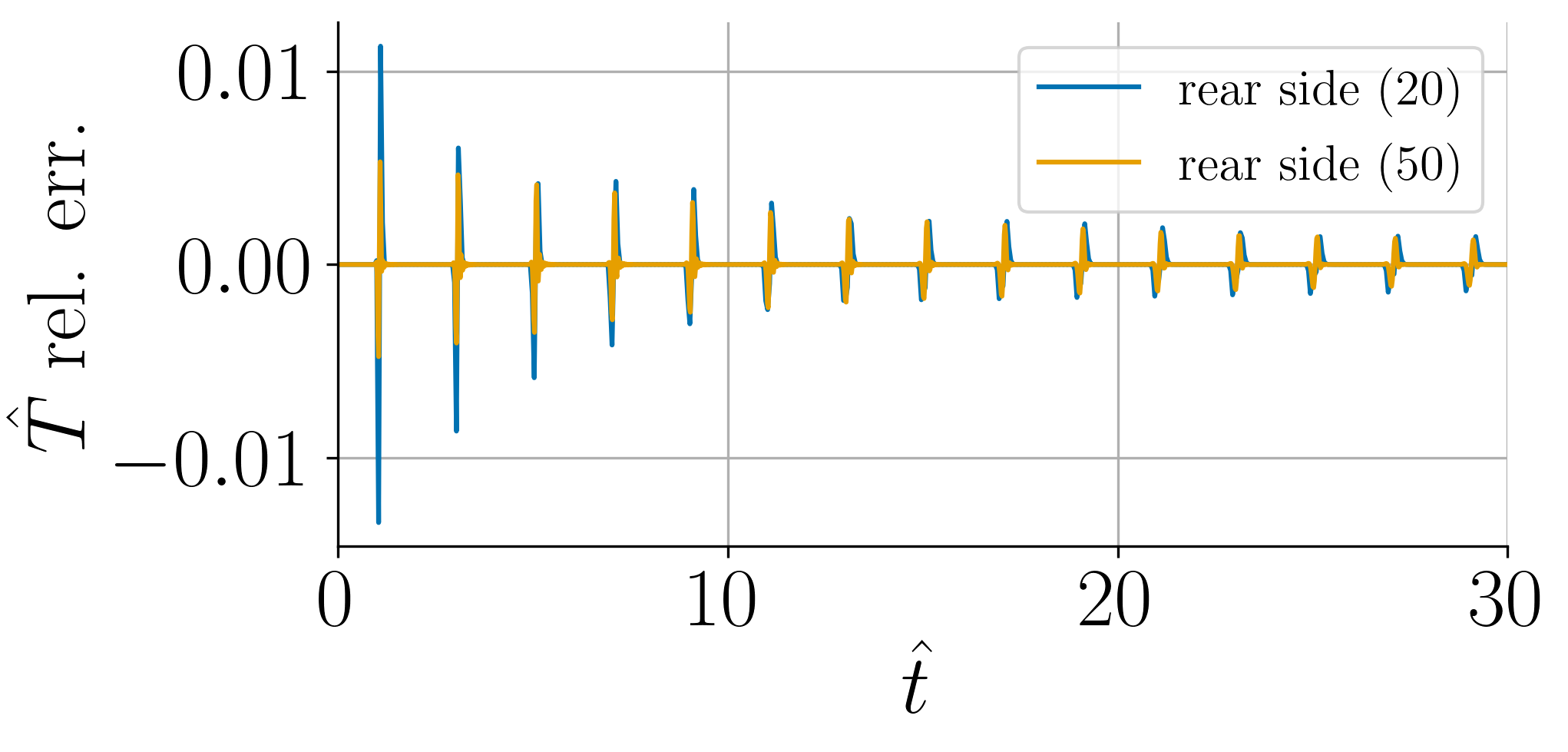}\\
    \includegraphics[width=0.4\textwidth]{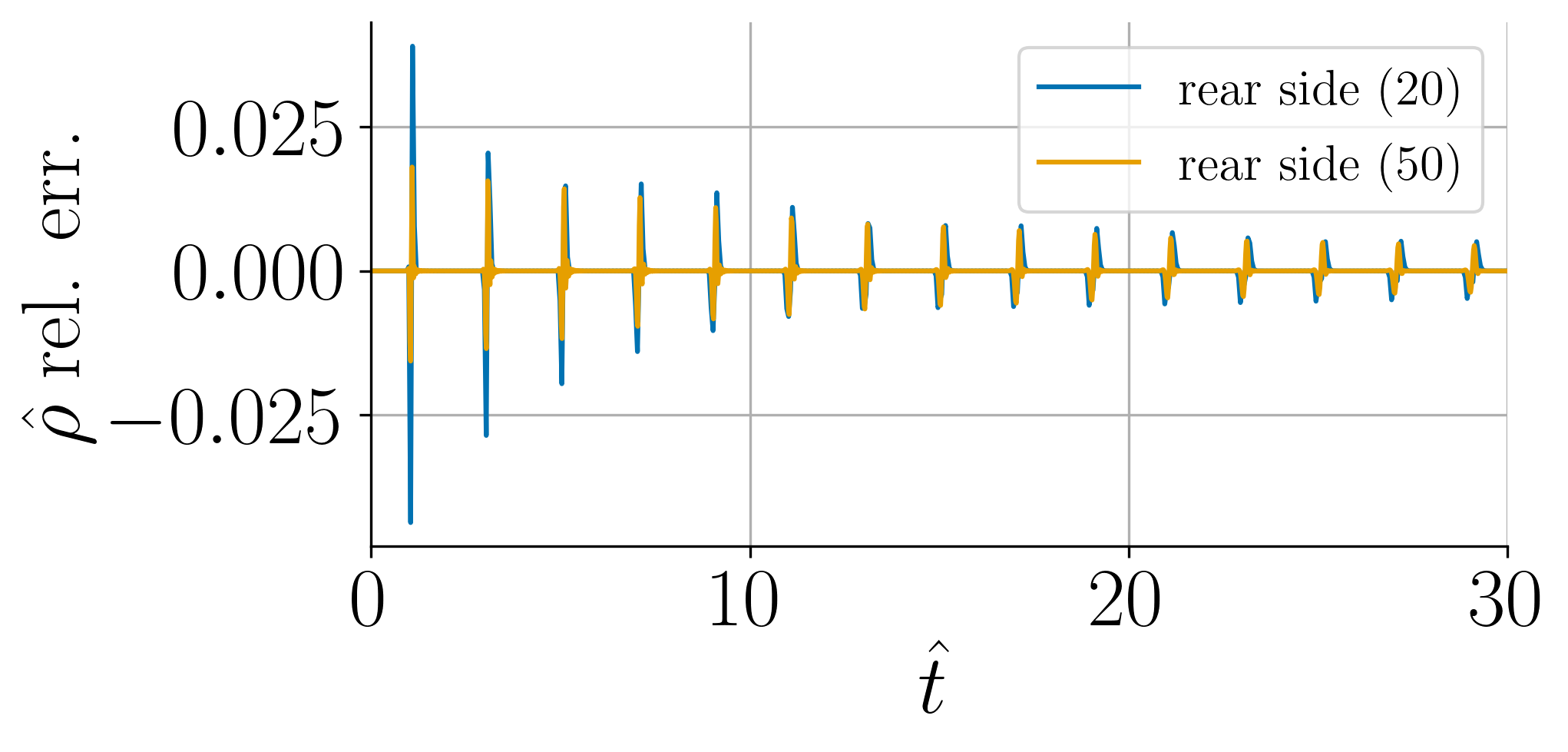}
    \caption{Grid dependence of time evolution of temperature and density at the rear side of the sample.}
    \label{fig:grid-ind-wave}
\end{figure}

We have also performed a grid independence test. The relative deviation of solutions produced by applying 20 and 50 cells are compared to the solution produced by 100 cells, which are presented on the rear side of the sample in Fig.~\ref{fig:grid-ind-wave}. Decreasing the number of applied cells from 100 to 20 causes less than three \% deviation for all fields.

\begin{figure}[!htb]\centering
    \includegraphics[width=0.4\textwidth]{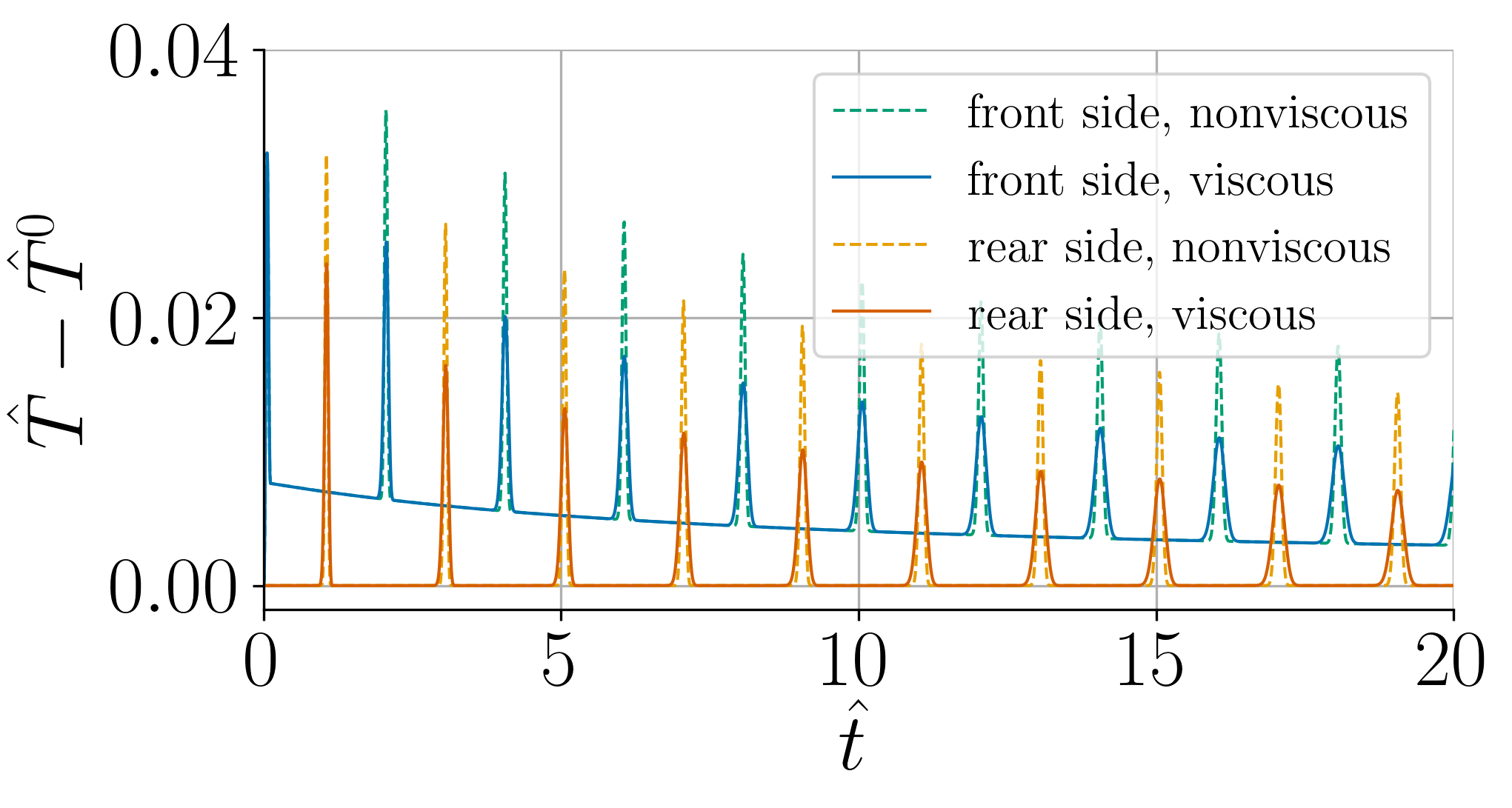}
    \caption{The effect of viscous momentum transport on the temperature field.}
    \label{fig:Trp-t-vis}
\end{figure}

Finally, we present the effect of viscosity, which in the investigated state results in $ \qRea = \qPea / \qPr = 17226.957 $. Viscous momentum transport leads to further dissipation in addition to heat conduction, therefore, resulting in stronger attenuation of the amplitudes and in faster and more intense dispersive spreading of the waveform, presented in Fig.~\ref{fig:Trp-t-vis}.

\subsection{The thermoacoustic approximation of the piston effect} \label{subsec:piston}

This section is devoted to analyze the piston effect in supercritical fluids via the developed scheme. In this case the duration of excitation is much longer than the acoustic time scale. Therefore, wave propagation occurs only at the beginning of the process until the initial pressure disturbance reaches the other end of the pipe, correspondingly, no reflection takes place. In order to prove the intensive effect of thermal expansion near the critical point, we present the same phenomenon with same dimensional excitation parameters for an ideal gas state, too.

Based on \cite{straub1995dynamic}, we have defined a cross-section specific heat of $ q = 30 \ {\rm \frac{J}{m^2}} $, to which in the above defined supercritical state the dimensionless value $ \hat{q} = 0.00035 $ corresponds, the non-dimensional pulse duration $ \hht_{\rm P} = 50 $ and P\'eclet number $ \qPea = 10^7 $ are chosen, which latter implies $ \qRea = 1722695.711 $. Now, viscous momentum transport is also taken into account. Since wave propagation is suppressed, we have reached a stable solution with the choice of $ \mathpzc{Co} = 1 $, which, with the same settings, would result in an unstable solution for the damped wave propagation investigated previously.

\begin{figure}[!htb]\centering
    \includegraphics[width=0.4\textwidth]{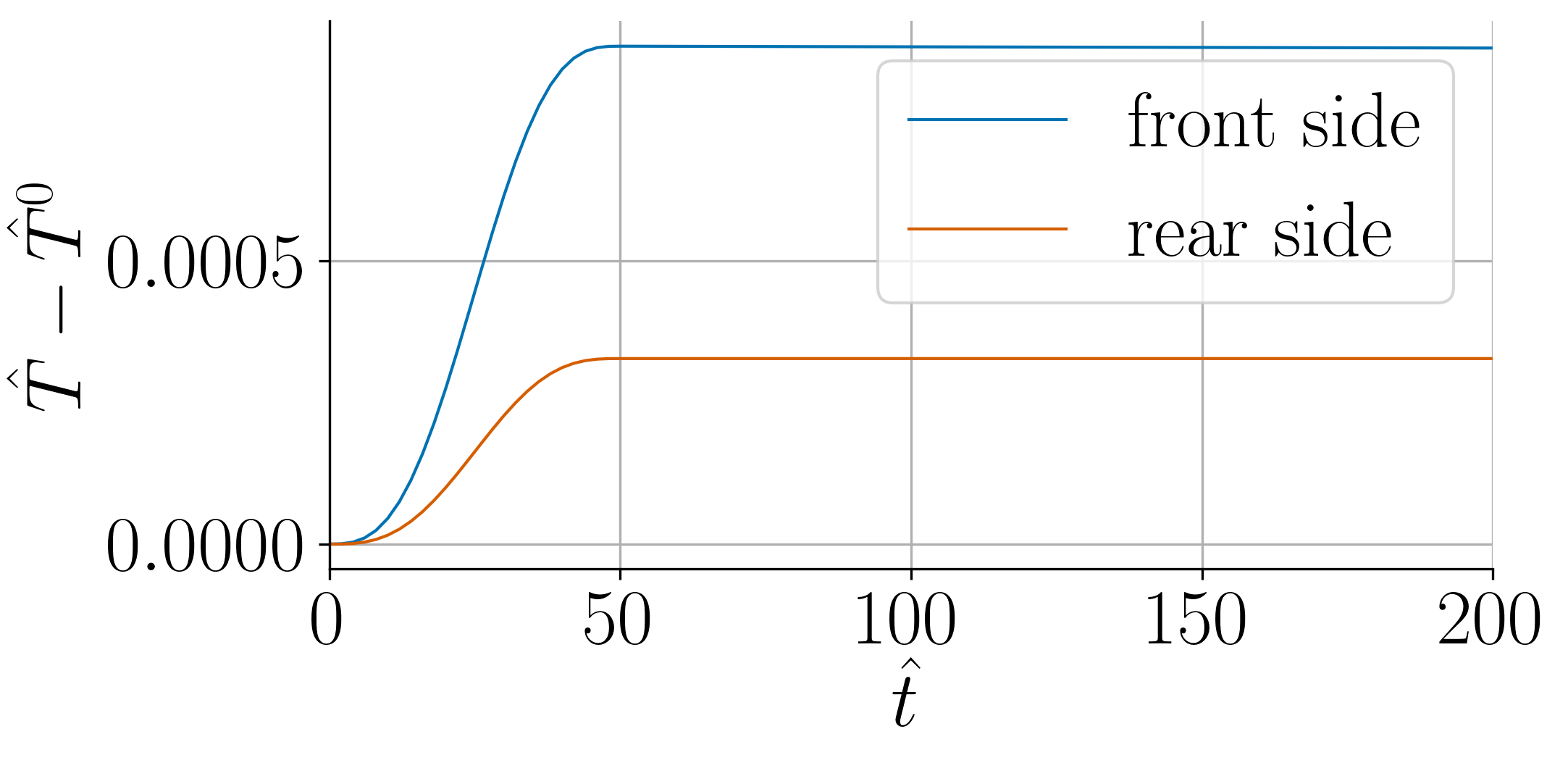}\\
    \hskip -1.3ex \includegraphics[width=0.4\textwidth]{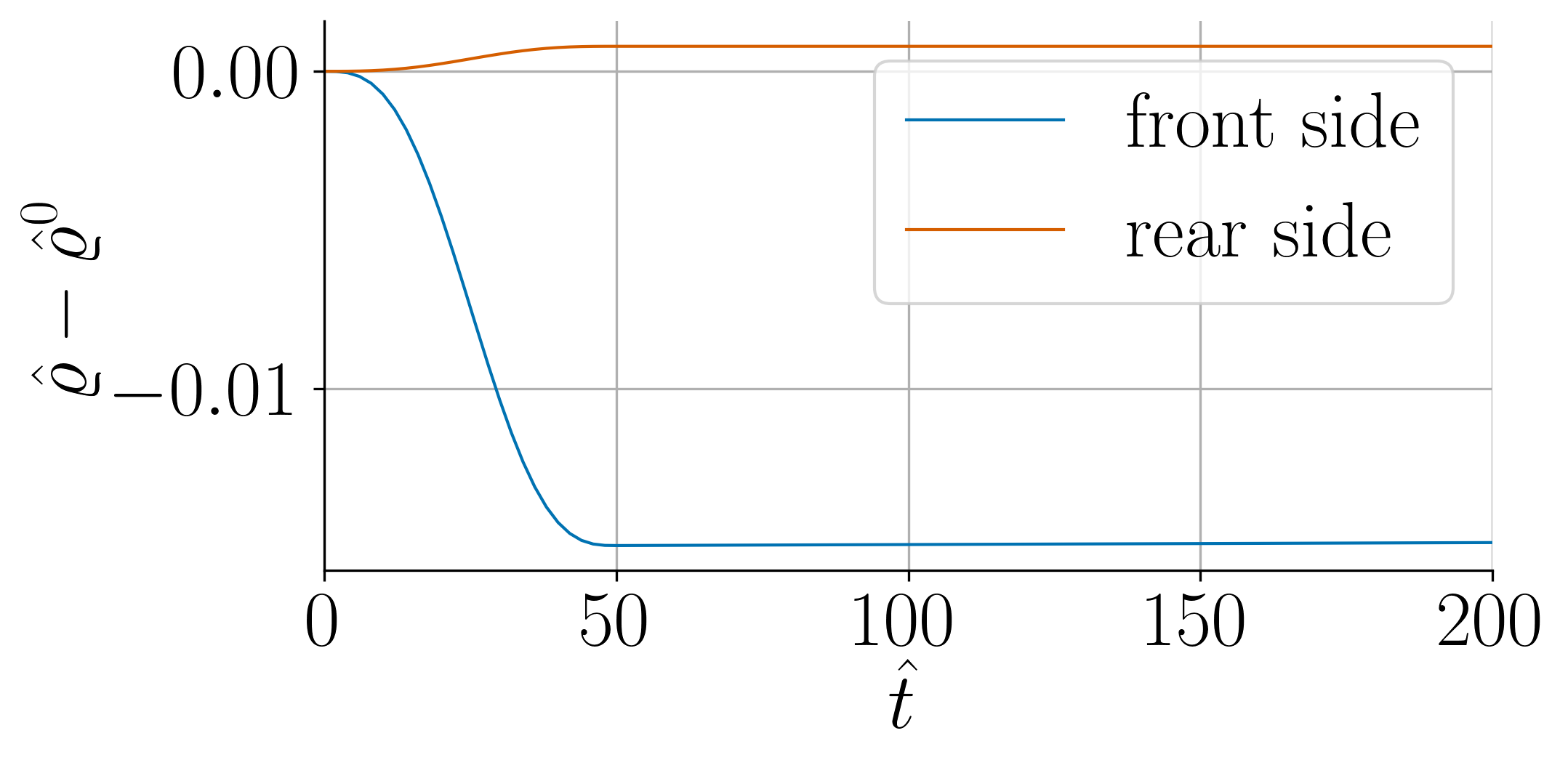}\\
    \includegraphics[width=0.4\textwidth]{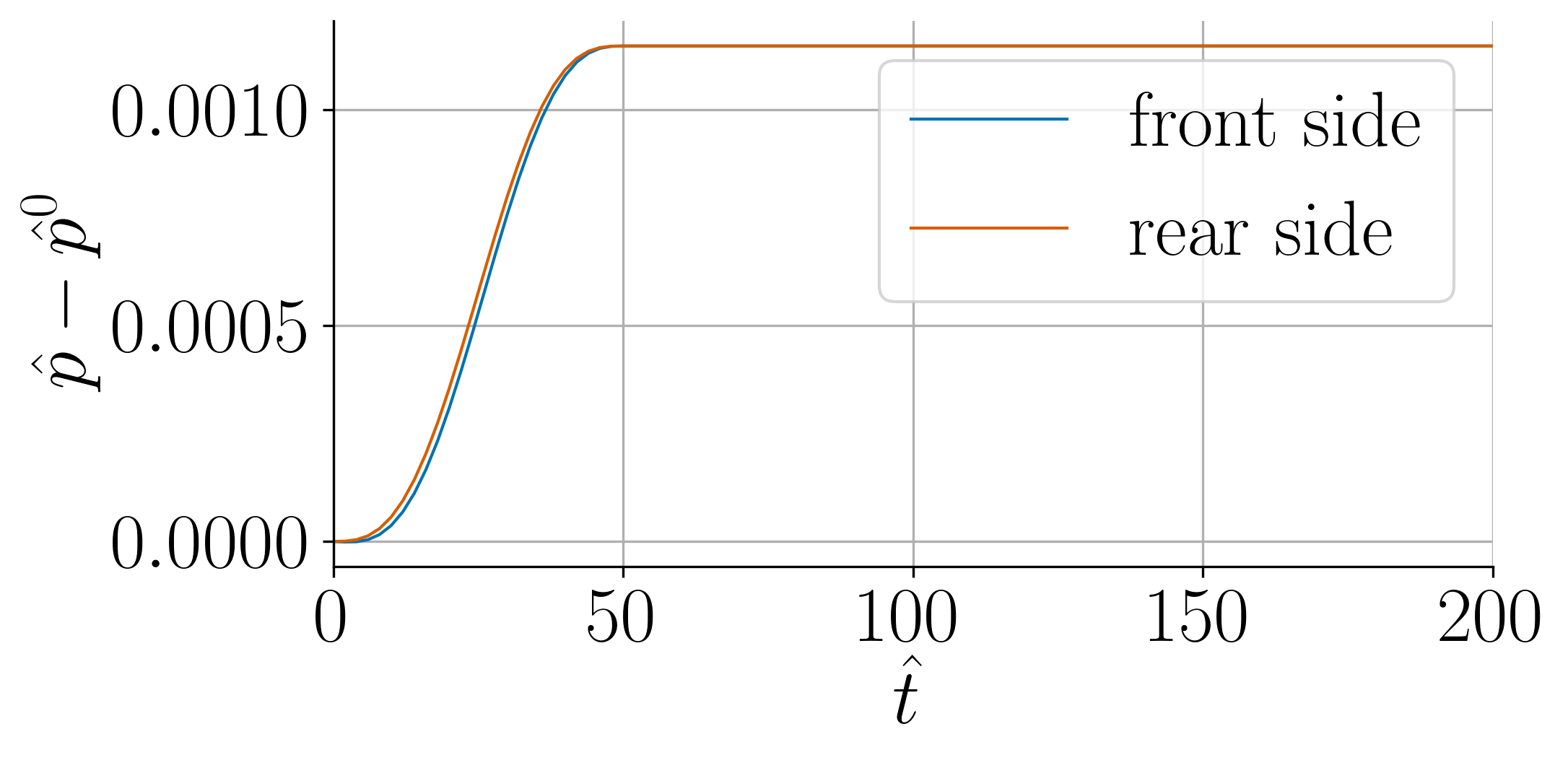}\\
    \caption{Time evolution of temperature, density, and pressure at the front side and at the rear side of the sample in the piston effect in supercritical fluid state.}
    \label{fig:sc-piston}
\end{figure}

\begin{figure}[!htb]\centering
    \includegraphics[width=0.4\textwidth]{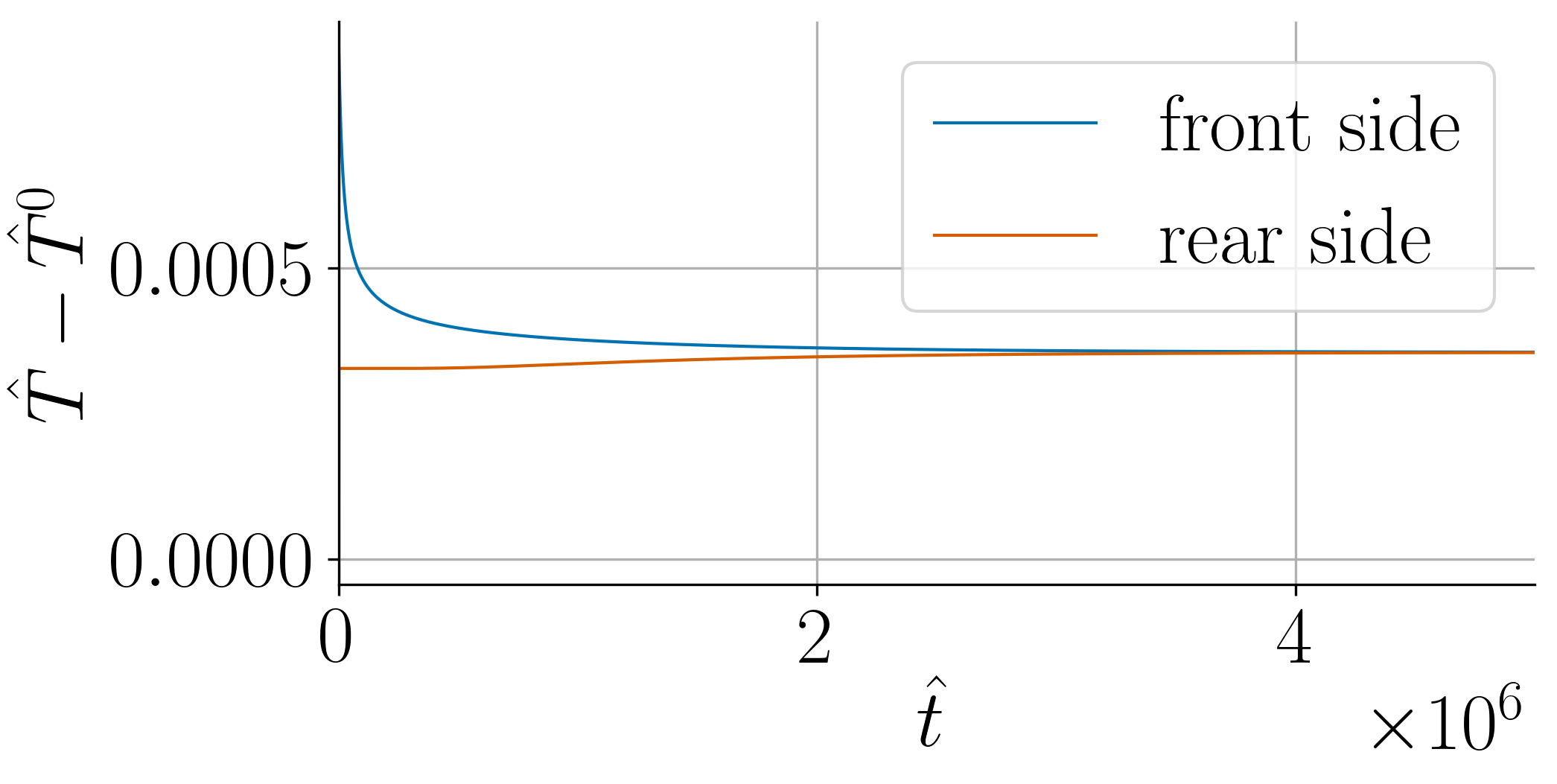}\\
    \hskip -1.3ex \includegraphics[width=0.4\textwidth]{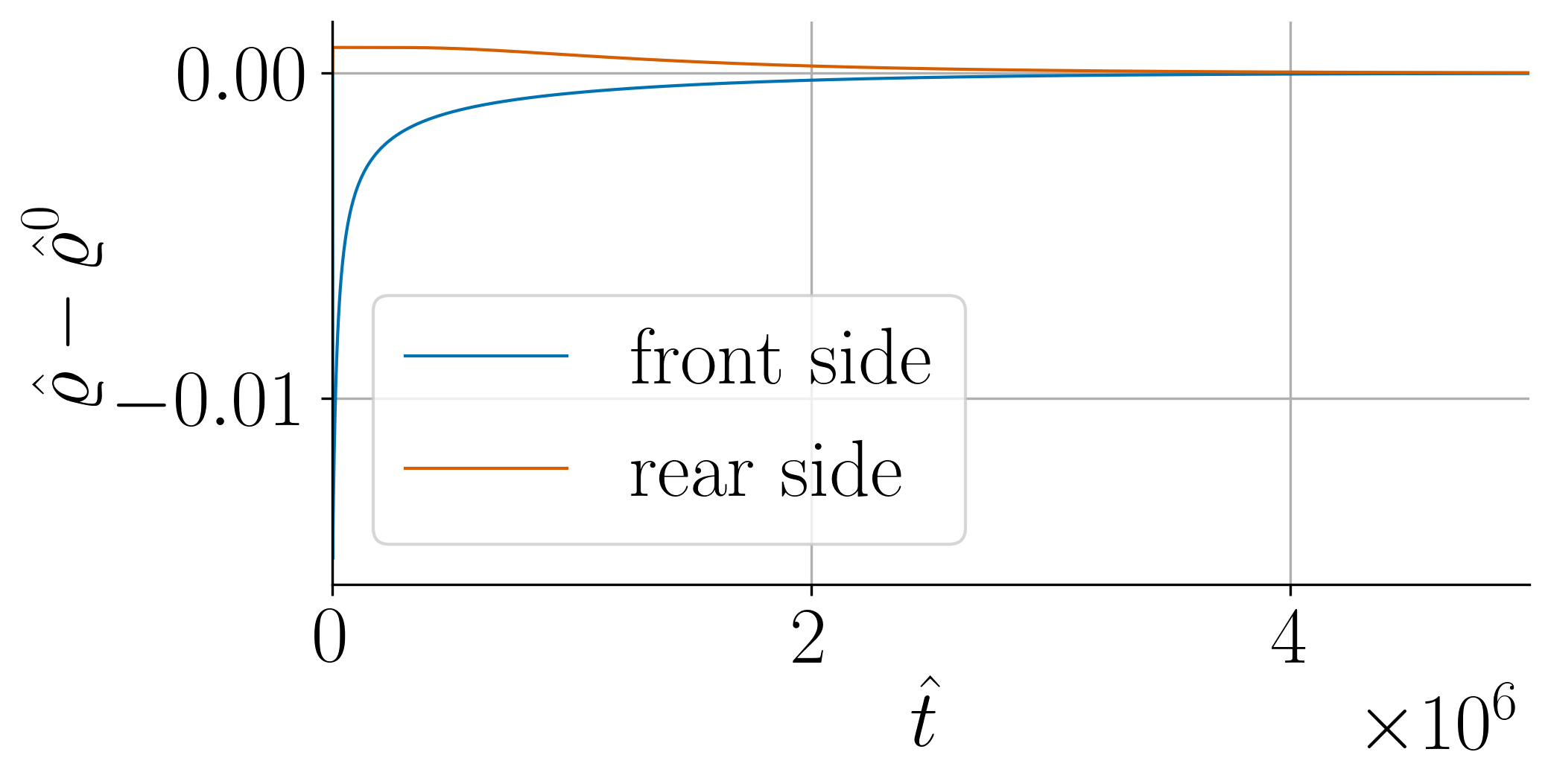}\\
    \caption{Time evolution of temperature, density, and pressure at the front side and at the rear side of the sample in the piston effect in supercritical fluid state.}
    \label{fig:sc-piston-long}
\end{figure}

Essentially, the pressure is homogeneous during the entire process and changes only during the heat pulse. During this time period, relative to the duration of the entire homogenization process, temperature and density at the rear side of the sample suddenly increase. In fact, most of the temperature change takes place during the heat pulse, temperature change during the diffusion is almost negligible. 
{At first glance, one might think that the initial rapid changes go beyond the limits of the linear framework but, compared to the acoustic time scale (\ie the time required to measure a signal on one end of the sample that is excited on the other side), these changes are not that steep.} The sudden changes at the beginning of the process are presented in Fig.~\ref{fig:sc-piston}, while long-time behavior, \ie homogenization of temperature and density, are shown in Fig.~\ref{fig:sc-piston-long}.

To enlighten the difference between heat conduction in supercritical and in an ideal gas state, we realized the previous calculation by considering the same $ \bT = 305 $ K temperature but with modified pressure chosen to be $ p_0 = 0.1 $ MPa, the corresponding density is $ \brho = 1.744 \ {\rm \frac{kg}{m^3}} $. In this state, carbon dioxide can be considered an ideal gas, since the ideal gas deviation factor is $ Z = p_0 / \left( R \brho \bT \right) = 0.995 \approx 1 $ with the specific gas constant $ R = 188.964 \ {\rm \frac{J}{kg \ K}} $ of CO$_2$. In this state, the required dimensionless parameters are $ \bgam = 1.291 $, $ \qB = 1.015 $ (which also indicates that the investigated state is nearly an ideal gas state), $ \qEca = 0.283 $, $ \qPr = 0.762 $ and $\bReta \approx 0.4 $ \cite{wang2019bulk}. Furthermore, keeping dimensional quantities constants imply $ \hat{q} = 0.125 $, $ \hht_{\rm P} = 73.694 $, $ \qPea = 16134.536 $ and $ \qRea = 21185.723 $.

\begin{figure}[!htb]\centering
    \includegraphics[width=0.4\textwidth]{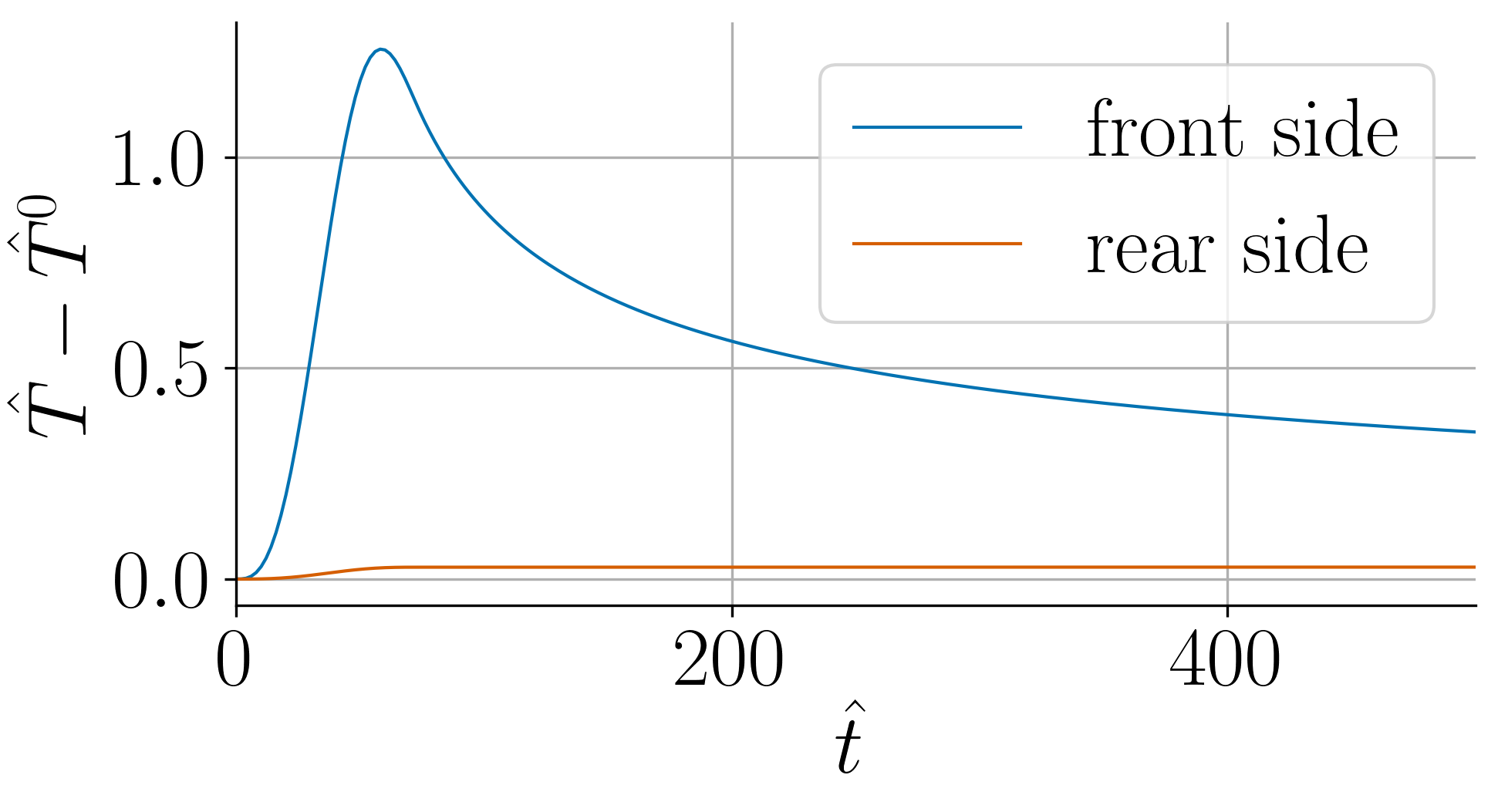}\\
    \hskip -1.3ex \includegraphics[width=0.4\textwidth]{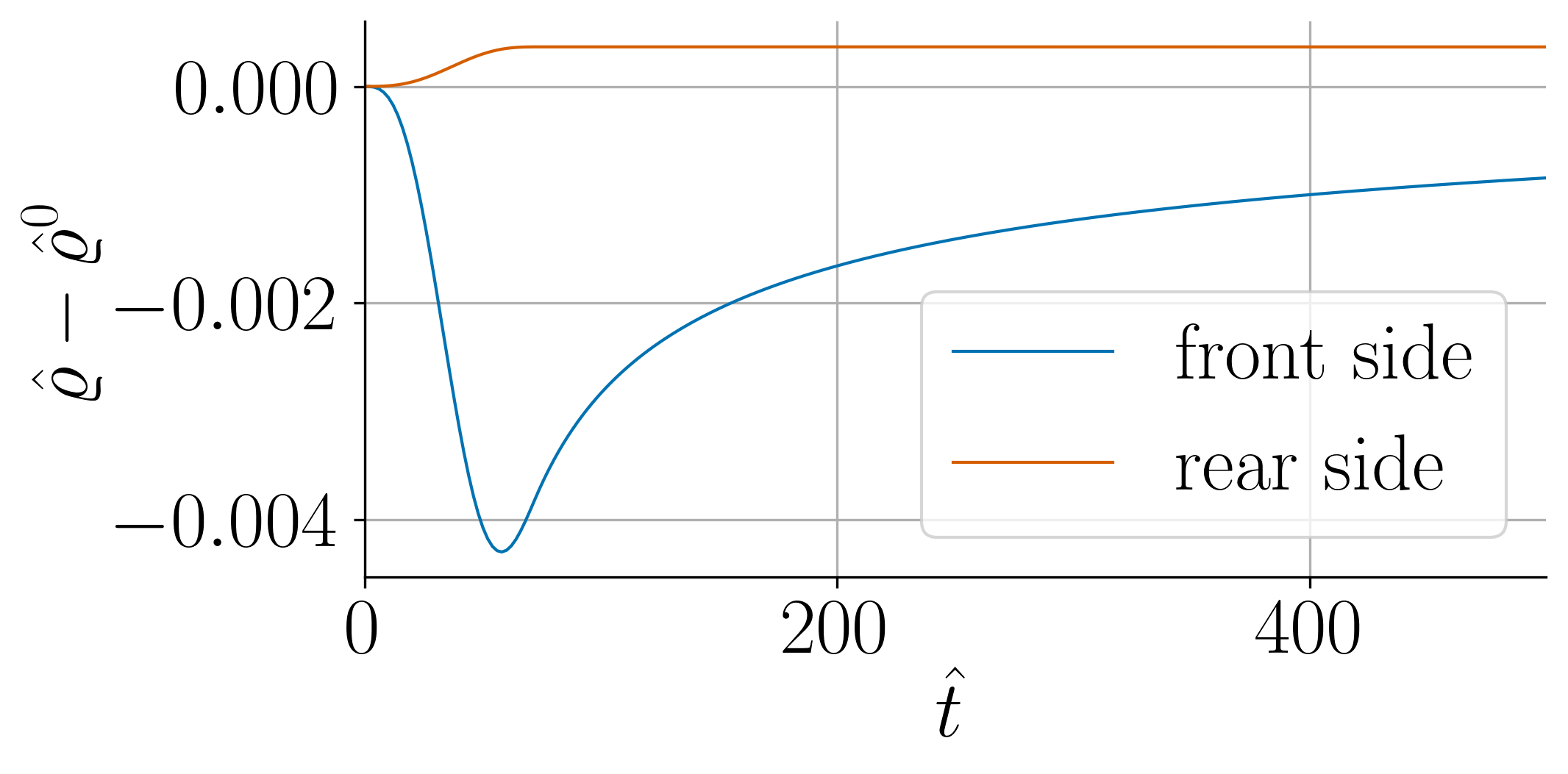}\\
    \includegraphics[width=0.4\textwidth]{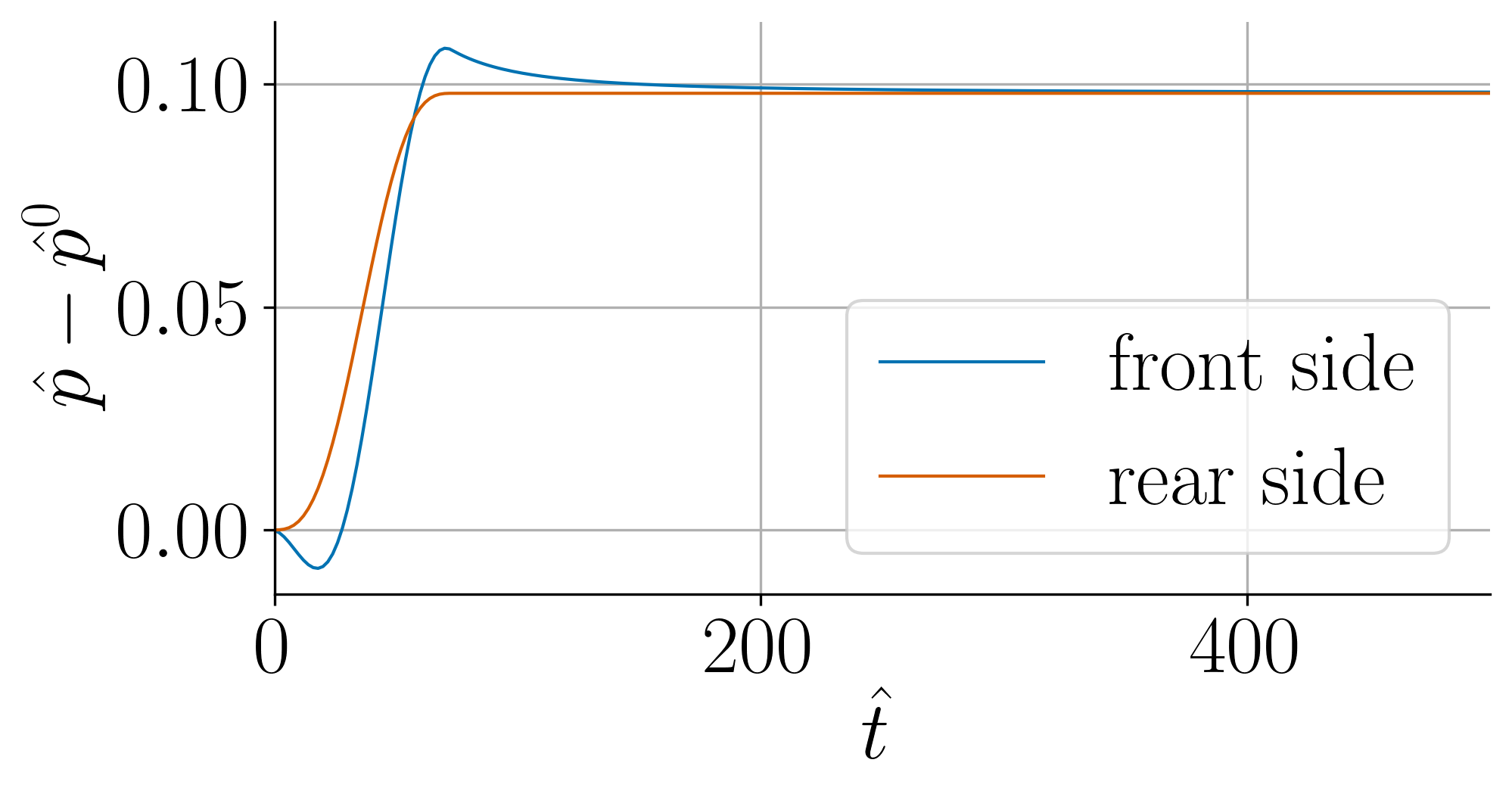}\\
    \caption{Time evolution of temperature, density, and pressure at the front side and at the rear side of the sample in the piston effect in an ideal gas state.}
    \label{fig:id-piston}
\end{figure}

\begin{figure}[!htb]\centering
    \includegraphics[width=0.4\textwidth]{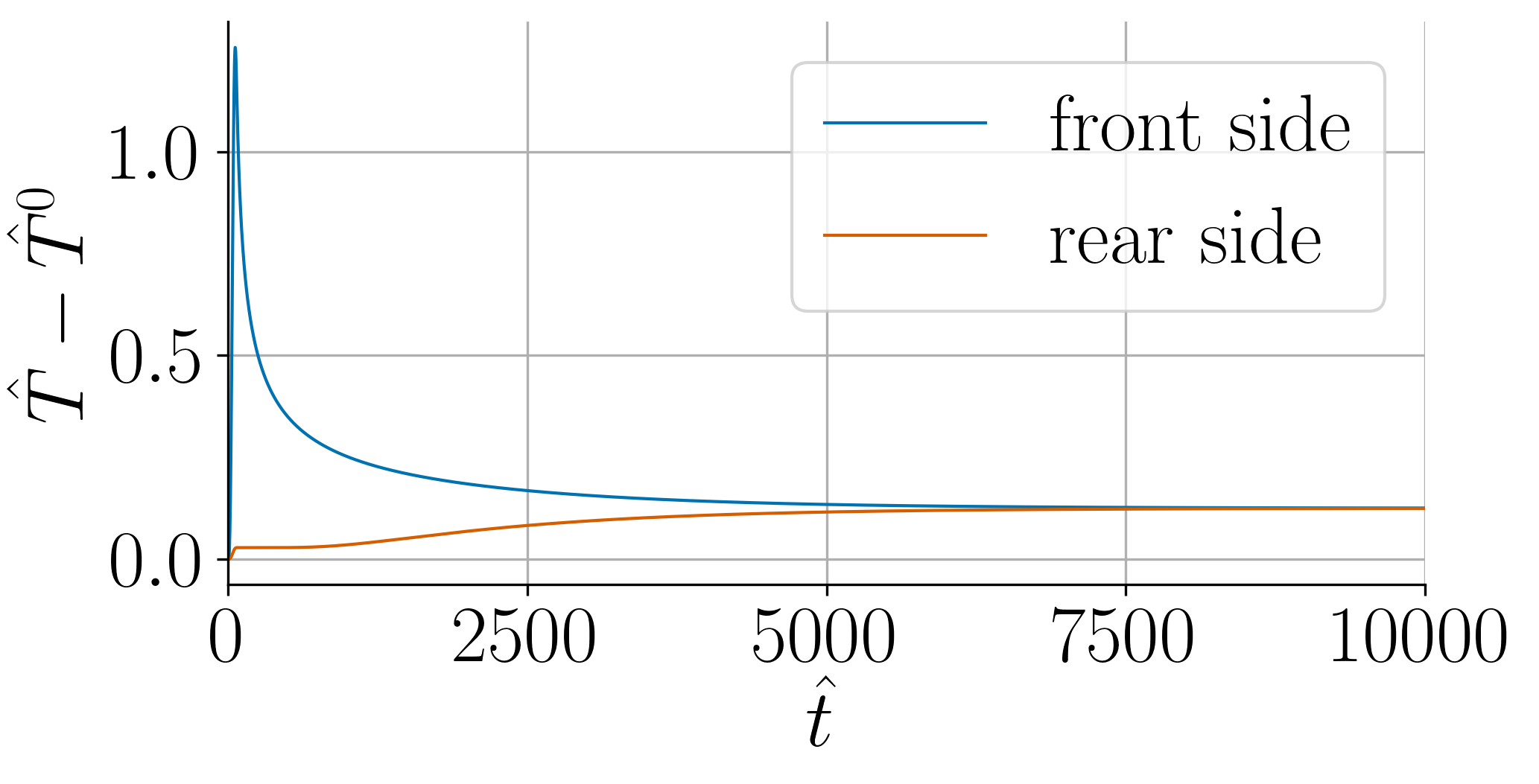}\\
    \hskip -1.3ex \includegraphics[width=0.4\textwidth]{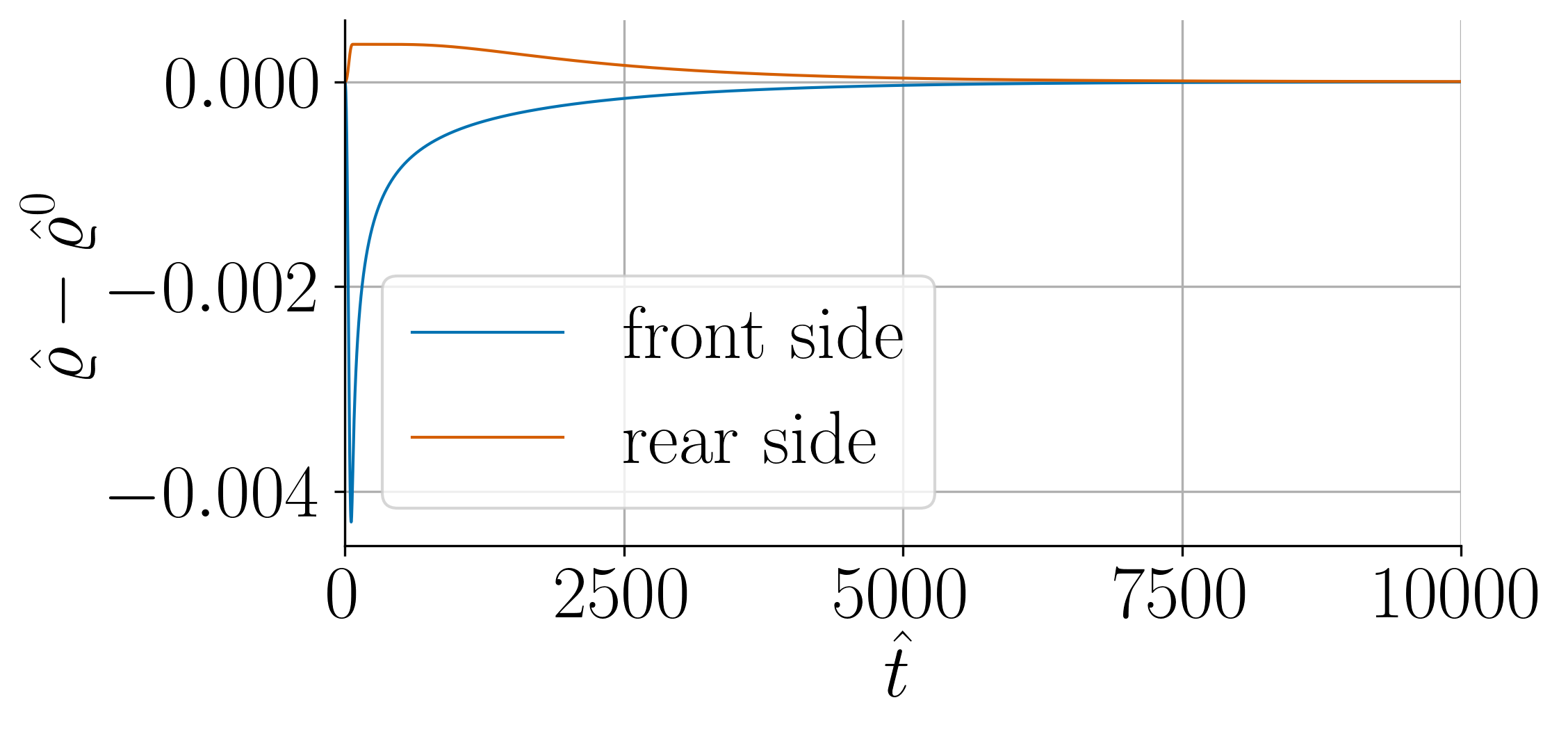}\\
    \caption{Time evolution of temperature, density, and pressure at the front side and at the rear side of the sample in the piston effect in an ideal gas state.}
    \label{fig:id-piston-long}
\end{figure}

The obtained short-time behavior for temperature, density, and pressure are shown in Fig.~\ref{fig:id-piston}, while time evolution along the entire process is presented in Fig.~\ref{fig:id-piston-long}. In an ideal-gas state, at the front side of the sample, an initial pressure decrease can be observed, and homogenization of pressure occurs delayed compared to the end of the heat pulse. The sudden increase of rear side temperature and density is also visible, however, due to the smaller value of the thermal expansion coefficient, the initial temperature increase is not as large as in supercritical state. Most of the temperature change takes place during the diffusion-dominated part of the process.

\begin{figure}[!htb]\centering
    \includegraphics[width=0.4\textwidth]{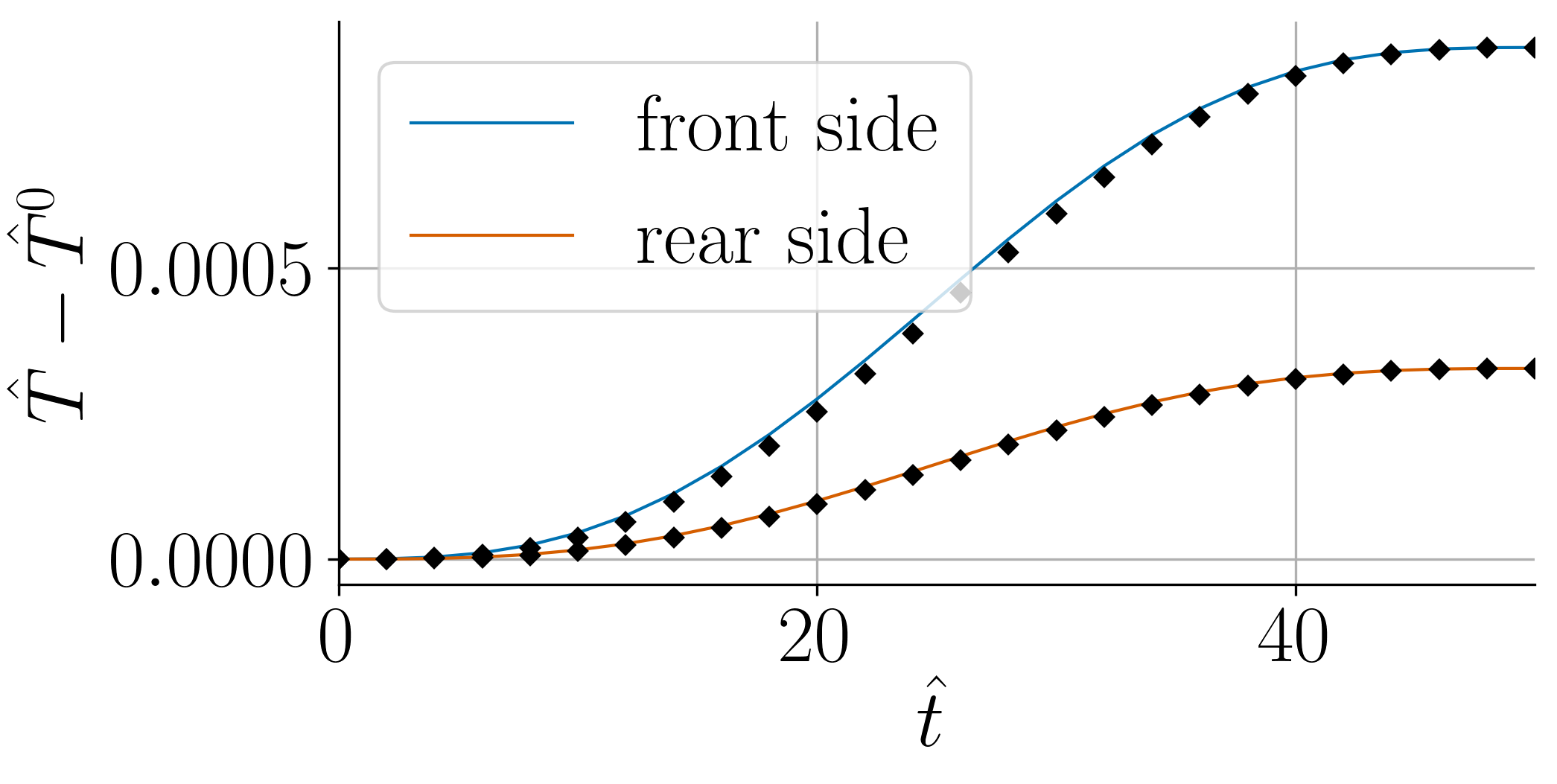}\\
    \hskip -1.3ex \includegraphics[width=0.4\textwidth]{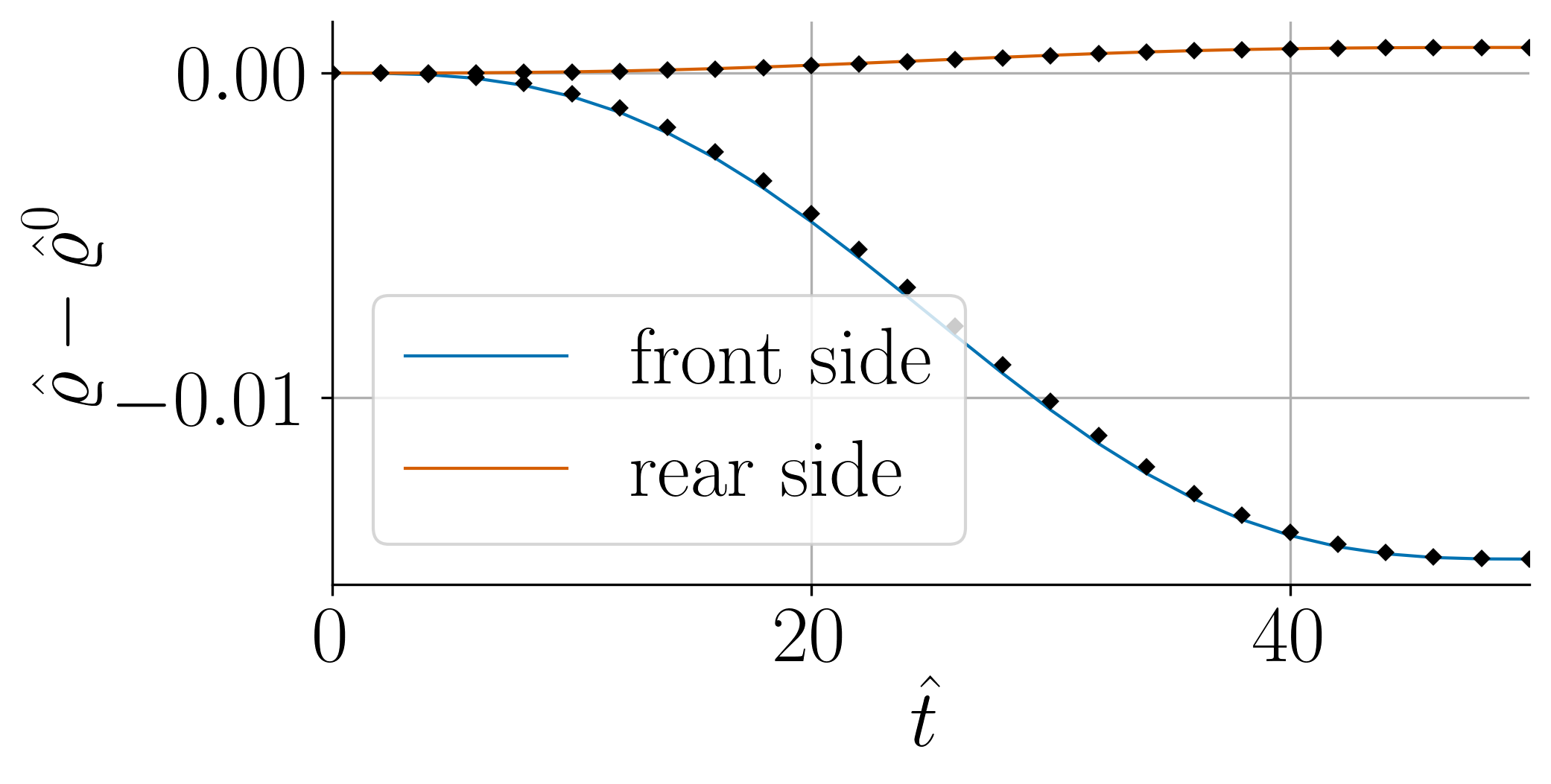}\\
    \includegraphics[width=0.4\textwidth]{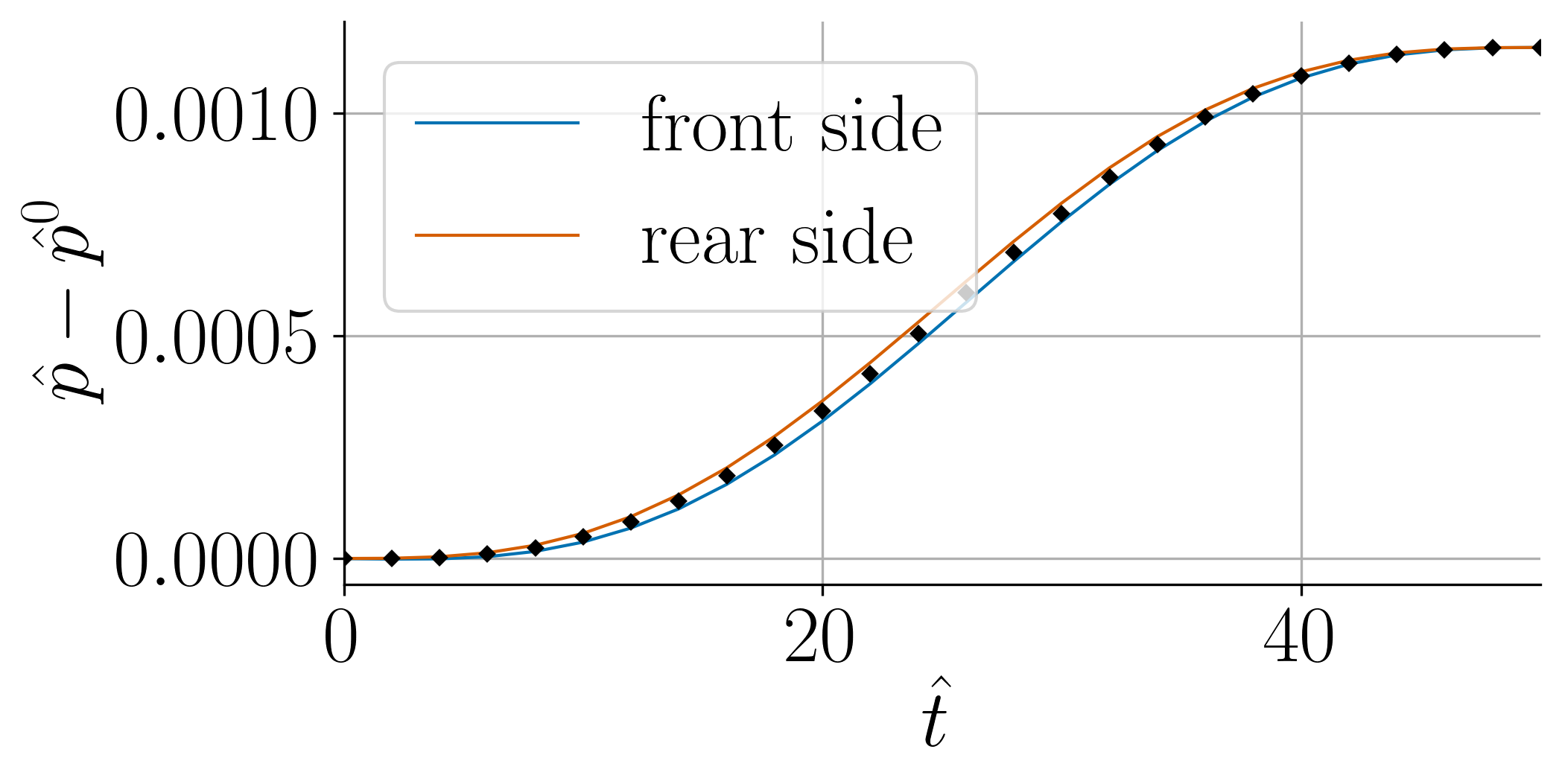}\\
    \caption{
    {Comparison of the thermoacoustic (continuous lines) and Boukari's approximation (rhomboids) of the piston effect in supercritical fluid state.}}
    \label{fig:bou-vs-ac-sc}
\end{figure}

\begin{figure}[!htb]\centering
    \includegraphics[width=0.4\textwidth]{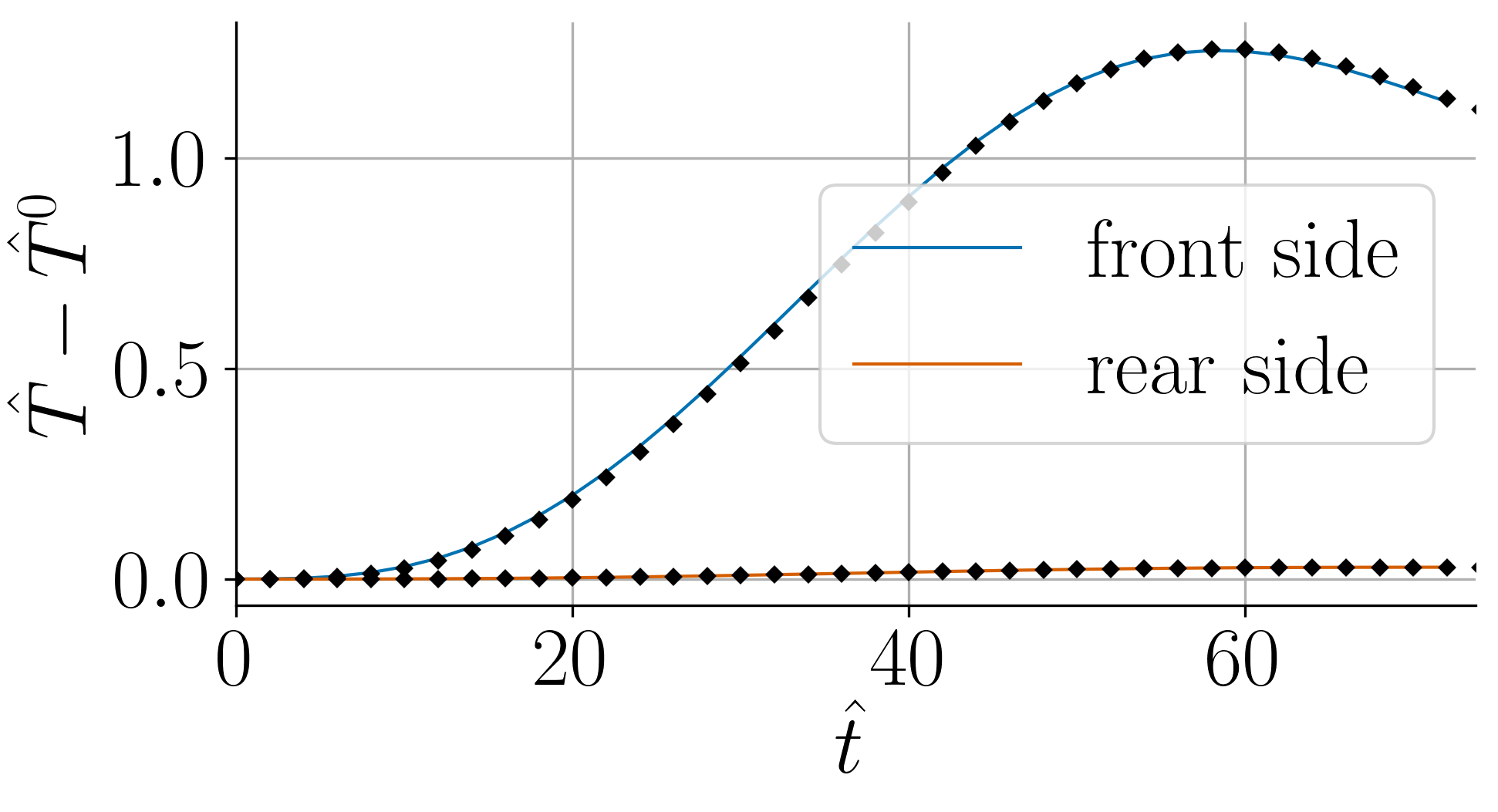}\\
    \hskip -1.3ex \includegraphics[width=0.4\textwidth]{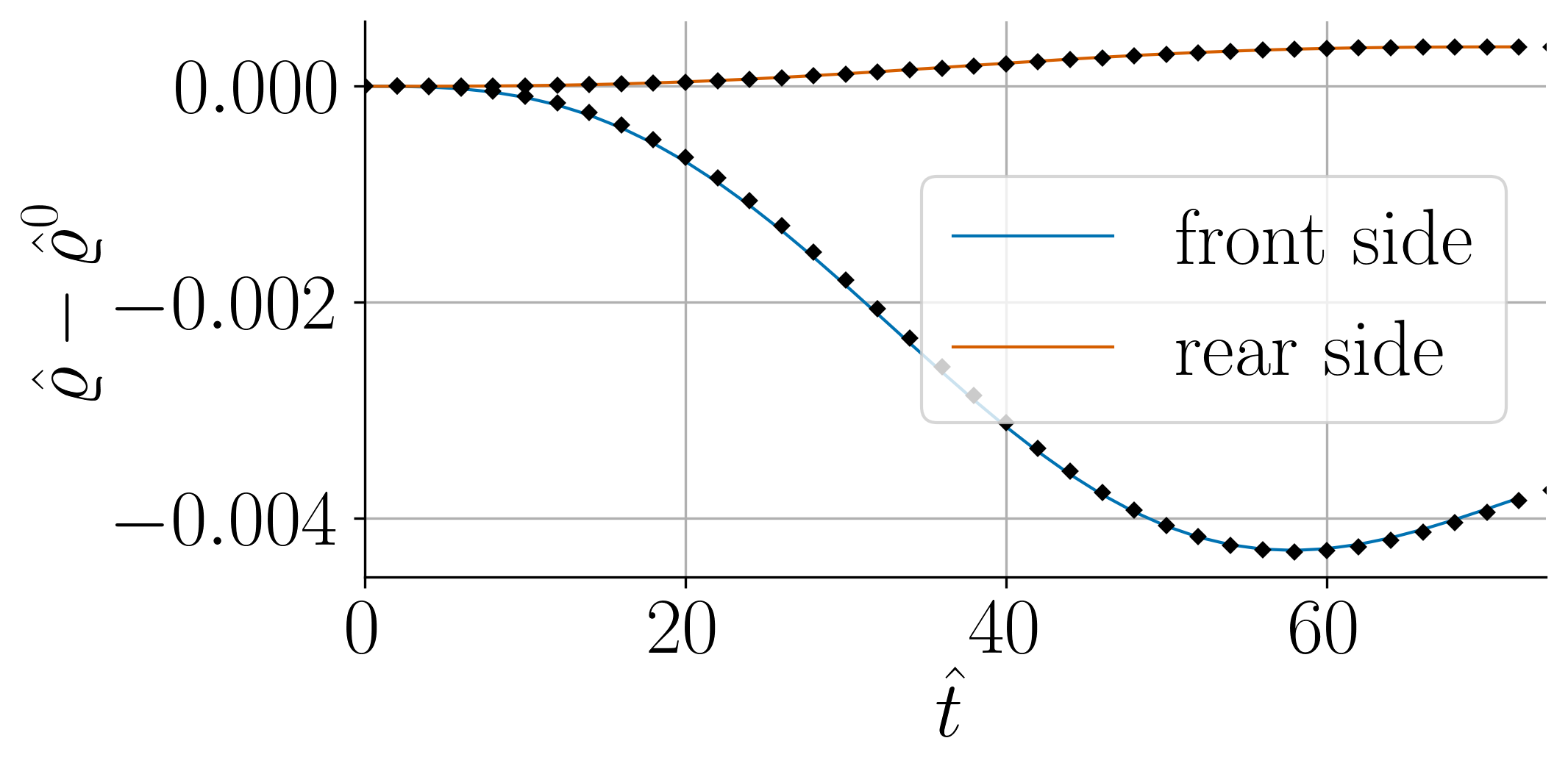}\\
    \includegraphics[width=0.4\textwidth]{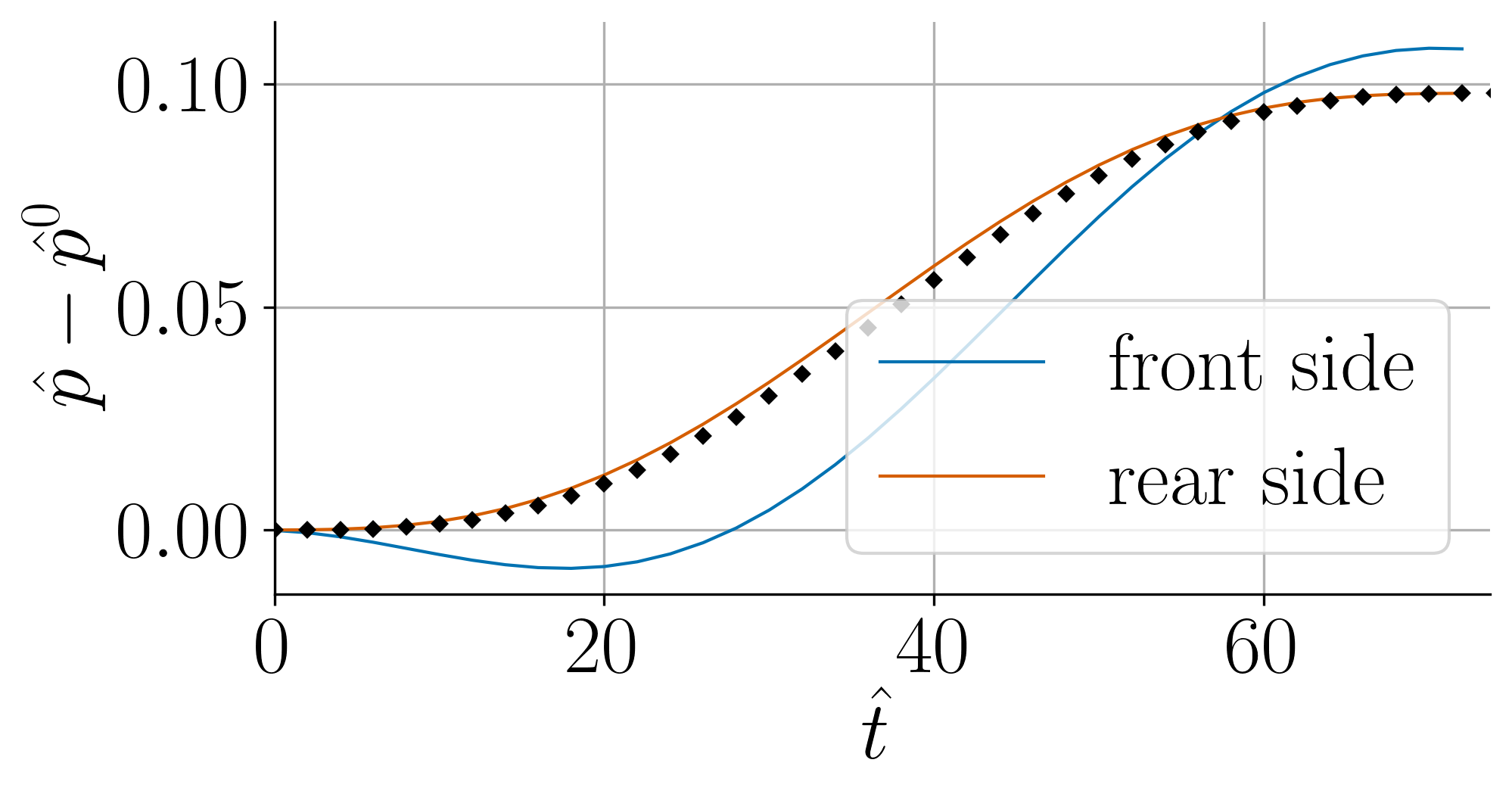}\\
    \caption{
    {Comparison of the thermoacoustic (continuous lines) and Boukari's approximation (rhomboids) of the piston effect in nearly ideal gas state.}}
    \label{fig:bou-vs-ac-ig}
\end{figure}

{As a known benchmark, we have compared the results obtained with our developed scheme to the results of Boukari's thermal equilibration model with homogeneous pressure \cite{boukari1990critical}, which effectively describes thermal equilibration of strongly compressible pure fluids near the liquid-vapor critical point. This approximation neglects velocity field and gravitational acceleration, accordingly, pressure is homogeneous during the entire process. Since for supercritical fluids intensive pressure changes (appearing in the beginning of the process) play an important role in equilibration of temperature, it is implemented in the heat equation as a homogeneous but time-dependent source term. A brief review of Boukari's approximation is given in Appendix~\ref{sec:app-B}. The results obtained with the two  examined approaches differ only slightly from each other, mainly in the beginning part of the process, therefore, the comparisons are presented only during the duration of the heat pulse. Fig.~\ref{fig:bou-vs-ac-sc} displays the results obtained for supercritical fluid state, while Fig.~\ref{fig:bou-vs-ac-ig} for nearly ideal gas state. One can conclude that, in general, the two approximations give nearly the same results. The most spectacular difference is presented in the pressure field, namely, in the supercritical state, Boukari's approximation provides an approximate spatial average of the pressure field, while in the near-ideal gas state this quasi-homogeneous approximation rather characterizes the back-wall pressure.}

\section{Conclusions and outlook} \label{sec:concl}

The increasing spread of technologies operating with supercritical fluids requires a deeper understanding and more accurate calculations of the appearing coupled thermomechanical phenomena. In supercritical fluids, due to the relative large thermal expansion, mechanical and thermal processes can be intensively coupled, causing both diffusive and dispersive wave propagation in the medium. To develop and test novel numerical methods supported by physics, the piston effect proves to be a well-suited benchmark problem. Although
several built-in time integration algorithms exist in commercial simulation tools, reliable numerical treatment of multiple time scale problems is a serious challenge for them, raising the need for novel self-developed numerical schemes.

During preparing the numerical procedure introduced here, we have reviewed the governing equations of linearized thermoacoustics, addressing thermodynamical questions such as entropy-preserving and entropy-increasing contributions. We have then developed a second-order accurate, fully explicit numerical time integration scheme for the linearized equations of thermoacoustics. The method is based on the splitting of the generating vector field of the spatially semi-discretized equations into reversible and irreversible -- \ie entropy-preserving and entropy-increasing -- contributions. The method is applicable for nonlinear equations, too. Within the framework of linear partial differential equations, this procedure leads to the numerical splitting of the governing equations into wave and heat equations, therefore, these separated parts can be numerically integrated according to their own individual mathematical properties.

{Our scheme -- and the associated physical background, the reversible-irreversible separation provided by GENERIC -- does not rely on linear properties, hence, a further consequent step would be to apply the method to a fuller nonlinear version of the model.} Since the presented method is fully explicit, treating material nonlinearities does not cause any particular complications, material properties should be considered as functions of the state variables and should be evaluated at the appropriate time instant, when values of the state variables are already known. 
{Based on our previous experience reported in \cite{takacs2024thermodynamically}, the omitted quadratic mechanical power terms are not expected to create problems, either. Nevertheless, phenomena like shock waves and strong dispersion caused by rapid changes in the speed of sound (\eg slightly above the critical point) may require careful attention.} 
Nonlinear treatment of thermofluid dynamical problems leads to the emergence of the convective terms, however, we do not expect any extra complications at this time, an accurate scheme developed to solve the advection equation 
{is expected to be implementable} into the sequence of splitting steps, 
{as our preliminary experience has shown (work in progress)}.

The generalization of the method to 3 spatial dimensions can be done directly along the lines of \cite{pozsar2020four}.
{
Namely, although the number of equations to be solved increases, our developed scheme does not rely on 1 spatial dimensional specialties so the only emerging task is the proper placement of the different tensorial order quantities in the discrete lattice. In detail, scalars (including the trace and scalar part of second-order tensors) should be placed at the center of a discrete cell, components of vectors at the centers of the faces, while off-diagonal components (belonging to the deviatoric part) of second-order tensors should be represented at the edges of a discrete cell.

A further opportunity to be investigated is the direct numerical exploitation of the Helmholtz decomposition (see, for example, \cite{schoder2020helmholtz,schoder2020postprocessing}). This procedure may be particularly useful for dealing with 3-dimensional problems with constant coefficients, since then the acoustic and thermal processes are decoupled from vortex evolution.}

A further advantage of the developed scheme is that the discrete problem can be formulated by finite difference equations, which methods are successfully implemented on GPUs and offer an attractive tool for fluid dynamical simulations \cite{pekkila2017methods}. Implementing finite difference methods on GPU can offer about 5 times speed-up compared to implementing finite element methods on the same GPU \cite{arduino2021gpuaccelerated}.

The approach presented here may be generalized to other dynamic problems where a coupled process takes place at significantly different time scales. In such cases, dispersion and dissipation error-free numerical solutions are required to reliably predict the time evolution of some physical quantities. One such problem is condensation induced water hammer phenomenon, which causes significant thermohydraulic transients (\ie high-amplitude pressure waves) in the primary loop of nuclear power reactors, which can lead to deformations or even to the breakage of pipelines. The study presented by Barna \etal \cite{barna2010experimental} (see also \cite{imre2010theoretical}) reports unphysical pressure peaks -- dispersion error -- obtained by numerical calculations and inapplicability of the commercial thermohydraulic tools RELAP5 and CATHARE. Figures presented in \cite{datta2016numerical} reveal similar dispersion errors. Condensation induced water hammer phenomenon is modeled via two-phase flow equations, therefore, the corresponding generalization of our presented method is required to treat this problem via our scheme. Another seminal area where the approach presented here is expected to be beneficial is sonochemistry (see, \eg \cite{alawamleh2024sonohydrogen}).

{Finally, we note that the scheme may also be applied with minor modifications to more complex wave phenomena. For instance, in superfluid Helium, two independent wave propagations are observed. One is the classical acoustic wave -- in this context, the first sound -- and the other the so-called second sound, a thermal wave originated from quantum effects \cite{landau1987fluid}. The one-fluid model of Extended Irreversible Thermodynamics for superfluid Helium describes the phenomenon via a modified Maxwell--Cattaneo--Vernotte heat flux constitutive equation \cite{jou2002second,mongiovi2018non}. In this case, to obtain a reliable dispersion error free numerical solution, the reversible part of the governing equations may be further splitted. The concrete numerical treatment, as well as clarification of an appropriate physical background (including objectivity of heat flux \cite{christov2005heat,christov2009frame}) of fluids with Maxwell--Cattaneo--Vernotte heat conduction, require a separate study. A further question is whether a phenomenon similar to the piston effect can also be observed in such fluids.}

\begin{acknowledgments}
    The authors thank Tamás Környey for the numerous useful and enlightening discussions.
    Project no.~TKP-6-6/PALY-2021 has been implemented with the support provided by the Ministry of Culture and Innovation of Hungary from the National Research, Development and Innovation Fund, financed under the TKP2021-NVA funding scheme.
    {The research was partially supported by the Sustainable Development and Technologies National Programme of the Hungarian Academy of Sciences (FFT NP FTA) and by the Hungarian Scientific Research Fund under Grant agreement NKKP-Advanced 150038 (Sz.M.).}
    The research was supported by the János Bolyai Research Scholarship of the Hungarian Academy of Sciences (K.R.).
    The work was supported by the University Researcher Scholarship Program (EKÖP) of the National Research, Development and Innovation Office, which is also supported by the Ministry of Culture and Innovation (T.D.).
    We also acknowledge KIFÜ (Governmental Agency for IT Development, Hungary) for granting us access to the Komondor HPC facility based in Hungary.
\end{acknowledgments}

\section*{\lat{In Memoriam:} Professor Tamás Környey}

Tamás Környey, our beloved friend and professor of the Department of Energy Engineering, Faculty of Mechanical Engineering, Budapest University of Technology and Economics, from whom we learned a lot about the thermodynamics of fluids, passed away on September 16, 2024, at the age of 84.

Professor Környey was born in Szeged in 1940. He obtained his degree in mechanical engineering in 1963, and after that, he worked at the Department of Energy Engineering and its predecessor departments until his retirement. He obtained his Candidate's degree in 1981 for his results on heat exchangers.

Tamás Környey joined Professor Heller (who invented the Heller--Forg\'o dry cooling tower system) in teaching Thermodynamics and Heat Transfer. He played a major role in modernizing the curriculum of these two subjects to meet modern technological requirements, as well as to convey traditional knowledge. For decades, he gave high-quality lectures on the fundamentals of thermodynamics and heat transfer, and he wrote extensive notes for both subjects, which are still used as handbooks by generations of engineers in Hungary.

Since the appearance of computer science, Professor Környey actively contributed to the field of numerical heat transport methods and calculations. He solved complex problems such as solidification in continuous casting of aluminum and heat and mass transfer during salami drying. In the field of thermodynamics, he dealt with empirical equations of the state of multiphase systems.

The authors dedicate the present work to the memory of Professor Környey. May he rest in peace.

\begin{figure}[H]
\centering
    \includegraphics[width=0.25\textwidth]{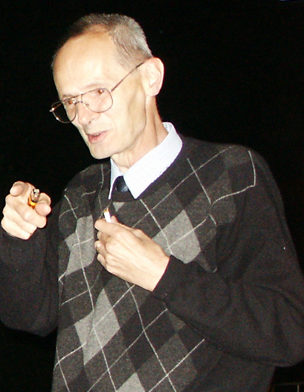}
    \\ {Tam\'as K\"ornyey (1940--2024).}
\end{figure}

\section*{Author contributions}

D.~M.~Takács: development of the numerical scheme, analytical calculations regarding the numerical aspects, numerical calculations, figures, manuscript text. T.~Fülöp: development of the numerical scheme, manuscript text. R.~Kovács: kinetic theory conclusions, manuscript text. M.~Szücs: conceptualization, theoretical background, dispersion relations, development of the numerical scheme, numerical calculations, figures, manuscript text.

\appendix
\section{The fully discretized numerical scheme}
\label{sec:app-A}

\begin{figure}[!htb]
    \centering
    \includegraphics{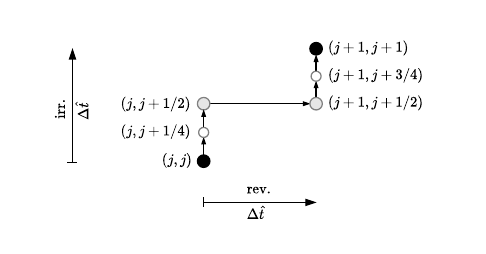}
    \caption{Illustration of substeps in time and the associated notation, for $\hv$. For other fields, the structure of the substeps is analogous, but might be staggered by $\hDt/2$, as indicated.}
    \label{fig:substeps}
\end{figure}

The reversible--irreversible vector field splitting time integration scheme \re{eq:split} yields the fully discretized numerical scheme
{
\footnotesize
\begin{align}
    \parbox{4 em}{{action of} $\displaystyle
    \psiirrev(\hDt/2)
    $}
    & \left\{
    \begin{aligned}
        \hv_{n}^{(j, j+\quart)} &= \hv_{n}^{j} - \frac{\hDt}{4} \frac{1}{\hrho^0} \pDhx{\hPi_{n+\half}^{j} - \hPi_{n-\half}^{j}} \\
        \hT_{n+\half}^{(j-\half,j+\quart)}
        & = \hT_{n+\half}^{j} - \frac{\hDt}{4} \frac{1}{\hrho^0} \pDhx{\hq_{n+1}^{j} - \hq_{n}^{j}} 
        \\
        \hq_{n}^{j+\quart} &= - \bgam \hrho^0 \frac{1}{\qPr \qRea} \pDhx{\hT_{n+\half}^{(j-\half,j+\quart)} - \hT_{n-\half}^{(j-\half,j+\quart)}} \\
        \hPi_{n+\half}^{j+\quart} &= - \hrho^0 \frac{1}{\qRea} \left( \bReta + \frac{4}{3} \right) \pDhx{\hv_{n+1}^{(j, j+\quart)} - \hv_{n}^{(j, j+\quart)}} \\
        \hv_{n}^{(j, j+\half)} &= \hv_{n}^{j} - \frac{\hDt}{2} \frac{1}{\hrho^0} \pDhx{\hPi_{n+\half}^{j+\quart} - \hPi_{n-\half}^{j+\quart}} \\
        \hT_{n+\half}^{(j-\half,j+\half)} &= \hT_{n+\half}^{j} - \frac{\hDt}{2} \frac{1}{\hrho^0} \pDhx{\hq_{n+1}^{j+\quart} - \hq_{n}^{j+\quart}}
    \end{aligned}
    \right. 
    \\[2ex]
    \parbox{4 em}{{action of} $\displaystyle
    \phirev(\hDt
    \!\:
    )
    $}
    & \left\{
    \begin{aligned}
        \hrho_{n+\half}^{j+\half} &= \hrho_{n+\half}^{j-\half} - \hDt \hrho^0 \pDhx{\hv_{n+1}^{(j, j+\half)} - \hv_{n}^{(j, j+\half)}} \\
        \hT_{n+\half}^{(j+\half, j+\half)} &= \hT_{n+\half}^{(j-\half, j+\half)} - \hDt \qB \qEca \pDhx{\hv_{n+1}^{(j, j+\half)} - \hv_{n}^{(j, j+\half)}} \\
        \hv_{n}^{(j+1, j+\half)} &= \hv_n^{(j, j+\half)}
        -
        \FTminus  
        \hDt \frac{\qB}{\bgam \hT^0} \pDhx{\hT_{n+\half}^{(j+\half, j+\half)}
        \FTminus  
        - \hT_{n-\half}^{(j+\half, j+\half)}} - \hDt \frac{1}{\bgam \hrho^0} \pDhx{\hrho_{n+\half}^{j+\half} - \hrho_{n-\half}^{j+\half}}  \\
        \hq_{n}^{j+\half} &= - \bgam \hrho^0 \frac{1}{\qPr \qRea} \pDhx{\hT_{n+\half}^{(j+\half,j+\half)} - \hT_{n-\half}^{(j+\half,j+\half)}} \\
        \hPi_{n+\half}^{j+\half} &= - \hrho^0 \frac{1}{\qRea} \left( \bReta + \frac{4}{3} \right) \pDhx{\hv_{n+1}^{(j+1, j+\half)} - \hv_{n}^{(j+1, j+\half)}} \\
    \end{aligned}
    \right. 
    \\[2ex]
    \parbox{4 em}{{action of} $\displaystyle
    \psiirrev(\hDt/2)
    $}
    & \left\{
    \begin{aligned}
        \hv_{n}^{(j+1, j+\tquart)} &= \hv_{n}^{(j+1, j+\half)} - \frac{\hDt}{4} \frac{1}{\hrho^0} \pDhx{\hPi_{n+\half}^{j+\half} - \hPi_{n-\half}^{j+\half}} \\
        \hT_{n+\half}^{(j+\half, j+\tquart)} &= T_{n+\half}^{(j+\half, j+\half)} - \frac{\hDt}{4} \frac{1}{\hrho^0} \pDhx{\hq_{n+1}^{j+\half} - \hq_{n}^{j+\half}} \\
        \hq_{n}^{j+\tquart} &= - \bgam \hrho^0 \frac{1}{\qPr \qRea} \pDhx{\hT_{n+\half}^{(j+\half,j+\tquart)} - \hT_{n-\half}^{(j+\half,j+\tquart)}} \\
        \hPi_{n+\half}^{j+\tquart} &= - \hrho^0 \frac{1}{\qRea} \left( \bReta + \frac{4}{3} \right) \pDhx{\hv_{n+1}^{(j+1, j+\tquart)} - \hv_{n}^{(j+1, j+\tquart)}} \\
        \hv_{n}^{j+1} &= \hv_{n}^{(j+1, j+\half)} - \frac{\hDt}{2} \frac{1}{\hrho^0} \pDhx{\hPi_{n+\half}^{j+\tquart} - \hPi_{n-\half}^{j+\tquart}} \\
        \hT_{n+\half}^{j+1} &= \hT_{n+\half}^{(j+\half, j+\half)} - \frac{\hDt}{2} \frac{1}{\hrho^0} \pDhx{\hq_{n+1}^{j+\tquart} - \hq_{n}^{j+\tquart}} \\
        \hq_{n}^{j+1} &= - \bgam \hrho^0 \frac{1}{\qPr \qRea} \pDhx{\hT_{n+\half}^{j+1} - \hT_{n-\half}^{j+1}} \\
        \hPi_{n+\half}^{j+1} &= - \hrho^0 \frac{1}{\qRea} \left( \bReta + \frac{4}{3} \right) \pDhx{\hv_{n+1}^{j+1} - \hv_{n}^{j+1}}
    \end{aligned}
    \right. 
\end{align}
}
where the superscript $j$ denotes the $j$th time step, that is, the numerical solution of a certain field at instant $\hht=j \hDt$. The values $\hrho_{n+\half}^{j+\half}$, $\hv_{n}^{j+1}$, $\hT_{n-\half}^{j+1}$, $\hq_{n}^{j+1}$ and $\hPi_{n+\half}^{j+1}$ are the solution values, while the others are intermediate quantities. This means that this scheme is staggered in both space and time. The double superscript indices denote that for some quantities where an irreversible and reversible part of the generating vector field also exists, the time-stepping of the two different parts is performed in steps or substeps consistent with those superscripts: the same idea appears in \cite{shang2020structurepreserving} but with a different, less straightforward notation. For a more concise presentation, the shorthand $\hT^j \equiv \hT^{(j-\half, j)}$ has been used in the case of the temperature field. For a visual illustration of the substeps, see Fig.~\ref{fig:substeps}.

Also note that for $\hrho$ and $\hT$, the initial condition has to be extended to $t=-\hDt/2$ using the reversible part of the vector field: this is straightforward in the case of homogeneous or equilibrium initial condition; for a general initial condition, a second-order accurate half-step has to be made backwards in time (\eg using the explicit midpoint method) to generate the appropriate initial conditions.

\section{
Boukari's thermal equilibration model with homogeneous pressure}
\label{sec:app-B}

{
Boukari \etal provided a thermodynamic-based explanation for the rapid temperature equilibration near the liquid-vapor critical point \cite{boukari1990critical}. In their work, temperature response to thermal excitation of a strongly compressible pure fluid closed in a constant volume is investigated. Although flow processes are neglected, the intensive pressure changes occurring in the beginning of the process are taken into account. In this case, it is more convenient to express the balance of internal energy in terms of the variables temperature and pressure. Via expressing relative volume change as
\begin{align}
    \label{eq:Dv_rho-Tp}
    \frac{\Dv \qrho}{\qrho} = - \bp \Dv T + \kT \Dv p ,
\end{align}
the balance \re{eq:heat-eq-Trho} of internal energy
can be reformulated as 
\begin{align}
    \label{eq:heat-eq-Tp}
    \qrho \left( \cv + \frac{T}{\qrho} \frac{\bp^2}{\kT} \right) \Dv T &= - \nab \cdot \qqq + T \bp \Dv p - \qqPi : \Sym{\left( \qqv \otimes \nab \right)} ,
\end{align}
where the isobaric specific heat capacity appears in the l.h.s. of \re{eq:heat-eq-Tp} [\cf \re{eq:cp}]. The model neglects gravity and assumes zero flow velocity, therefore, material time derivative coincides with the partial time derivative, the Navier--Stokes equation reduces to
\begin{align}
    \label{eq:Bou-NS}
    \nab p = 0 ,
\end{align}
and viscous dissipation disappears from \re{eq:heat-eq-Tp}. The reduced Navier--Stokes equation \re{eq:Bou-NS} expresses that the pressure field is homogeneous during the entire process, however, it may depend on time. Assuming constant material properties and 1 spatial dimensional propagation along the axial direction of a pipe with length $ X $, \re{eq:heat-eq-Tp} reduces to
\begin{align}
    \label{eq:Bou-1}
    \brho \bcp \pdt{T} &= - \pdx{\qq} + \bT \bbp \dt{p} .
\end{align}
Consequently, time evolution of pressure is also required to solve the problem. Since the investigated pipe segment is closed, the mass of the filled fluid is constant, therefore, in the linear approximation,
\begin{align}
    0 = \dt{m} = \dt{} \int_{\mathcal{V}} \qrho \dd V = \int_0^X \pdt{\qrho} A \dd x \stackrel{\re{eq:Dv_rho-Tp}}{=} - A \qrho^0 \bbp \int_0^X \pdt{T} \dd x + A X \qrho^0 \bkT \dt{p} ,
\end{align}
from which -- via \re{eq:as} -- the time evolution of the pressure reads
\begin{align}
    \label{eq:Bou-2}
    \dt{p} = \frac{1}{X} \brho \bbp \frac{\left( \bas \right)^2}{\bgam} \int\limits_0^X \pdt{T} \dd x .
\end{align}
}

{
Via the units defined in Sec.~\ref{sec:num-scheme} the non-dimensional forms of \re{eq:Bou-1} and \re{eq:Bou-2} are
\begin{align}
    \label{Bou-3}
    \pdht{\hT} &= - \frac{1}{\bgam \hrho^0} \pdhx{\hq} + \frac{\qB \qEca}{\hrho^0} \dht{\hat{p}} , \\
    \label{Bou-4}
    \dht{\hat{p}} &= \frac{\hrho^0}{\bgam} \frac{\qB}{\hT^0} \int\limits_0^1 \pdht{\hT} \dd \hx .
\end{align}
These, together with the non-dimensional Fourier's law \re{eq:dimless4}, are solved at the end of Sec.~\ref{subsec:piston}, for the same initial and boundary conditions as applied throughout Sec.~\ref{subsec:piston}.
We have solved the system of equations iteratively via the explicit Euler method with a sufficiently small time step to ensure stability.}

\bibliography{bibs_merged}

\begin{thebibliography}{69}%
\makeatletter
\providecommand \@ifxundefined [1]{%
 \@ifx{#1\undefined}
}%
\providecommand \@ifnum [1]{%
 \ifnum #1\expandafter \@firstoftwo
 \else \expandafter \@secondoftwo
 \fi
}%
\providecommand \@ifx [1]{%
 \ifx #1\expandafter \@firstoftwo
 \else \expandafter \@secondoftwo
 \fi
}%
\providecommand \natexlab [1]{#1}%
\providecommand \enquote  [1]{``#1''}%
\providecommand \bibnamefont  [1]{#1}%
\providecommand \bibfnamefont [1]{#1}%
\providecommand \citenamefont [1]{#1}%
\providecommand \href@noop [0]{\@secondoftwo}%
\providecommand \href [0]{\begingroup \@sanitize@url \@href}%
\providecommand \@href[1]{\@@startlink{#1}\@@href}%
\providecommand \@@href[1]{\endgroup#1\@@endlink}%
\providecommand \@sanitize@url [0]{\catcode `\\12\catcode `\$12\catcode `\&12\catcode `\#12\catcode `\^12\catcode `\_12\catcode `\%12\relax}%
\providecommand \@@startlink[1]{}%
\providecommand \@@endlink[0]{}%
\providecommand \url  [0]{\begingroup\@sanitize@url \@url }%
\providecommand \@url [1]{\endgroup\@href {#1}{\urlprefix }}%
\providecommand \urlprefix  [0]{URL }%
\providecommand \Eprint [0]{\href }%
\providecommand \doibase [0]{http://dx.doi.org/}%
\providecommand \selectlanguage [0]{\@gobble}%
\providecommand \bibinfo  [0]{\@secondoftwo}%
\providecommand \bibfield  [0]{\@secondoftwo}%
\providecommand \translation [1]{[#1]}%
\providecommand \BibitemOpen [0]{}%
\providecommand \bibitemStop [0]{}%
\providecommand \bibitemNoStop [0]{.\EOS\space}%
\providecommand \EOS [0]{\spacefactor3000\relax}%
\providecommand \BibitemShut  [1]{\csname bibitem#1\endcsname}%
\let\auto@bib@innerbib\@empty
\bibitem [{\citenamefont {Carl{\`{e}}s}(2010)}]{carles2010brief}%
  \BibitemOpen
  \bibfield  {author} {\bibinfo {author} {\bibfnamefont {P.}~\bibnamefont {Carl{\`{e}}s}},\ }\href {\doibase 10.1016/j.supflu.2010.02.017} {\bibfield  {journal} {\bibinfo  {journal} {The Journal of Supercritical Fluids}\ }\textbf {\bibinfo {volume} {53}},\ \bibinfo {pages} {2} (\bibinfo {year} {2010})}\BibitemShut {NoStop}%
\bibitem [{\citenamefont {Imre}\ \emph {et~al.}(2019)\citenamefont {Imre}, \citenamefont {Groniewsky}, \citenamefont {Gy{\"{o}}rke}, \citenamefont {Katona},\ and\ \citenamefont {Velmovszki}}]{imre2019anomalous}%
  \BibitemOpen
  \bibfield  {author} {\bibinfo {author} {\bibfnamefont {A.~R.}\ \bibnamefont {Imre}}, \bibinfo {author} {\bibfnamefont {A.}~\bibnamefont {Groniewsky}}, \bibinfo {author} {\bibfnamefont {G.}~\bibnamefont {Gy{\"{o}}rke}}, \bibinfo {author} {\bibfnamefont {A.}~\bibnamefont {Katona}}, \ and\ \bibinfo {author} {\bibfnamefont {D.}~\bibnamefont {Velmovszki}},\ }\href {\doibase 10.3311/ppch.12905} {\bibfield  {journal} {\bibinfo  {journal} {Periodica Polytechnica Chemical Engineering}\ }\textbf {\bibinfo {volume} {63}},\ \bibinfo {pages} {276} (\bibinfo {year} {2019})}\BibitemShut {NoStop}%
\bibitem [{\citenamefont {Carl{\`{e}}s}(1998)}]{carles1998effect}%
  \BibitemOpen
  \bibfield  {author} {\bibinfo {author} {\bibfnamefont {P.}~\bibnamefont {Carl{\`{e}}s}},\ }\href@noop {} {\bibfield  {journal} {\bibinfo  {journal} {Physics of Fluids}\ }\textbf {\bibinfo {volume} {10}},\ \bibinfo {pages} {2164} (\bibinfo {year} {1998})}\BibitemShut {NoStop}%
\bibitem [{\citenamefont {Onuki}(1997)}]{onuki1997dynamic}%
  \BibitemOpen
  \bibfield  {author} {\bibinfo {author} {\bibfnamefont {A.}~\bibnamefont {Onuki}},\ }\href {\doibase 10.1103/PhysRevE.55.403} {\bibfield  {journal} {\bibinfo  {journal} {Physical Review E}\ }\textbf {\bibinfo {volume} {55}},\ \bibinfo {pages} {403} (\bibinfo {year} {1997})}\BibitemShut {NoStop}%
\bibitem [{\citenamefont {Hasan}\ and\ \citenamefont {Farouk}(2012)}]{hasan2012thermoacoustic}%
  \BibitemOpen
  \bibfield  {author} {\bibinfo {author} {\bibfnamefont {N.}~\bibnamefont {Hasan}}\ and\ \bibinfo {author} {\bibfnamefont {B.}~\bibnamefont {Farouk}},\ }\href {\doibase 10.1016/j.supflu.2012.04.007} {\bibfield  {journal} {\bibinfo  {journal} {The Journal of Supercritical Fluids}\ }\textbf {\bibinfo {volume} {68}},\ \bibinfo {pages} {13} (\bibinfo {year} {2012})}\BibitemShut {NoStop}%
\bibitem [{\citenamefont {Lemmon}\ \emph {et~al.}(2024)\citenamefont {Lemmon}, \citenamefont {Bell}, \citenamefont {Huber},\ and\ \citenamefont {McLinden}}]{lemmon2024thermophysical}%
  \BibitemOpen
  \bibfield  {author} {\bibinfo {author} {\bibfnamefont {E.~W.}\ \bibnamefont {Lemmon}}, \bibinfo {author} {\bibfnamefont {I.~H.}\ \bibnamefont {Bell}}, \bibinfo {author} {\bibfnamefont {M.~L.}\ \bibnamefont {Huber}}, \ and\ \bibinfo {author} {\bibfnamefont {M.~O.}\ \bibnamefont {McLinden}},\ }in\ \href {\doibase 10.18434/T4D303} {\emph {\bibinfo {booktitle} {{NIST} {C}hemistry {W}eb{B}ook, {NIST} {S}tandard {R}eference {D}atabase {N}umber 69}}},\ \bibinfo {editor} {edited by\ \bibinfo {editor} {\bibfnamefont {P.~J.}\ \bibnamefont {Linstrom}}\ and\ \bibinfo {editor} {\bibfnamefont {W.~G.}\ \bibnamefont {Mallard}}}\ (\bibinfo  {publisher} {National Institute of Standards and Technology},\ \bibinfo {address} {Gaithersburg MD, 20899, USA},\ \bibinfo {year} {2024})\ \bibinfo {note} {https://webbook.nist.gov/chemistry/fluid/ (retrieved November 24, 2024)}\BibitemShut {NoStop}%
\bibitem [{\citenamefont {Knez}\ \emph {et~al.}(2014)\citenamefont {Knez}, \citenamefont {Marko{\v{c}}i{\v{c}}}, \citenamefont {Leitgeb}, \citenamefont {Primo{\v{z}}i{\v{c}}}, \citenamefont {Knez~Hrn{\v{c}}i{\v{c}}},\ and\ \citenamefont {{\v{S}}kerget}}]{knez2014industrial}%
  \BibitemOpen
  \bibfield  {author} {\bibinfo {author} {\bibfnamefont {{\v{Z}}.}~\bibnamefont {Knez}}, \bibinfo {author} {\bibfnamefont {E.}~\bibnamefont {Marko{\v{c}}i{\v{c}}}}, \bibinfo {author} {\bibfnamefont {M.}~\bibnamefont {Leitgeb}}, \bibinfo {author} {\bibfnamefont {M.}~\bibnamefont {Primo{\v{z}}i{\v{c}}}}, \bibinfo {author} {\bibfnamefont {M.}~\bibnamefont {Knez~Hrn{\v{c}}i{\v{c}}}}, \ and\ \bibinfo {author} {\bibfnamefont {M.}~\bibnamefont {{\v{S}}kerget}},\ }\href {\doibase 10.1016/j.energy.2014.07.044} {\bibfield  {journal} {\bibinfo  {journal} {Energy}\ }\textbf {\bibinfo {volume} {77}},\ \bibinfo {pages} {235} (\bibinfo {year} {2014})}\BibitemShut {NoStop}%
\bibitem [{\citenamefont {Reinsch}\ \emph {et~al.}(2017)\citenamefont {Reinsch}, \citenamefont {Dobson}, \citenamefont {Asanuma}, \citenamefont {Huenges}, \citenamefont {Poletto},\ and\ \citenamefont {Sanjuan}}]{reinsch2017utilizing}%
  \BibitemOpen
  \bibfield  {author} {\bibinfo {author} {\bibfnamefont {T.}~\bibnamefont {Reinsch}}, \bibinfo {author} {\bibfnamefont {P.}~\bibnamefont {Dobson}}, \bibinfo {author} {\bibfnamefont {H.}~\bibnamefont {Asanuma}}, \bibinfo {author} {\bibfnamefont {E.}~\bibnamefont {Huenges}}, \bibinfo {author} {\bibfnamefont {F.}~\bibnamefont {Poletto}}, \ and\ \bibinfo {author} {\bibfnamefont {B.}~\bibnamefont {Sanjuan}},\ }\href@noop {} {\bibfield  {journal} {\bibinfo  {journal} {Geothermal Energy}\ }\textbf {\bibinfo {volume} {5}},\ \bibinfo {pages} {1} (\bibinfo {year} {2017})}\BibitemShut {NoStop}%
\bibitem [{\citenamefont {Dobson}\ \emph {et~al.}(2017)\citenamefont {Dobson}, \citenamefont {Asanuma}, \citenamefont {Huenges}, \citenamefont {Poletto}, \citenamefont {Reinsch},\ and\ \citenamefont {Sanjuan}}]{dobson2017supercritical}%
  \BibitemOpen
  \bibfield  {author} {\bibinfo {author} {\bibfnamefont {P.}~\bibnamefont {Dobson}}, \bibinfo {author} {\bibfnamefont {H.}~\bibnamefont {Asanuma}}, \bibinfo {author} {\bibfnamefont {E.}~\bibnamefont {Huenges}}, \bibinfo {author} {\bibfnamefont {F.}~\bibnamefont {Poletto}}, \bibinfo {author} {\bibfnamefont {T.}~\bibnamefont {Reinsch}}, \ and\ \bibinfo {author} {\bibfnamefont {B.}~\bibnamefont {Sanjuan}},\ }in\ \href {https://hal.science/hal-01497951} {\emph {\bibinfo {booktitle} {{42nd Workshop on Geothermal Reservoir Engineering}}}},\ \bibinfo {series and number} {Proceedings '' 42nd Workshop on Geothermal Reservoir Engineering''}\ (\bibinfo {address} {Stanford, CA, United States},\ \bibinfo {year} {2017})\BibitemShut {NoStop}%
\bibitem [{\citenamefont {Rahman}\ \emph {et~al.}(2020)\citenamefont {Rahman}, \citenamefont {Dongxu}, \citenamefont {Jahan}, \citenamefont {Salvatores},\ and\ \citenamefont {Zhao}}]{rahman2020design}%
  \BibitemOpen
  \bibfield  {author} {\bibinfo {author} {\bibfnamefont {M.~M.}\ \bibnamefont {Rahman}}, \bibinfo {author} {\bibfnamefont {J.}~\bibnamefont {Dongxu}}, \bibinfo {author} {\bibfnamefont {N.}~\bibnamefont {Jahan}}, \bibinfo {author} {\bibfnamefont {M.}~\bibnamefont {Salvatores}}, \ and\ \bibinfo {author} {\bibfnamefont {J.}~\bibnamefont {Zhao}},\ }\href {\doibase 10.1016/j.pnucene.2020.103320} {\bibfield  {journal} {\bibinfo  {journal} {Progress in Nuclear Energy}\ }\textbf {\bibinfo {volume} {124}},\ \bibinfo {pages} {103320} (\bibinfo {year} {2020})}\BibitemShut {NoStop}%
\bibitem [{\citenamefont {Wu}\ \emph {et~al.}(2022)\citenamefont {Wu}, \citenamefont {Ren}, \citenamefont {Feng}, \citenamefont {Shan}, \citenamefont {Huang},\ and\ \citenamefont {Yang}}]{wu2022review}%
  \BibitemOpen
  \bibfield  {author} {\bibinfo {author} {\bibfnamefont {P.}~\bibnamefont {Wu}}, \bibinfo {author} {\bibfnamefont {Y.}~\bibnamefont {Ren}}, \bibinfo {author} {\bibfnamefont {M.}~\bibnamefont {Feng}}, \bibinfo {author} {\bibfnamefont {J.}~\bibnamefont {Shan}}, \bibinfo {author} {\bibfnamefont {Y.}~\bibnamefont {Huang}}, \ and\ \bibinfo {author} {\bibfnamefont {W.}~\bibnamefont {Yang}},\ }\href {\doibase 10.1016/j.pnucene.2022.104409} {\bibfield  {journal} {\bibinfo  {journal} {Progress in Nuclear Energy}\ }\textbf {\bibinfo {volume} {153}},\ \bibinfo {pages} {104409} (\bibinfo {year} {2022})}\BibitemShut {NoStop}%
\bibitem [{\citenamefont {Theofanous}\ \emph {et~al.}(2002{\natexlab{a}})\citenamefont {Theofanous}, \citenamefont {Tu}, \citenamefont {Dinh},\ and\ \citenamefont {Dinh}}]{theofanous2002boiling1}%
  \BibitemOpen
  \bibfield  {author} {\bibinfo {author} {\bibfnamefont {T.~G.}\ \bibnamefont {Theofanous}}, \bibinfo {author} {\bibfnamefont {J.~P.}\ \bibnamefont {Tu}}, \bibinfo {author} {\bibfnamefont {A.~T.}\ \bibnamefont {Dinh}}, \ and\ \bibinfo {author} {\bibfnamefont {T.-N.}\ \bibnamefont {Dinh}},\ }\href@noop {} {\bibfield  {journal} {\bibinfo  {journal} {Experimental Thermal and Fluid Science}\ }\textbf {\bibinfo {volume} {26}},\ \bibinfo {pages} {775} (\bibinfo {year} {2002}{\natexlab{a}})}\BibitemShut {NoStop}%
\bibitem [{\citenamefont {Theofanous}\ \emph {et~al.}(2002{\natexlab{b}})\citenamefont {Theofanous}, \citenamefont {Dinh}, \citenamefont {Tu},\ and\ \citenamefont {Dinh}}]{theofanous2002boiling2}%
  \BibitemOpen
  \bibfield  {author} {\bibinfo {author} {\bibfnamefont {T.~G.}\ \bibnamefont {Theofanous}}, \bibinfo {author} {\bibfnamefont {T.-N.}\ \bibnamefont {Dinh}}, \bibinfo {author} {\bibfnamefont {J.~P.}\ \bibnamefont {Tu}}, \ and\ \bibinfo {author} {\bibfnamefont {A.~T.}\ \bibnamefont {Dinh}},\ }\href@noop {} {\bibfield  {journal} {\bibinfo  {journal} {Experimental Thermal and Fluid Science}\ }\textbf {\bibinfo {volume} {26}},\ \bibinfo {pages} {793} (\bibinfo {year} {2002}{\natexlab{b}})}\BibitemShut {NoStop}%
\bibitem [{\citenamefont {Longmire}\ and\ \citenamefont {Banuti}(2022)}]{longmire2022onset}%
  \BibitemOpen
  \bibfield  {author} {\bibinfo {author} {\bibfnamefont {N.}~\bibnamefont {Longmire}}\ and\ \bibinfo {author} {\bibfnamefont {D.~T.}\ \bibnamefont {Banuti}},\ }\href@noop {} {\bibfield  {journal} {\bibinfo  {journal} {International Journal of Heat and Mass Transfer}\ }\textbf {\bibinfo {volume} {193}},\ \bibinfo {pages} {122957} (\bibinfo {year} {2022})}\BibitemShut {NoStop}%
\bibitem [{Note1()}]{Note1}%
  \BibitemOpen
  \bibinfo {note} {Which is correlated to high compressibility, via the Imre ellipse \cite {takacs2024leading}}\BibitemShut {NoStop}%
\bibitem [{\citenamefont {Carl{\`{e}}s}(2006)}]{carles2006thermoacoustic}%
  \BibitemOpen
  \bibfield  {author} {\bibinfo {author} {\bibfnamefont {P.}~\bibnamefont {Carl{\`{e}}s}},\ }\href@noop {} {\bibfield  {journal} {\bibinfo  {journal} {Physics of Fluids}\ }\textbf {\bibinfo {volume} {18}} (\bibinfo {year} {2006})}\BibitemShut {NoStop}%
\bibitem [{\citenamefont {Zappoli}(2003)}]{zappoli2003nearcritical}%
  \BibitemOpen
  \bibfield  {author} {\bibinfo {author} {\bibfnamefont {B.}~\bibnamefont {Zappoli}},\ }\href {\doibase 10.1016/j.crme.2003.05.001} {\bibfield  {journal} {\bibinfo  {journal} {Comptes Rendus M{\'{e}}canique}\ }\textbf {\bibinfo {volume} {331}},\ \bibinfo {pages} {713} (\bibinfo {year} {2003})}\BibitemShut {NoStop}%
\bibitem [{\citenamefont {Zappoli}\ \emph {et~al.}(2015)\citenamefont {Zappoli}, \citenamefont {Beysens},\ and\ \citenamefont {Garrabos}}]{zappoli2015heat}%
  \BibitemOpen
  \bibfield  {author} {\bibinfo {author} {\bibfnamefont {B.}~\bibnamefont {Zappoli}}, \bibinfo {author} {\bibfnamefont {D.}~\bibnamefont {Beysens}}, \ and\ \bibinfo {author} {\bibfnamefont {Y.}~\bibnamefont {Garrabos}},\ }\href {\doibase 10.1007/978-94-017-9187-8} {\emph {\bibinfo {title} {Heat transfers and related effects in supercritical fluids}}}\ (\bibinfo  {publisher} {Springer},\ \bibinfo {year} {2015})\BibitemShut {NoStop}%
\bibitem [{\citenamefont {Straub}\ \emph {et~al.}(1995)\citenamefont {Straub}, \citenamefont {Eicher},\ and\ \citenamefont {Haupt}}]{straub1995dynamic}%
  \BibitemOpen
  \bibfield  {author} {\bibinfo {author} {\bibfnamefont {J.}~\bibnamefont {Straub}}, \bibinfo {author} {\bibfnamefont {L.}~\bibnamefont {Eicher}}, \ and\ \bibinfo {author} {\bibfnamefont {A.}~\bibnamefont {Haupt}},\ }\href@noop {} {\bibfield  {journal} {\bibinfo  {journal} {Physical Review E}\ }\textbf {\bibinfo {volume} {51}},\ \bibinfo {pages} {5556} (\bibinfo {year} {1995})}\BibitemShut {NoStop}%
\bibitem [{\citenamefont {Onuki}\ \emph {et~al.}(1990)\citenamefont {Onuki}, \citenamefont {Hao},\ and\ \citenamefont {Ferrell}}]{onuki1990fast}%
  \BibitemOpen
  \bibfield  {author} {\bibinfo {author} {\bibfnamefont {A.}~\bibnamefont {Onuki}}, \bibinfo {author} {\bibfnamefont {H.}~\bibnamefont {Hao}}, \ and\ \bibinfo {author} {\bibfnamefont {R.~A.}\ \bibnamefont {Ferrell}},\ }\href@noop {} {\bibfield  {journal} {\bibinfo  {journal} {Physical Review A}\ }\textbf {\bibinfo {volume} {41}},\ \bibinfo {pages} {2256} (\bibinfo {year} {1990})}\BibitemShut {NoStop}%
\bibitem [{\citenamefont {Zappoli}\ \emph {et~al.}(1990)\citenamefont {Zappoli}, \citenamefont {Bailly}, \citenamefont {Garrabos}, \citenamefont {Le~Neindre}, \citenamefont {Guenoun},\ and\ \citenamefont {Beysens}}]{zappoli1990anomalous}%
  \BibitemOpen
  \bibfield  {author} {\bibinfo {author} {\bibfnamefont {B.}~\bibnamefont {Zappoli}}, \bibinfo {author} {\bibfnamefont {D.}~\bibnamefont {Bailly}}, \bibinfo {author} {\bibfnamefont {Y.}~\bibnamefont {Garrabos}}, \bibinfo {author} {\bibfnamefont {B.}~\bibnamefont {Le~Neindre}}, \bibinfo {author} {\bibfnamefont {P.}~\bibnamefont {Guenoun}}, \ and\ \bibinfo {author} {\bibfnamefont {D.}~\bibnamefont {Beysens}},\ }\href {\doibase 10.1103/PhysRevA.41.2264} {\bibfield  {journal} {\bibinfo  {journal} {Physical Review A}\ }\textbf {\bibinfo {volume} {41}},\ \bibinfo {pages} {2264} (\bibinfo {year} {1990})}\BibitemShut {NoStop}%
\bibitem [{\citenamefont {Garrabos}\ \emph {et~al.}(1998)\citenamefont {Garrabos}, \citenamefont {Bonetti}, \citenamefont {Beysens}, \citenamefont {Perrot}, \citenamefont {Fr{\"{o}}hlich}, \citenamefont {Carl{\`{e}}s},\ and\ \citenamefont {Zappoli}}]{garrabos1998relaxation}%
  \BibitemOpen
  \bibfield  {author} {\bibinfo {author} {\bibfnamefont {Y.}~\bibnamefont {Garrabos}}, \bibinfo {author} {\bibfnamefont {M.}~\bibnamefont {Bonetti}}, \bibinfo {author} {\bibfnamefont {D.}~\bibnamefont {Beysens}}, \bibinfo {author} {\bibfnamefont {F.}~\bibnamefont {Perrot}}, \bibinfo {author} {\bibfnamefont {T.}~\bibnamefont {Fr{\"{o}}hlich}}, \bibinfo {author} {\bibfnamefont {P.}~\bibnamefont {Carl{\`{e}}s}}, \ and\ \bibinfo {author} {\bibfnamefont {B.}~\bibnamefont {Zappoli}},\ }\href@noop {} {\bibfield  {journal} {\bibinfo  {journal} {Physical Review E}\ }\textbf {\bibinfo {volume} {57}},\ \bibinfo {pages} {5665} (\bibinfo {year} {1998})}\BibitemShut {NoStop}%
\bibitem [{\citenamefont {Boukari}\ \emph {et~al.}(1990)\citenamefont {Boukari}, \citenamefont {Shaumeyer}, \citenamefont {Briggs},\ and\ \citenamefont {Gammon}}]{boukari1990critical}%
  \BibitemOpen
  \bibfield  {author} {\bibinfo {author} {\bibfnamefont {H.}~\bibnamefont {Boukari}}, \bibinfo {author} {\bibfnamefont {J.~N.}\ \bibnamefont {Shaumeyer}}, \bibinfo {author} {\bibfnamefont {M.~E.}\ \bibnamefont {Briggs}}, \ and\ \bibinfo {author} {\bibfnamefont {R.~W.}\ \bibnamefont {Gammon}},\ }\href@noop {} {\bibfield  {journal} {\bibinfo  {journal} {Physical Review A}\ }\textbf {\bibinfo {volume} {41}},\ \bibinfo {pages} {2260} (\bibinfo {year} {1990})}\BibitemShut {NoStop}%
\bibitem [{\citenamefont {Wagner}\ \emph {et~al.}(2001)\citenamefont {Wagner}, \citenamefont {Hos},\ and\ \citenamefont {Bayazitoglu}}]{wagner2001variable}%
  \BibitemOpen
  \bibfield  {author} {\bibinfo {author} {\bibfnamefont {H.}~\bibnamefont {Wagner}}, \bibinfo {author} {\bibfnamefont {P.}~\bibnamefont {Hos}}, \ and\ \bibinfo {author} {\bibfnamefont {Y.}~\bibnamefont {Bayazitoglu}},\ }\href {\doibase 10.2514/2.6639} {\bibfield  {journal} {\bibinfo  {journal} {Journal of Thermophysics and Heat Transfer}\ }\textbf {\bibinfo {volume} {15}},\ \bibinfo {pages} {497} (\bibinfo {year} {2001})}\BibitemShut {NoStop}%
\bibitem [{\citenamefont {F{\"u}l{\"o}p}\ \emph {et~al.}(2020)\citenamefont {F{\"u}l{\"o}p}, \citenamefont {Kov{\'{a}}cs}, \citenamefont {Sz{\"u}cs},\ and\ \citenamefont {Fawaier}}]{fulop2020thermodynamical}%
  \BibitemOpen
  \bibfield  {author} {\bibinfo {author} {\bibfnamefont {T.}~\bibnamefont {F{\"u}l{\"o}p}}, \bibinfo {author} {\bibfnamefont {R.}~\bibnamefont {Kov{\'{a}}cs}}, \bibinfo {author} {\bibfnamefont {M.}~\bibnamefont {Sz{\"u}cs}}, \ and\ \bibinfo {author} {\bibfnamefont {M.}~\bibnamefont {Fawaier}},\ }\href {\doibase 10.3390/e22020155} {\bibfield  {journal} {\bibinfo  {journal} {Entropy}\ }\textbf {\bibinfo {volume} {22}},\ \bibinfo {pages} {155} (\bibinfo {year} {2020})}\BibitemShut {NoStop}%
\bibitem [{\citenamefont {Pozs{\'{a}}r}\ \emph {et~al.}(2020)\citenamefont {Pozs{\'{a}}r}, \citenamefont {Sz{\"u}cs}, \citenamefont {Kov{\'{a}}cs},\ and\ \citenamefont {F{\"u}l{\"o}p}}]{pozsar2020four}%
  \BibitemOpen
  \bibfield  {author} {\bibinfo {author} {\bibfnamefont {{\'{A}}.}~\bibnamefont {Pozs{\'{a}}r}}, \bibinfo {author} {\bibfnamefont {M.}~\bibnamefont {Sz{\"u}cs}}, \bibinfo {author} {\bibfnamefont {R.}~\bibnamefont {Kov{\'{a}}cs}}, \ and\ \bibinfo {author} {\bibfnamefont {T.}~\bibnamefont {F{\"u}l{\"o}p}},\ }\href {\doibase 10.3390/e22121376} {\bibfield  {journal} {\bibinfo  {journal} {Entropy}\ }\textbf {\bibinfo {volume} {22}},\ \bibinfo {pages} {1376} (\bibinfo {year} {2020})}\BibitemShut {NoStop}%
\bibitem [{\citenamefont {Grmela}\ and\ \citenamefont {{\"{O}}ttinger}(1997)}]{grmela1997dynamics}%
  \BibitemOpen
  \bibfield  {author} {\bibinfo {author} {\bibfnamefont {M.}~\bibnamefont {Grmela}}\ and\ \bibinfo {author} {\bibfnamefont {H.~C.}\ \bibnamefont {{\"{O}}ttinger}},\ }\href {\doibase 10.1103/physreve.56.6620} {\bibfield  {journal} {\bibinfo  {journal} {Physical Review E}\ }\textbf {\bibinfo {volume} {56}},\ \bibinfo {pages} {6620} (\bibinfo {year} {1997})}\BibitemShut {NoStop}%
\bibitem [{\citenamefont {Shang}\ and\ \citenamefont {{\"{O}}ttinger}(2020)}]{shang2020structurepreserving}%
  \BibitemOpen
  \bibfield  {author} {\bibinfo {author} {\bibfnamefont {X.}~\bibnamefont {Shang}}\ and\ \bibinfo {author} {\bibfnamefont {H.~C.}\ \bibnamefont {{\"{O}}ttinger}},\ }\href {\doibase 10.1098/rspa.2019.0446} {\bibfield  {journal} {\bibinfo  {journal} {Proceedings of the Royal Society A: Mathematical, Physical and Engineering Sciences}\ }\textbf {\bibinfo {volume} {476}},\ \bibinfo {pages} {20190446} (\bibinfo {year} {2020})}\BibitemShut {NoStop}%
\bibitem [{Note2()}]{Note2}%
  \BibitemOpen
  \bibinfo {note} {In order to avoid creating the impression that this splitting is applicable to linear equations only, instead of the frequently used term `operator splitting' we refer to this approach as \protect \emph {vector field splitting}.}\BibitemShut {Stop}%
\bibitem [{\citenamefont {Anderson}\ and\ \citenamefont {Wendt}(1995)}]{anderson1995computational}%
  \BibitemOpen
  \bibfield  {author} {\bibinfo {author} {\bibfnamefont {J.~D.}\ \bibnamefont {Anderson}}\ and\ \bibinfo {author} {\bibfnamefont {J.}~\bibnamefont {Wendt}},\ }\href@noop {} {\emph {\bibinfo {title} {Computational Fluid Dynamics}}}\ (\bibinfo  {publisher} {Springer},\ \bibinfo {year} {1995})\BibitemShut {NoStop}%
\bibitem [{\citenamefont {Versteeg}\ and\ \citenamefont {Malalasekera}(2007)}]{versteeg2007introduction}%
  \BibitemOpen
  \bibfield  {author} {\bibinfo {author} {\bibfnamefont {H.~K.}\ \bibnamefont {Versteeg}}\ and\ \bibinfo {author} {\bibfnamefont {W.}~\bibnamefont {Malalasekera}},\ }\href@noop {} {\emph {\bibinfo {title} {An Introduction to Computational Fluid Dynamics: The Finite Volume Method}}}\ (\bibinfo  {publisher} {Pearson Education Limited},\ \bibinfo {year} {2007})\BibitemShut {NoStop}%
\bibitem [{\citenamefont {Anderson}(2001)}]{anderson2001fundamentals}%
  \BibitemOpen
  \bibfield  {author} {\bibinfo {author} {\bibfnamefont {J.~D.}\ \bibnamefont {Anderson}},\ }\href@noop {} {\emph {\bibinfo {title} {Fundamentals of Aerodynamics}}}\ (\bibinfo  {publisher} {McGraw-Hill},\ \bibinfo {year} {2001})\BibitemShut {NoStop}%
\bibitem [{Note3()}]{Note3}%
  \BibitemOpen
  \bibinfo {note} {Let us make a short remark on the tensorial notations. If $ \protect \boldsymbol {\varphi } $ denotes a quantity with tensorial order $ N \ge 0 $ then its components w.r.t.\ the corresponding Cartesian basis vectors $ \protect \mathbf {e}_j , \ j = 1,2,3 $ are denoted as $ \varphi _{i_{1} \protect \dots i_{N}} , \ i_{1} , \protect \dots , i_{N} = 1,2,3 $. In accordance with the usual notation and usage of the nabla operator in continuum mechanics, \protect \textit {right} gradient, \protect \textit {right} divergence, and \protect \textit {right} curl (indicated by subscript R) are defined as \begin {align*} \protect \operatorname {grad}_{\protect \rm R} \protect \boldsymbol {\varphi } &:= \protect \boldsymbol {\varphi } \otimes \nabla = \partial _j \varphi _{i_{1} \protect \dots i_{N}} \protect \mathbf {e}_{i_1} \otimes \protect \dots \otimes \protect \mathbf {e}_{i_N} \otimes \protect \mathbf {e}_{j}, \\ \protect \operatorname {div}_{\protect \rm R} \protect \boldsymbol {\varphi } &:=
  \protect \boldsymbol {\varphi } \cdot \nabla = \partial _j \varphi _{i_{1} \protect \dots i_{N - 1} j} \protect \mathbf {e}_{i_1} \otimes \protect \dots \otimes \protect \mathbf {e}_{i_{N-1}} , \\ \protect \operatorname {curl}_{\protect \rm R} \protect \boldsymbol {\varphi } &:= \protect \boldsymbol {\varphi } \times \nabla = \epsilon _{i_{N} j k} \partial _j \varphi _{i_{1} \protect \dots i_{N}} \protect \mathbf {e}_{i_1} \otimes \protect \dots \otimes \protect \mathbf {e}_{i_{N -1}} \otimes \protect \mathbf {e}_{k} , \end {align*} where Einstein summation is applied and \begin {align*} \epsilon _{ijk} = \begin {cases} 1 , & \protect \text {if } \left ( i , j , k \right ) \protect \text { is an even permutation} , \\ - 1 , & \protect \text {if } \left ( i , j , k \right ) \protect \text { is an odd permutation} , \\ 0 , & \protect \text {otherwise} \end {cases} \end {align*} is the Levi-Civita symbol. Especially, for an arbitrary scalar field $ a $ \begin {align*} a \otimes \nabla = \nabla \otimes a =:
  \nabla a ; \end {align*} in this case, divergence and curl are not interpreted. For an arbitrary vector field $ \protect \mathbf {v} $ \begin {align*} \protect \mathbf {v} \otimes \nabla &= \left ( \nabla \otimes \protect \mathbf {v} \right )^{\protect \rm T} , & \protect \mathbf {v} \cdot \nabla &= \nabla \cdot \protect \mathbf {v} , & \protect \mathbf {v} \times \nabla &= - \nabla \times \protect \mathbf {v} . \end {align*} In what follows, wherever needed, parentheses indicate whether $\nabla $ acts to the left or to the right [$(a \nabla ) b$ vs. $a (\nabla b)$].}\BibitemShut {Stop}%
\bibitem [{\citenamefont {de~Groot}\ and\ \citenamefont {Mazur}(1962)}]{degroot1962nonequilibrium}%
  \BibitemOpen
  \bibfield  {author} {\bibinfo {author} {\bibfnamefont {S.~R.}\ \bibnamefont {de~Groot}}\ and\ \bibinfo {author} {\bibfnamefont {P.}~\bibnamefont {Mazur}},\ }\href@noop {} {\emph {\bibinfo {title} {Non-equilibrium thermodynamics}}}\ (\bibinfo  {publisher} {Dover Publications},\ \bibinfo {address} {Amsterdam},\ \bibinfo {year} {1962})\BibitemShut {NoStop}%
\bibitem [{\citenamefont {Gyarmati}(1970)}]{gyarmati1970nonequilibrium}%
  \BibitemOpen
  \bibfield  {author} {\bibinfo {author} {\bibfnamefont {I.}~\bibnamefont {Gyarmati}},\ }\href {\doibase 10.1007/978-3-642-51067-0} {\emph {\bibinfo {title} {Non-equilibrium Thermodynamics}}}\ (\bibinfo  {publisher} {Springer-Verlag},\ \bibinfo {address} {Berlin Heidelberg},\ \bibinfo {year} {1970})\BibitemShut {NoStop}%
\bibitem [{\citenamefont {Bejan}(2016)}]{bejan2016advanced}%
  \BibitemOpen
  \bibfield  {author} {\bibinfo {author} {\bibfnamefont {A.}~\bibnamefont {Bejan}},\ }\href@noop {} {\emph {\bibinfo {title} {Advanced engineering thermodynamics}}}\ (\bibinfo  {publisher} {John Wiley \& Sons},\ \bibinfo {year} {2016})\BibitemShut {NoStop}%
\bibitem [{\citenamefont {Grigull}(1964)}]{grigull1964prinzip}%
  \BibitemOpen
  \bibfield  {author} {\bibinfo {author} {\bibfnamefont {U.}~\bibnamefont {Grigull}},\ }\href {\doibase 10.1016/0017-9310(64)90020-1} {\bibfield  {journal} {\bibinfo  {journal} {International Journal of Heat and Mass Transfer}\ }\textbf {\bibinfo {volume} {7}},\ \bibinfo {pages} {23} (\bibinfo {year} {1964})}\BibitemShut {NoStop}%
\bibitem [{\citenamefont {Matolcsi}(2004)}]{matolcsi2004ordinary}%
  \BibitemOpen
  \bibfield  {author} {\bibinfo {author} {\bibfnamefont {T.}~\bibnamefont {Matolcsi}},\ }\href {http://montavid.hu/materials/Matolcsi_Ordinary_Thermodynamics_2017-04-26.pdf} {\emph {\bibinfo {title} {Ordinary thermodynamics}}}\ (\bibinfo  {publisher} {Aka\-d{\'{e}}\-mi\-ai Kiad{\'{o}} (Publishing House of the Hungarian Academy of Sciences)},\ \bibinfo {address} {Budapest},\ \bibinfo {year} {2004})\BibitemShut {NoStop}%
\bibitem [{Note4()}]{Note4}%
  \BibitemOpen
  \bibinfo {note} {Such a hierarchy is a typical outcome whenever one eliminates some of the degrees of freedom, see, \protect \textit {e.g.,\ }\cite {van2015thermodynamic,fulop2018emergence}.}\BibitemShut {Stop}%
\bibitem [{\citenamefont {Holmes}\ \emph {et~al.}(2011)\citenamefont {Holmes}, \citenamefont {Parker},\ and\ \citenamefont {Povey}}]{holmes2011temperature}%
  \BibitemOpen
  \bibfield  {author} {\bibinfo {author} {\bibfnamefont {M.~J.}\ \bibnamefont {Holmes}}, \bibinfo {author} {\bibfnamefont {N.~G.}\ \bibnamefont {Parker}}, \ and\ \bibinfo {author} {\bibfnamefont {M.~J.~W.}\ \bibnamefont {Povey}},\ }in\ \href {\doibase 10.1088/1742-6596/269/1/012011} {\emph {\bibinfo {booktitle} {Journal of Physics: Conference Series}}},\ Vol.\ \bibinfo {volume} {269}\ (\bibinfo {organization} {IOP Publishing},\ \bibinfo {year} {2011})\ p.\ \bibinfo {pages} {012011}\BibitemShut {NoStop}%
\bibitem [{\citenamefont {Wang}\ \emph {et~al.}(2019)\citenamefont {Wang}, \citenamefont {Ubachs},\ and\ \citenamefont {Van De~Water}}]{wang2019bulk}%
  \BibitemOpen
  \bibfield  {author} {\bibinfo {author} {\bibfnamefont {Y.}~\bibnamefont {Wang}}, \bibinfo {author} {\bibfnamefont {W.}~\bibnamefont {Ubachs}}, \ and\ \bibinfo {author} {\bibfnamefont {W.}~\bibnamefont {Van De~Water}},\ }\href@noop {} {\bibfield  {journal} {\bibinfo  {journal} {The Journal of Chemical Physics}\ }\textbf {\bibinfo {volume} {150}} (\bibinfo {year} {2019})}\BibitemShut {NoStop}%
\bibitem [{\citenamefont {Tak{\'{a}}cs}\ \emph {et~al.}(2024)\citenamefont {Tak{\'{a}}cs}, \citenamefont {Pozs{\'{a}}r},\ and\ \citenamefont {F{\"{u}}l{\"{o}}p}}]{takacs2024thermodynamically}%
  \BibitemOpen
  \bibfield  {author} {\bibinfo {author} {\bibfnamefont {D.~M.}\ \bibnamefont {Tak{\'{a}}cs}}, \bibinfo {author} {\bibfnamefont {{\'{A}}.}~\bibnamefont {Pozs{\'{a}}r}}, \ and\ \bibinfo {author} {\bibfnamefont {T.}~\bibnamefont {F{\"{u}}l{\"{o}}p}},\ }\href {\doibase 10.1007/s00161-024-01280-w} {\bibfield  {journal} {\bibinfo  {journal} {Continuum Mechanics and Thermodynamics}\ }\textbf {\bibinfo {volume} {36}},\ \bibinfo {pages} {525} (\bibinfo {year} {2024})}\BibitemShut {NoStop}%
\bibitem [{\citenamefont {McLachlan}\ and\ \citenamefont {Quispel}(2002)}]{mclahlan2002splitting}%
  \BibitemOpen
  \bibfield  {author} {\bibinfo {author} {\bibfnamefont {R.~I.}\ \bibnamefont {McLachlan}}\ and\ \bibinfo {author} {\bibfnamefont {G.~R.~W.}\ \bibnamefont {Quispel}},\ }\href {\doibase 10.1017/s0962492902000053} {\bibfield  {journal} {\bibinfo  {journal} {Acta Numerica}\ }\textbf {\bibinfo {volume} {11}},\ \bibinfo {pages} {341} (\bibinfo {year} {2002})}\BibitemShut {NoStop}%
\bibitem [{\citenamefont {Strang}(1968)}]{strang1968construction}%
  \BibitemOpen
  \bibfield  {author} {\bibinfo {author} {\bibfnamefont {G.}~\bibnamefont {Strang}},\ }\href@noop {} {\bibfield  {journal} {\bibinfo  {journal} {{SIAM} Journal on Numerical Analysis}\ }\textbf {\bibinfo {volume} {5}},\ \bibinfo {pages} {506} (\bibinfo {year} {1968})}\BibitemShut {NoStop}%
\bibitem [{\citenamefont {Marchuk}(1968)}]{marchuk1968some}%
  \BibitemOpen
  \bibfield  {author} {\bibinfo {author} {\bibfnamefont {G.~I.}\ \bibnamefont {Marchuk}},\ }\href@noop {} {\bibfield  {journal} {\bibinfo  {journal} {Aplikace Matematiky}\ }\textbf {\bibinfo {volume} {13}},\ \bibinfo {pages} {103} (\bibinfo {year} {1968})}\BibitemShut {NoStop}%
\bibitem [{\citenamefont {Leimkuhler}\ and\ \citenamefont {Reich}(2005)}]{leimkuhler2005simulating}%
  \BibitemOpen
  \bibfield  {author} {\bibinfo {author} {\bibfnamefont {B.}~\bibnamefont {Leimkuhler}}\ and\ \bibinfo {author} {\bibfnamefont {S.}~\bibnamefont {Reich}},\ }\href {\doibase 10.1017/cbo9780511614118} {\emph {\bibinfo {title} {{Simulating Hamiltonian Dynamics}}}}\ (\bibinfo  {publisher} {Cambridge University Press},\ \bibinfo {year} {2005})\BibitemShut {NoStop}%
\bibitem [{\citenamefont {Hairer}\ \emph {et~al.}(2006)\citenamefont {Hairer}, \citenamefont {Lubich},\ and\ \citenamefont {Wanner}}]{hairer2006geometric}%
  \BibitemOpen
  \bibfield  {author} {\bibinfo {author} {\bibfnamefont {E.}~\bibnamefont {Hairer}}, \bibinfo {author} {\bibfnamefont {C.}~\bibnamefont {Lubich}}, \ and\ \bibinfo {author} {\bibfnamefont {G.}~\bibnamefont {Wanner}},\ }\href@noop {} {\emph {\bibinfo {title} {Geometric numerical integration}}},\ \bibinfo {edition} {2nd}\ ed.,\ \bibinfo {series} {Springer Series in Computational Mathematics}, Vol.~\bibinfo {volume} {31}\ (\bibinfo  {publisher} {Springer-Verlag},\ \bibinfo {address} {Berlin},\ \bibinfo {year} {2006})\ pp.\ \bibinfo {pages} {xviii+644}\BibitemShut {NoStop}%
\bibitem [{\citenamefont {Yoshida}(1990)}]{yoshida1990construction}%
  \BibitemOpen
  \bibfield  {author} {\bibinfo {author} {\bibfnamefont {H.}~\bibnamefont {Yoshida}},\ }\href {\doibase 10.1016/0375-9601(90)90092-3} {\bibfield  {journal} {\bibinfo  {journal} {Physics Letters A}\ }\textbf {\bibinfo {volume} {150}},\ \bibinfo {pages} {262} (\bibinfo {year} {1990})}\BibitemShut {NoStop}%
\bibitem [{\citenamefont {Suzuki}(1990)}]{suzuki1990fractal}%
  \BibitemOpen
  \bibfield  {author} {\bibinfo {author} {\bibfnamefont {M.}~\bibnamefont {Suzuki}},\ }\href {\doibase https://doi.org/10.1016/0375-9601(90)90962-N} {\bibfield  {journal} {\bibinfo  {journal} {Physics Letters A}\ }\textbf {\bibinfo {volume} {146}},\ \bibinfo {pages} {319} (\bibinfo {year} {1990})}\BibitemShut {NoStop}%
\bibitem [{\citenamefont {McLachlan}(1995)}]{mclahlan1995numerical}%
  \BibitemOpen
  \bibfield  {author} {\bibinfo {author} {\bibfnamefont {R.~I.}\ \bibnamefont {McLachlan}},\ }\href {\doibase 10.1137/0916010} {\bibfield  {journal} {\bibinfo  {journal} {SIAM Journal on Scientific Computing}\ }\textbf {\bibinfo {volume} {16}},\ \bibinfo {pages} {151} (\bibinfo {year} {1995})}\BibitemShut {NoStop}%
\bibitem [{\citenamefont {Reich}(1999)}]{reich1999backward}%
  \BibitemOpen
  \bibfield  {author} {\bibinfo {author} {\bibfnamefont {S.}~\bibnamefont {Reich}},\ }\href {\doibase 10.1137/S003614299732979} {\bibfield  {journal} {\bibinfo  {journal} {SIAM Journal on Numerical Analysis}\ }\textbf {\bibinfo {volume} {36}},\ \bibinfo {pages} {1549} (\bibinfo {year} {1999})}\BibitemShut {NoStop}%
\bibitem [{Note5()}]{Note5}%
  \BibitemOpen
  \bibinfo {note} {A numerical scheme is called consistent if, as $\Delta \protect \hat {t}\to 0$, the numerical solution converges to the exact one.}\BibitemShut {Stop}%
\bibitem [{\citenamefont {Lee}(2012)}]{lee2012introduction}%
  \BibitemOpen
  \bibfield  {author} {\bibinfo {author} {\bibfnamefont {J.~M.}\ \bibnamefont {Lee}},\ }\href {\doibase 10.1007/978-1-4419-9982-5} {\emph {\bibinfo {title} {Introduction to Smooth Manifolds}}}\ (\bibinfo  {publisher} {Springer New York},\ \bibinfo {year} {2012})\BibitemShut {NoStop}%
\bibitem [{\citenamefont {Schoder}\ \emph {et~al.}(2020{\natexlab{a}})\citenamefont {Schoder}, \citenamefont {Roppert},\ and\ \citenamefont {Kaltenbacher}}]{schoder2020helmholtz}%
  \BibitemOpen
  \bibfield  {author} {\bibinfo {author} {\bibfnamefont {S.}~\bibnamefont {Schoder}}, \bibinfo {author} {\bibfnamefont {K.}~\bibnamefont {Roppert}}, \ and\ \bibinfo {author} {\bibfnamefont {M.}~\bibnamefont {Kaltenbacher}},\ }\href@noop {} {\bibfield  {journal} {\bibinfo  {journal} {SN Partial Differential Equations and Applications}\ }\textbf {\bibinfo {volume} {1}},\ \bibinfo {pages} {1} (\bibinfo {year} {2020}{\natexlab{a}})}\BibitemShut {NoStop}%
\bibitem [{\citenamefont {Schoder}\ \emph {et~al.}(2020{\natexlab{b}})\citenamefont {Schoder}, \citenamefont {Roppert},\ and\ \citenamefont {Kaltenbacher}}]{schoder2020postprocessing}%
  \BibitemOpen
  \bibfield  {author} {\bibinfo {author} {\bibfnamefont {S.}~\bibnamefont {Schoder}}, \bibinfo {author} {\bibfnamefont {K.}~\bibnamefont {Roppert}}, \ and\ \bibinfo {author} {\bibfnamefont {M.}~\bibnamefont {Kaltenbacher}},\ }\href@noop {} {\bibfield  {journal} {\bibinfo  {journal} {AIAA Journal}\ }\textbf {\bibinfo {volume} {58}},\ \bibinfo {pages} {3019} (\bibinfo {year} {2020}{\natexlab{b}})}\BibitemShut {NoStop}%
\bibitem [{\citenamefont {Pekkil{\"{a}}}\ \emph {et~al.}(2017)\citenamefont {Pekkil{\"{a}}}, \citenamefont {V{\"{a}}is{\"{a}}l{\"{a}}}, \citenamefont {K{\"{a}}pyl{\"{a}}}, \citenamefont {K{\"{a}}pyl{\"{a}}},\ and\ \citenamefont {Anjum}}]{pekkila2017methods}%
  \BibitemOpen
  \bibfield  {author} {\bibinfo {author} {\bibfnamefont {J.}~\bibnamefont {Pekkil{\"{a}}}}, \bibinfo {author} {\bibfnamefont {M.~S.}\ \bibnamefont {V{\"{a}}is{\"{a}}l{\"{a}}}}, \bibinfo {author} {\bibfnamefont {M.~J.}\ \bibnamefont {K{\"{a}}pyl{\"{a}}}}, \bibinfo {author} {\bibfnamefont {P.~J.}\ \bibnamefont {K{\"{a}}pyl{\"{a}}}}, \ and\ \bibinfo {author} {\bibfnamefont {O.}~\bibnamefont {Anjum}},\ }\href@noop {} {\bibfield  {journal} {\bibinfo  {journal} {Computer Physics Communications}\ }\textbf {\bibinfo {volume} {217}},\ \bibinfo {pages} {11} (\bibinfo {year} {2017})}\BibitemShut {NoStop}%
\bibitem [{\citenamefont {Arduino}\ \emph {et~al.}(2021)\citenamefont {Arduino}, \citenamefont {Bottauscio},\ and\ \citenamefont {Zilberti}}]{arduino2021gpuaccelerated}%
  \BibitemOpen
  \bibfield  {author} {\bibinfo {author} {\bibfnamefont {A.}~\bibnamefont {Arduino}}, \bibinfo {author} {\bibfnamefont {O.}~\bibnamefont {Bottauscio}}, \ and\ \bibinfo {author} {\bibfnamefont {L.}~\bibnamefont {Zilberti}},\ }in\ \href@noop {} {\emph {\bibinfo {booktitle} {2021 XXXIVth General Assembly and Scientific Symposium of the International Union of Radio Science (URSI GASS)}}}\ (\bibinfo {organization} {IEEE},\ \bibinfo {year} {2021})\ p.\ \bibinfo {pages} {9560241}\BibitemShut {NoStop}%
\bibitem [{\citenamefont {Barna}\ \emph {et~al.}(2010)\citenamefont {Barna}, \citenamefont {Imre}, \citenamefont {Baranyai},\ and\ \citenamefont {{\'{E}}zs{\"{o}}l}}]{barna2010experimental}%
  \BibitemOpen
  \bibfield  {author} {\bibinfo {author} {\bibfnamefont {I.~F.}\ \bibnamefont {Barna}}, \bibinfo {author} {\bibfnamefont {A.~R.}\ \bibnamefont {Imre}}, \bibinfo {author} {\bibfnamefont {G.}~\bibnamefont {Baranyai}}, \ and\ \bibinfo {author} {\bibfnamefont {G.}~\bibnamefont {{\'{E}}zs{\"{o}}l}},\ }\href@noop {} {\bibfield  {journal} {\bibinfo  {journal} {Nuclear Engineering and Design}\ }\textbf {\bibinfo {volume} {240}},\ \bibinfo {pages} {146} (\bibinfo {year} {2010})}\BibitemShut {NoStop}%
\bibitem [{\citenamefont {Imre}\ \emph {et~al.}(2010)\citenamefont {Imre}, \citenamefont {Barna}, \citenamefont {{\'{E}}zs{\"{o}}l}, \citenamefont {H{\'{a}}zi},\ and\ \citenamefont {Kraska}}]{imre2010theoretical}%
  \BibitemOpen
  \bibfield  {author} {\bibinfo {author} {\bibfnamefont {A.~R.}\ \bibnamefont {Imre}}, \bibinfo {author} {\bibfnamefont {I.~F.}\ \bibnamefont {Barna}}, \bibinfo {author} {\bibfnamefont {G.}~\bibnamefont {{\'{E}}zs{\"{o}}l}}, \bibinfo {author} {\bibfnamefont {G.}~\bibnamefont {H{\'{a}}zi}}, \ and\ \bibinfo {author} {\bibfnamefont {T.}~\bibnamefont {Kraska}},\ }\href {\doibase 10.1016/j.nucengdes.2010.03.008} {\bibfield  {journal} {\bibinfo  {journal} {Nuclear Engineering and Design}\ }\textbf {\bibinfo {volume} {240}},\ \bibinfo {pages} {1569} (\bibinfo {year} {2010})}\BibitemShut {NoStop}%
\bibitem [{\citenamefont {Datta}\ \emph {et~al.}(2016)\citenamefont {Datta}, \citenamefont {Chakravarty}, \citenamefont {Ghosh}, \citenamefont {Mukhopadhyay}, \citenamefont {Sen}, \citenamefont {Dutta},\ and\ \citenamefont {Goyal}}]{datta2016numerical}%
  \BibitemOpen
  \bibfield  {author} {\bibinfo {author} {\bibfnamefont {P.}~\bibnamefont {Datta}}, \bibinfo {author} {\bibfnamefont {A.}~\bibnamefont {Chakravarty}}, \bibinfo {author} {\bibfnamefont {K.}~\bibnamefont {Ghosh}}, \bibinfo {author} {\bibfnamefont {A.}~\bibnamefont {Mukhopadhyay}}, \bibinfo {author} {\bibfnamefont {S.}~\bibnamefont {Sen}}, \bibinfo {author} {\bibfnamefont {A.}~\bibnamefont {Dutta}}, \ and\ \bibinfo {author} {\bibfnamefont {P.}~\bibnamefont {Goyal}},\ }\href@noop {} {\bibfield  {journal} {\bibinfo  {journal} {Nuclear Engineering and Design}\ }\textbf {\bibinfo {volume} {304}},\ \bibinfo {pages} {50} (\bibinfo {year} {2016})}\BibitemShut {NoStop}%
\bibitem [{\citenamefont {Al-Awamleh}\ and\ \citenamefont {Heged{\H{u}}s}(2024)}]{alawamleh2024sonohydrogen}%
  \BibitemOpen
  \bibfield  {author} {\bibinfo {author} {\bibfnamefont {A.}~\bibnamefont {Al-Awamleh}}\ and\ \bibinfo {author} {\bibfnamefont {F.}~\bibnamefont {Heged{\H{u}}s}},\ }\href {\doibase 10.3311/PPme.37299} {\bibfield  {journal} {\bibinfo  {journal} {Periodica Polytechnica Mechanical Engineering}\ }\textbf {\bibinfo {volume} {68}},\ \bibinfo {pages} {254} (\bibinfo {year} {2024})}\BibitemShut {NoStop}%
\bibitem [{\citenamefont {Landau}\ and\ \citenamefont {Lifshitz}(1987)}]{landau1987fluid}%
  \BibitemOpen
  \bibfield  {author} {\bibinfo {author} {\bibfnamefont {L.~D.}\ \bibnamefont {Landau}}\ and\ \bibinfo {author} {\bibfnamefont {E.~M.}\ \bibnamefont {Lifshitz}},\ }\href@noop {} {\emph {\bibinfo {title} {Fluid Mechanics: Volume 6 of {C}ourse of theoretical physics}}},\ \bibinfo {edition} {2nd}\ ed.\ (\bibinfo  {publisher} {Pergamon Press},\ \bibinfo {year} {1987})\BibitemShut {NoStop}%
\bibitem [{\citenamefont {Jou}\ \emph {et~al.}(2002)\citenamefont {Jou}, \citenamefont {Lebon},\ and\ \citenamefont {Mongiov{\`\i}}}]{jou2002second}%
  \BibitemOpen
  \bibfield  {author} {\bibinfo {author} {\bibfnamefont {D.}~\bibnamefont {Jou}}, \bibinfo {author} {\bibfnamefont {G.}~\bibnamefont {Lebon}}, \ and\ \bibinfo {author} {\bibfnamefont {M.~S.}\ \bibnamefont {Mongiov{\`\i}}},\ }\href@noop {} {\bibfield  {journal} {\bibinfo  {journal} {Physical Review B}\ }\textbf {\bibinfo {volume} {66}},\ \bibinfo {pages} {224509} (\bibinfo {year} {2002})}\BibitemShut {NoStop}%
\bibitem [{\citenamefont {Mongiov{\`\i}}\ \emph {et~al.}(2018)\citenamefont {Mongiov{\`\i}}, \citenamefont {Jou},\ and\ \citenamefont {Sciacca}}]{mongiovi2018non}%
  \BibitemOpen
  \bibfield  {author} {\bibinfo {author} {\bibfnamefont {M.~S.}\ \bibnamefont {Mongiov{\`\i}}}, \bibinfo {author} {\bibfnamefont {D.}~\bibnamefont {Jou}}, \ and\ \bibinfo {author} {\bibfnamefont {M.}~\bibnamefont {Sciacca}},\ }\href@noop {} {\bibfield  {journal} {\bibinfo  {journal} {Physics Reports}\ }\textbf {\bibinfo {volume} {726}},\ \bibinfo {pages} {1} (\bibinfo {year} {2018})}\BibitemShut {NoStop}%
\bibitem [{\citenamefont {Christov}\ and\ \citenamefont {Jordan}(2005)}]{christov2005heat}%
  \BibitemOpen
  \bibfield  {author} {\bibinfo {author} {\bibfnamefont {C.}~\bibnamefont {Christov}}\ and\ \bibinfo {author} {\bibfnamefont {P.}~\bibnamefont {Jordan}},\ }\href@noop {} {\bibfield  {journal} {\bibinfo  {journal} {Physical Review Letters}\ }\textbf {\bibinfo {volume} {94}},\ \bibinfo {pages} {154301} (\bibinfo {year} {2005})}\BibitemShut {NoStop}%
\bibitem [{\citenamefont {Christov}(2009)}]{christov2009frame}%
  \BibitemOpen
  \bibfield  {author} {\bibinfo {author} {\bibfnamefont {C.}~\bibnamefont {Christov}},\ }\href@noop {} {\bibfield  {journal} {\bibinfo  {journal} {Mechanics Research Communications}\ }\textbf {\bibinfo {volume} {36}},\ \bibinfo {pages} {481} (\bibinfo {year} {2009})}\BibitemShut {NoStop}%
\bibitem [{\citenamefont {Tak{\'a}cs}\ \emph {et~al.}(2024)\citenamefont {Tak{\'a}cs}, \citenamefont {F{\"u}l{\"o}p},\ and\ \citenamefont {Imre}}]{takacs2024leading}%
  \BibitemOpen
  \bibfield  {author} {\bibinfo {author} {\bibfnamefont {D.~M.}\ \bibnamefont {Tak{\'a}cs}}, \bibinfo {author} {\bibfnamefont {T.}~\bibnamefont {F{\"u}l{\"o}p}}, \ and\ \bibinfo {author} {\bibfnamefont {A.~R.}\ \bibnamefont {Imre}},\ }\href {\doibase 10.1016/j.supflu.2024.106216} {\bibfield  {journal} {\bibinfo  {journal} {The Journal of Supercritical Fluids}\ }\textbf {\bibinfo {volume} {208}},\ \bibinfo {pages} {106216} (\bibinfo {year} {2024})}\BibitemShut {NoStop}%
\bibitem [{\citenamefont {V{\'a}n}\ \emph {et~al.}(2015)\citenamefont {V{\'a}n}, \citenamefont {Kov{\'a}cs},\ and\ \citenamefont {F{\"u}l{\"o}p}}]{van2015thermodynamic}%
  \BibitemOpen
  \bibfield  {author} {\bibinfo {author} {\bibfnamefont {P.}~\bibnamefont {V{\'a}n}}, \bibinfo {author} {\bibfnamefont {R.}~\bibnamefont {Kov{\'a}cs}}, \ and\ \bibinfo {author} {\bibfnamefont {T.}~\bibnamefont {F{\"u}l{\"o}p}},\ }\href {\doibase 10.3176/PROC.2015.3S.09} {\bibfield  {journal} {\bibinfo  {journal} {Proceedings of the Estonian Academy of Sciences}\ }\textbf {\bibinfo {volume} {64}},\ \bibinfo {pages} {389} (\bibinfo {year} {2015})}\BibitemShut {NoStop}%
\bibitem [{\citenamefont {F{\"u}l{\"o}p}\ \emph {et~al.}(2018)\citenamefont {F{\"u}l{\"o}p}, \citenamefont {Kov{\'a}cs}, \citenamefont {Lovas}, \citenamefont {Rieth}, \citenamefont {Fodor}, \citenamefont {Sz{\"u}cs}, \citenamefont {V{\'a}n},\ and\ \citenamefont {Gr{\'o}f}}]{fulop2018emergence}%
  \BibitemOpen
  \bibfield  {author} {\bibinfo {author} {\bibfnamefont {T.}~\bibnamefont {F{\"u}l{\"o}p}}, \bibinfo {author} {\bibfnamefont {R.}~\bibnamefont {Kov{\'a}cs}}, \bibinfo {author} {\bibfnamefont {{\'A}.}~\bibnamefont {Lovas}}, \bibinfo {author} {\bibfnamefont {{\'A}.}~\bibnamefont {Rieth}}, \bibinfo {author} {\bibfnamefont {T.}~\bibnamefont {Fodor}}, \bibinfo {author} {\bibfnamefont {M.}~\bibnamefont {Sz{\"u}cs}}, \bibinfo {author} {\bibfnamefont {P.}~\bibnamefont {V{\'a}n}}, \ and\ \bibinfo {author} {\bibfnamefont {G.}~\bibnamefont {Gr{\'o}f}},\ }\href {\doibase 10.3390/e20110832} {\bibfield  {journal} {\bibinfo  {journal} {Entropy}\ }\textbf {\bibinfo {volume} {20}},\ \bibinfo {pages} {832} (\bibinfo {year} {2018})}\BibitemShut {NoStop}%
\end{thebibliography}%

\end{document}